\documentclass[12pt]{amsart}
\usepackage{amssymb, enumerate}
\usepackage{xcolor}

\theoremstyle{plain}
\newtheorem{theorem}{Theorem}[section]
\newtheorem{prop}[theorem]{Proposition}
\newtheorem{proposition}[theorem]{Proposition}
\newtheorem{lemma}[theorem]{Lemma}
\newtheorem{corollary}[theorem]{Corollary}

\numberwithin{equation}{section}

\theoremstyle{definition}
\newtheorem{definition}[theorem]{Definition}
\newtheorem{remark}[theorem]{Remark}

\newcommand\D{\mathcal{D}}

\newcommand\R{\mathbf{R}}

\newcommand\dd{\,{\rm d}}

\newcommand{\Sc}{\mathcal S_c}
\newcommand{\supp}{\operatorname{supp}}
\newcommand{\norm}[1]{\lVert#1\rVert}

\def\d{\rho}

\allowdisplaybreaks
\usepackage[colorlinks, citecolor=blue,pagebackref,hypertexnames=false]{hyperref}

\hypersetup{pdftitle={Poisson functions on manifolds with ends and applications, v23}}

\begin{document}

	\title[Poisson functions on manifolds with ends]{Poisson area and non-tangential maximal functions on manifolds with ends and applications}
	\date{}

 	\author{Peng Chen}
 	\address{
 		Peng Chen, Department of
 		Mathematics, Sun Yat-sen University
 		Guangzhou, 510275, China}
 	\email{chenpeng3@mail.sysu.edu.cn}

 	\author{Ji Li}
 	
 	\address{
 		J. Li, School of
 		Mathematical and Physical Sciences, Macquarie
 		University, NSW, 2109, Australia}
 	\email{ji.li@mq.edu.au}
 	
 	\author{Lixin Yan}
 	\address{
 		Lixin Yan, Department of
 		Mathematics, Sun Yat-sen University
 		Guangzhou, 510275, China}
 	\email{mcsylx@mail.sysu.edu.cn}

	\begin{abstract}
		Let $M=\mathbb R^m\sharp(\mathbb R^n\times\mathbb S^{m-n})$,
		where $m>n\geq3$. The full-cone Poisson area and non-tangential
		maximal functions fail to be of weak type $(p,p)$ for $1<p<2$.
		Let $K$ be the compact part joining the two ends and put
		$\rho(x)=1+d(x,K)$. For fixed small $c>0$, we restrict the horizontal
		distance from the vertex to $c\rho(x)$ and retain every positive time.
		The resulting maximal function is of weak type $(1,1)$ and
		bounded on $L^p$ for $1<p\leq\infty$. The corresponding area
		function uses the full spatial and time gradient. It satisfies
		a Fefferman--Stein distribution inequality. Its $L^p$ norm is
		equivalent to that of $f$ for $1<p\leq2$.
		This area function defines a global Hardy space. We prove its $\Delta$-atomic description, its injective embedding
		into $L^1$, its characterization by the restricted maximal
		function, and its real interpolation identity with $L^2$. Direct estimates for spectral-band kernels
		then give a Mihlin multiplier theorem from this Hardy space to
		$L^1$. Interpolation and duality give $L^p$ boundedness for every
		$1<p<\infty$.
	\end{abstract}
	
	\maketitle
	
	\section{Introduction}
	
	In his 1970 book \cite{Stein1970}, Stein
	developed Littlewood--Paley theory for semigroups on general
	measure spaces. This theory includes maximal function estimates,
	$g$-function estimates and a multiplier theorem.
	
	Let $\{T_t\}_{t\geq0}$ be a semigroup on a measure space.
	Assume that it is strongly continuous on $L^2$ and that each $T_t$
	is self-adjoint on $L^2$, preserves nonnegative functions, and satisfies
	$$
	T_t1=1,\qquad
	\|T_tf\|_q\leq\|f\|_q,\qquad 1\leq q\leq\infty.
	$$
	These are the symmetric diffusion semigroups in Stein's terminology.
	In \cite[Chapters III and IV]{Stein1970},	Stein obtained the   maximal theorem, 
	$$
	\left\|\sup_{t>0}|T_tf|\right\|_p
	\leq C_p\|f\|_p,\qquad 1<p\leq\infty;
	$$	
 He also proved the	Littlewood--Paley type inequality:  
 	$$
\left\|\left (\int_0^\infty
	t\big|\partial_tT_tf\big|^2\, dt  \right)^{1/2} \right\|_p \leq C_p\|f\|_p,\qquad 1<p<\infty.
	$$
	In the same setting, let $\Delta$ be the nonnegative self-adjoint
	operator such that $T_t=e^{-t\Delta}$. 
	His multiplier theorem
applies to functions
	of Laplace transform type:
	$$
	\mathfrak m(\lambda)
	=\lambda\int_0^\infty e^{-t\lambda}h(t)\,dt,
	\qquad h\in L^\infty(0,\infty),\quad \lambda>0.
	$$
	Setting $\mathfrak m(0)=0$, one has
	$
	\|\mathfrak m(\Delta)f\|_p
	\leq C_p\|h\|_\infty\|f\|_p,\ 1<p<\infty.
	$

	These results depend on semigroup properties and require no
	metric, doubling condition or heat kernel estimate.

	In classical harmonic analysis,  the area and non-tangential maximal functions are two fundamental tools that  serve to characterize at once both non-tangential convergence and also $L^p$, $(1<p<\infty)$ and $H^p$,  $(0<p\leq1)$ norms (see \cite[Chapter VII]{St1},  \cite{Stein1993}). 
	These two functions use distance and
	volume in their definitions through the cone. However, these
	two functions, when defined using the classical cones, need not be
	bounded on $L^p$ on a manifold whose volume measure is not doubling.
	We prove this failure on a manifold with two ends and then
	restrict the cones.
	
	Let $m>n\geq3$, and put
	$
	\mathcal R^n=\mathbb R^n\times\mathbb S^{m-n},
	$
	where $\mathbb S^{m-n}$ is the unit sphere in
	$\mathbb R^{m-n+1}$. We consider a connected complete Riemannian
	manifold
	$$
	M=\mathbb R^m\sharp\mathcal R^n.
	$$
	More precisely, outside a compact set $K$ with smooth boundary,
	$M$ has two ends $E_m$ and $E_n$. They are isometric to the
	complements of compact sets in $\mathbb R^m$ and $\mathcal R^n$,
	respectively. Fix $K$ so that the two ends are the model regions
	$r>r_0$, where $r_0\geq1$ is fixed and $r$ is the Euclidean
	radial coordinate in the noncompact factor. Then
	$d(x,K)=r(x)-r_0$ on either end. Write $d$ for the geodesic
	distance, $\mu$ for Riemannian measure, and
	$
	V(y,t)=\mu(B(y,t)).
	$
	The local dimension is $m$ on both ends, whereas large annuli have
	volume of order $R^m$ on $E_m$ and $R^n$ on $E_n$.
	In particular, $\mu(M)=\infty$ and $\mu$ is not doubling on $M$.
	Grigor'yan and Saloff-Coste obtained two-sided heat kernel estimates
	on manifolds with ends \cite{GS}.
	
	Here $\Delta$ is the nonnegative self-adjoint Laplace--Beltrami
	operator on $L^2(M)$. For $t>0$, set
	$
	P_t=e^{-t\sqrt\Delta},\  u(y,t)=P_tf(y),
	\  \nabla_{y,t}u=(\nabla_yu,\partial_tu).
	$
	For $\beta>0$, define
	$$
	\mathcal N_P^\beta f(x)=\sup_{(y,t)\in\Gamma^\beta(x)}|u(y,t)|
	\qquad {\rm with}\qquad
	\Gamma^\beta(x)=\{(y,t)\in M\times(0,\infty):d(x,y)<\beta t\}.
	$$
	We write $\Gamma(x)=\Gamma^1(x)$ and define
	\begin{equation}\label{eq:area-definition}
		S(u)(x)=\left(\iint_{\Gamma(x)}
		|\nabla_{y,t}u(y,t)|^2\frac{t\,d\mu(y)\,dt}{V(y,t)}
		\right)^{1/2},\qquad x\in M.
	\end{equation}
	Both definitions use the full geodesic cone at every point of $M$.
	There is no restriction on the end containing $y$. The spectral theorem and Tonelli give
	\begin{equation}\label{eq:area-L2-intro}
		\|S(P_tf)\|_{L^2(M)}^2=\frac12\|f\|_{L^2(M)}^2.
	\end{equation}

	Our first main result is the following: both operators
	fail to be of weak type $(p,p)$ when $1<p<2$, as the following
	theorem shows.
	\begin{theorem}\label{thm:counterexamples}
		Let $1<p<2$ and $\beta>0$. There is a family of nonnegative functions
		$f_R\in C_c^\infty(M)$, $R\to\infty$, such that, with $u_R=P_tf_R$,
		\begin{equation}\label{eq:maximal-ratio}
			\frac{\|\mathcal N_P^\beta f_R\|_{L^{p,\infty}(M)}}
			{\|f_R\|_{L^p(M)}}\geq c_{p,\beta}R^{(m-n)/p}
		\end{equation}
		and
		\begin{equation}\label{eq:area-ratio}
			\frac{\|S(u_R)\|_{L^{p,\infty}(M)}}{\|f_R\|_{L^p(M)}}
			\geq c_pR^{(m-n)(1/p-1/2)}.
		\end{equation}
		Consequently, neither operator is of weak type $(p,p)$ or bounded
		on $L^p(M)$ for $1<p<2$.
	\end{theorem}

	The theorem shows that the classical definitions of area and
	non-tangential maximal functions need not give $L^p$ bounds on
	every manifold. Their definitions must be adjusted case by case, that is, the cones must be chosen according to a specific manifold.
	
	On this manifold with ends, we restrict the horizontal radius of each cone. For
	$0<c\leq1/100$, put
	$$
	\Gamma_c(x)=\{(y,t):t>0,\ d(x,y)<\min\{t,c\rho(x)\}\},
	\qquad \rho(x)=1+d(x,K).
	$$
	The associated non-tangential maximal and area functions are
	$$
	\mathcal N_cf(x)=\sup_{(y,t)\in\Gamma_c(x)}|P_tf(y)|,
	\qquad
	\Sc f(x)^2=\iint_{\Gamma_c(x)}
	|\nabla_{y,t}P_tf(y)|^2\frac{t\,d\mu(y)\,dt}{V(y,t)}.
	$$
	
	We first prove the following properties.
	\begin{theorem}\label{prop1}
		For $1<p<\infty$ and $0<c\leq1/100$, denote by $C_{p,c}$ a positive constant depending only on $p$ and $c$.
		Then we have
		
		(1) $\mathcal N_c$ is of weak type $(1,1)$ and
		$
		\|\mathcal N_cf\|_p\leq C_{p,c}\|f\|_p$ for
		$ 1<p\leq\infty.
		$
		
		(2) For every $\lambda>0$,
		$$
		\mu\{\Sc f>\lambda\}
		\leq C_c\mu\{\mathcal N_cf>\lambda\}
		+\frac{C_c}{\lambda^2}\int_0^\lambda
		s\,\mu\{\mathcal N_cf>s\}\,ds.
		$$
		
		(3) Thus, we further have
		$
		\|\Sc f\|_p\leq C_{p,c}\|f\|_p,
		$ for $1<p\leq2.
		$
		
		(4) There exists a positive constant $C$ such that for every $f\in L^2(M)$,
		\begin{equation}\label{eq:intro-L2-two-sided}
			C^{-1}\|f\|_2\leq\|\Sc f\|_2.
		\end{equation}
		For $1<p<2$ and $f\in L^p(M)$, we also have
		\begin{equation}\label{eq:intro-Lp-lower}
			\|f\|_p\leq C_{p,c}\|\Sc f\|_p.
		\end{equation}
	\end{theorem}
	To prove (1), we bound $\mathcal N_cf$ pointwise by the radial
	heat maximal function applied to $|f|$. We then apply \cite{DLS}
	to obtain the weak type $(1,1)$ and $L^p$ estimates.
	
	The estimate in (2) is a form of the Fefferman--Stein inequality
	from their paper on Euclidean Hardy spaces in
	\cite{FS}. Related estimates were later
	proved for non-symmetric elliptic equations by Hofmann, Kenig,
	Mayboroda and Pipher \cite{HKMP}, and on Shilov boundaries by
	Li \cite{L}. We use the argument in \cite{L} to prove (2). Integration gives (3).
	
	To prove (4), we use local gradient estimates and a radial
	Hardy inequality on both ends to bound the $L^2$ norm of the
	Poisson extension at a height comparable to $\rho(x)$ by
	$\|\Sc f\|_2$. Green's formula with a cutoff at this height
	then gives $\|f\|_2\leq C\|\Sc f\|_2$.
	For $1<p<2$, we test a local Green formula against Poisson
	extensions on the complete models. A tent space estimate and
	the $L^p$ radial Hardy inequality give \eqref{eq:intro-Lp-lower};
	see Theorem~\ref{thm:reverse-Lp}.

We now define the Hardy space associated with $\Delta$ and prove
its atomic decomposition and its characterizations by the area
and non-tangential maximal functions.
	
	\begin{definition}
	Define
		\begin{equation}\label{eq:H1}
			H^1_{S_c}(M)
			=\overline{\{f\in L^2(M):\Sc f\in L^1(M)\}}
			^{\,\|\Sc(\cdot)\|_1}.
		\end{equation}
		We also define the atoms: 
		$a=\Delta b$, with $b$ and $\Delta b$ supported in a ball $B(x,r)$,
		$r\leq\delta\rho(x)$, and
		$$
		\|b\|_2\leq r^2\mu(B)^{-1/2},\qquad
		\|\Delta b\|_2\leq\mu(B)^{-1/2}.
		$$
		Here $\delta>0$ is fixed sufficiently small.
	\end{definition}
	
	\begin{theorem}\label{thm:hardy-characterization}\label{thm:interpolation}
		Fix $0<c\leq1/100$ and a sufficiently small $0<\delta<1/16$.
		The following assertions hold.
		\begin{enumerate}[(i)]
			\item The completion $ H^1_{S_c}(M) $ embeds injectively into $L^1(M)$.
			Every $f\in H^1_{S_c}(M) $ has a representation
			\begin{equation}\label{eq:intro-atomic}
				f=\sum_{\nu\geq1}\lambda_\nu a_\nu,
				\qquad \sum_{\nu\geq1}|\lambda_\nu|<\infty,
				\qquad
				\|f\|_{ H^1_{S_c}(M) }\approx
				\inf_{f=\sum_\nu\lambda_\nu a_\nu}
				\sum_{\nu\geq1}|\lambda_\nu|.
			\end{equation}
			The series converges in $L^1(M)$ and in $ H^1_{S_c}(M) $.
			Conversely, every such atomic sum belongs to $ H^1_{S_c}(M) $.
			
			\item The restricted maximal function characterizes this space:
			\begin{equation}\label{eq:intro-maximal-Hardy}
				H^1_{S_c}(M) =\{f\in L^1(M):\mathcal N_cf\in L^1(M)\},
				\qquad
				\|f\|_{ H^1_{S_c}(M) }\approx\|\mathcal N_cf\|_1.
			\end{equation}
			
			\item For $1<p<2$,
			\begin{equation}\label{eq:intro-interpolation}
				( H^1_{S_c}(M) ,L^2(M))_{2-2/p,p}=L^p(M)
				\quad\text{with equivalent norms}.
			\end{equation}
			The same identity holds for the atomic space in (i).
		\end{enumerate}
		The constants may depend on $M,c,\delta$ and, in (iii), on $p$.
	\end{theorem}
	
	Section \ref{sec:interp-statement} proves the atomic decomposition,
	justifies its use on the completion, and proves the maximal
	characterization.
	The weighted Green estimate and the decomposition in the
	following sections prove the interpolation identity.
	Equation \eqref{eq:intro-L2-two-sided} identifies the $L^2$
	completion defined by $\Sc$.
	
	We apply this identity to spectral multipliers. The first step is
	a uniform $L^1$ estimate on every atom. The atomic description then
	gives a map from $ H^1_{S_c}(M) $ to $L^1$, consistent with the $L^2$
	spectral definition. Interpolation gives $1<p<2$, and duality gives
	$2<p<\infty$.
	
	For an integer $N\geq0$ and $\mathfrak m\in C^N((0,\infty))$, put
	\begin{equation}\label{eq:mihlin-4}
		A_N(\mathfrak m)=\max_{0\leq j\leq N}\sup_{\lambda>0}
		\lambda^j|\mathfrak m^{(j)}(\lambda)|.
	\end{equation}
	Here $\lambda$ is the spectral variable of $\Delta$.
	
	\begin{theorem}\label{thm:spectral-multiplier}
		Let $N$ be an integer with $N>(m+3)/2$, and let
		$\mathfrak m\in C^N((0,\infty))$ satisfy $A_N(\mathfrak m)<\infty$.
		Then $\mathfrak m(\Delta)$ has a unique bounded extension from
		$ H^1_{S_c}(M) $ to $L^1(M)$ such that
		\begin{equation}\label{eq:mihlin-42}
			\|\mathfrak m(\Delta)f\|_1
			\leq C_{N,\delta,c}A_N(\mathfrak m)\|f\|_{ H^1_{S_c}(M) }.
		\end{equation}
		This extension agrees with the $L^2$ spectral definition on
		$ H^1_{S_c}(M) \cap L^2$. Interpolation and duality give, for every
		$1<p<\infty$,
		\begin{equation}\label{eq:mihlin-43}
			\|\mathfrak m(\Delta)f\|_p
			\leq C_{p,N}A_N(\mathfrak m)\|f\|_p.
		\end{equation}
		No value or limit of $\mathfrak m$ at zero is required.
		One may take $N=\lfloor(m+3)/2\rfloor+1$.
	\end{theorem}
	
	The proof uses the real heat kernel estimates and finite
	propagation. Cutoffs on the exact product ends give a
	commutator supported in a fixed compact set. Its kernel
	estimate gives uniform $L^1$ bounds for spectral bands near
	zero and hence for the $\Delta$-atoms of
	Theorem \ref{thm:hardy-characterization}.
	
	Section \ref{sec:restricted-bounds} proves the maximal and area
	estimates and non-tangential boundary convergence.
	Section \ref{sec:rfs} proves the Fefferman--Stein inequality.
	Section \ref{sec:L2-statement} gives the two-sided $L^2$ estimate.
	Section \ref{sec:interp-statement} proves the atomic and maximal
	characterizations and the reverse $L^p$ estimate for $1<p<2$.
	The following sections prove the interpolation identity.
	Section \ref{sec:spectral-multipliers} proves
	Theorem \ref{thm:spectral-multiplier}.

	\begin{remark}\label{rem:atomic-cancellation}
		Our cones contain every positive height, and each atom has
		the form $a=\Delta b$. In particular,
		$$
		\int_M a\,\overline h\dd\mu=0
		\quad\text{for every bounded harmonic function }h.
		$$
		Indeed, choose $\zeta\in C_c^\infty(M)$ equal to $1$ near
		the atom's ball. Then
		$\langle \Delta b,\zeta h\rangle=\langle b,\Delta(\zeta h)\rangle=0$.
		The same cancellation holds for every $f\in H^1_{S_c}(M) $ by
		the $L^1$ convergence of its atomic series. The multiplier
		theorem requires no value or limit of $\mathfrak m$ at zero.
		In particular, it applies to $\mathfrak m(\lambda)=\lambda^{i\tau}$,
		$\tau\ne0$. Proposition \ref{prop:imaginary-cancellation}
		shows that the $L^2$ size condition on a unit ball, without
		cancellation, does not give a uniform $L^1$ estimate for this
		multiplier.
	\end{remark}
	
	Throughout the paper, constants may depend on $M,m,n$ and fixed
	cutoff functions. Constants for the restricted cones may also
	depend on $c$. Dependence on a fixed aperture or a fixed time
	multiple is indicated when needed. No constant in the estimates
	for $f_R$ depends on $R$. We write $A\lesssim B$ for $A\leq CB$
	and $A\approx B$ when both inequalities hold.

	\section{Geometry and estimates for the Poisson semigroup}
	
	\subsection{Annuli and volume bounds}
	Let $K$ and $r_0$ be as in the introduction. On each end
	$E_d=\{r>r_0\}$, write $r(x)$ for the Euclidean radial
	coordinate in the noncompact factor. The set $K$ is contained
	in a fixed connected compact neighbourhood. Whenever a radial cutoff is zero for
	$r\leq r_0$, we extend it by zero to the rest of $M$.
	We shall always take $R\geq R_0$, where $R_0>4r_0$ is fixed.
	Put
	$$
	X_R=\{x\in E_m:R<r(x)<2R\},\qquad
	Y_R=\{y\in E_n:R<r(y)<2R\}.
	$$
	Let $\omega_d$ be the volume of the unit ball in $\mathbb R^d$,
	and let $v_{m-n}$ be the volume of $\mathbb S^{m-n}$.
	The isometries with the model regions give the exact formulas
	$$
	\mu(X_R)=\omega_m(2^m-1)R^m,\qquad
	\mu(Y_R)=v_{m-n}\omega_n(2^n-1)R^n.
	$$
	In particular,
	\begin{equation}\label{eq:annuli-volume}
		\mu(X_R)\approx R^m,\qquad \mu(Y_R)\approx R^n.
	\end{equation}
	There is a constant $D\geq1$ such that
	\begin{equation}\label{eq:annuli-distance}
		d(x,y)\leq DR\quad(x\in X_R,\ y\in Y_R),\qquad
		d(y,z)\leq DR\quad(y,z\in Y_R).
	\end{equation}
	Indeed, radial paths from these points to a fixed bounded
	neighbourhood of $K$ have length at most a constant times $R$;
	paths inside that neighbourhood have bounded length.
	
	We also use the global upper bound
	\begin{equation}\label{eq:volume-upper}
		V(y,t)\leq C t^m,\qquad y\in M,\ t>0.
	\end{equation}
	The model coordinates and finitely many charts on the compact
	part give $V(y,t)\leq Ct^m$ for $0<t\leq t_*$, uniformly in $y$.
	The uniform bound on $V(y,1)$ extends this to $t\leq1$.
	For $t\geq1$, a ball disjoint from $K$ lies in one model end,
	where its volume is at most $Ct^m$.
	If $B(y,t)$ meets $K$, then $B(y,t)\subset\{z:d(z,K)<2t\}$, so
	$$
	V(y,t)\leq C(1+t^m+t^n)\leq C't^m.
	$$
	This proves \eqref{eq:volume-upper}.
	The exact metrics and finitely many charts on the compact part also give,
	for each fixed $R_0<\infty$,
	\begin{equation}\label{eq:bounded-volume-lower}
		V(x,R)\geq c_{R_0}R^m,\qquad x\in M,\quad0<R\leq R_0.
	\end{equation}
	Indeed, uniform normal coordinates give the bound for $R\leq r_*$.
	For $r_*<R\leq R_0$, monotonicity gives
	$V(x,R)\geq c r_*^m\geq c(r_*/R_0)^mR^m$.

	To check that $\mu$ is not doubling, choose $y_R\in E_n$ on a
	fixed radial ray with $r(y_R)=R$, and put $t_R=3R/4$.
	For large $R$, $d(y_R,K)>t_R$. Minimizing paths of length less
	than $t_R$ stay in $E_n$, so comparison with the product model gives
	$V(y_R,t_R)\leq CR^n$.
	Paths through a fixed compact neighbourhood of $K$ give
	$$
	\{x\in E_m:R/8<r(x)<R/4\}\subset B(y_R,2t_R),
	\qquad V(y_R,2t_R)\geq C^{-1}R^m.
	$$
	Hence
	$$
	\frac{V(y_R,2t_R)}{V(y_R,t_R)}\geq c'R^{m-n}\longrightarrow\infty.
	$$
	
	\subsection{The two kernel estimates used below}
	Let $h_s(y,z)$ denote the kernel of $e^{-s\Delta}$.
	The small-end heat kernel lower bound of \cite{GS}, in the form
	recalled in \cite[Theorem A, part 7]{BDLW}, gives
	\begin{equation}\label{eq:heat-lower}
		h_s(y,z)\geq C^{-1} s^{-n/2}
		\exp\left(-C\frac{d(y,z)^2}{s}\right),
		\qquad y,z\in E_n,\quad s>1.
	\end{equation}
	Only this term of that lower bound is needed.
	By \cite[Theorem 2.2, part 6, with $k=0$]{BDLW},
	\begin{equation}\label{eq:poisson-upper}
		0\leq P_t(y,z)\leq C\left[
		\frac{t}{(t+d(y,z))^{m+1}}
		+\frac{t}{(t+d(y,z))^{n+1}}\right],
		\quad y,z\in E_n,\quad t>0.
	\end{equation}
	The subordination formula is
	\begin{equation}\label{eq:subordination}
		P_t(y,z)=\frac{t}{2\sqrt\pi}\int_0^\infty
		e^{-t^2/(4s)}h_s(y,z)\,\frac{ds}{s^{3/2}}.
	\end{equation}
	The scalar subordination identity is
	$$
	e^{-t\sqrt\lambda}
	=\frac{t}{2\sqrt\pi}\int_0^\infty
	e^{-t^2/(4s)}e^{-s\lambda}\,\frac{ds}{s^{3/2}},
	\qquad \lambda\geq0.
	$$
	The spectral theorem and Tonelli give \eqref{eq:subordination}.
	Since $h_s\geq0$, restriction to $s\in[t^2,2t^2]$ gives a lower
	bound. The subordination weight has integral 1, so the sub-Markov
	property of $e^{-s\Delta}$ gives positivity and $L^\infty$
	contractivity of $P_t$.
	
	\begin{lemma}\label{lem:annulus-extension}
		Choose $\vartheta\in C_c^\infty((1/2,3))$ such that
		$0\leq\vartheta\leq1$ and $\vartheta=1$ on $[1,2]$.
		For sufficiently large $R$, define
		$$
		f_R(z)=
		\begin{cases}
			\vartheta(r(z)/R),&z\in E_n\text{ in the region }r(z)>r_0,\\
			0,&\text{elsewhere on }M.
		\end{cases}
		$$
		Then $f_R\in C_c^\infty(M)$ and, for every $0<p<\infty$,
		\begin{equation}\label{eq:fR-norm}
			c_pR^n\leq\|f_R\|_p^p\leq C_pR^n.
		\end{equation}
		For each fixed $A\geq1$, there is $c_A>0$ such that
		\begin{equation}\label{eq:u-lower}
			P_{AR}f_R(y)\geq c_A,\qquad y\in Y_R.
		\end{equation}
		There is a constant $C_1$, independent of $R$ and $B\geq1$, such that
		\begin{equation}\label{eq:u-upper}
			P_{BR}f_R(y)\leq C_1 B^{-n},\qquad y\in Y_R.
		\end{equation}
	\end{lemma}
	
	\begin{proof}
		Since $R/2>r_0$, the function vanishes near the boundary of $E_n$.
		Its extension by zero is smooth and has compact support.
		Polar coordinates give
		\begin{equation}\label{eq:fR-exact-norm}
			\|f_R\|_p^p
			=v_{m-n}n\omega_n R^n
			\int_{1/2}^{3}|\vartheta(s)|^p s^{n-1}\,ds.
		\end{equation}
		Since $\vartheta=1$ on $[1,2]$ and $0\leq\vartheta\leq1$,
		this proves \eqref{eq:fR-norm}.
		
		Fix $A\geq1$, put $t=AR$, and take $y,z\in Y_R$.
		On the interval $t^2\leq s\leq2t^2$, estimate
		\eqref{eq:heat-lower} and \eqref{eq:annuli-distance} give
		$$
		h_s(y,z)\geq C^{-1}(2A^2R^2)^{-n/2}
		e^{-CD^2/A^2}\geq c'_A R^{-n}.
		$$
		Using \eqref{eq:subordination} on this interval, we obtain
		$$
		P_{AR}(y,z)
		\geq\frac{c'_AR^{-n}}{2\sqrt\pi}
		t\int_{t^2}^{2t^2}e^{-t^2/(4s)}\frac{ds}{s^{3/2}}
		=\frac{c'_AR^{-n}}{2\sqrt\pi}
		\int_1^2 e^{-1/(4v)}\frac{dv}{v^{3/2}}
		\geq c''_A R^{-n}.
		$$
		Integrating over $z\in Y_R$, where $f_R=1$, proves
		\eqref{eq:u-lower}.
		
		For the upper bound, \eqref{eq:poisson-upper} gives
		$$
		P_{BR}f_R(y)
		\leq C\big((BR)^{-m}+(BR)^{-n}\big)\|f_R\|_1
		\leq C'\big(B^{-m}R^{n-m}+B^{-n}\big)
		\leq2C'B^{-n},
		$$
		since $R\geq1$, $B\geq1$, and $m>n$. This proves
		\eqref{eq:u-upper} with $C_1=2C'$.
	\end{proof}
	
	\subsection{Regularity and the \texorpdfstring{$L^2$}{L2} energy identity}
	We use the self-adjoint operator associated with the closure of
	$g\mapsto\int_M|\nabla g|^2\,d\mu$, initially defined on
	$C_c^\infty(M)$.
	For nonnegative integers $j,k$ and $t>0$, the spectral theorem gives
	$$
	\|\partial_t^j\Delta^kP_tf\|_2
	\leq\sup_{\lambda\geq0}
	\lambda^{k+j/2}e^{-t\sqrt\lambda}\,\|f\|_2
	\leq C_{j,k}t^{-2k-j}\|f\|_2.
	$$
	Here $C_{j,k}=\sup_{q\geq0}q^{2k+j}e^{-q}$, with $C_{0,0}=1$.
	Local elliptic regularity and Sobolev embedding give a smooth
	representative of $P_tf$ on $M\times(0,\infty)$, satisfying
	$$
	\partial_tP_tf=-\sqrt\Delta P_tf,\qquad
	\partial_t^2P_tf=\Delta P_tf.
	$$
	We use this representative throughout.
	For $g$ in the form domain,
	$$
	\|\sqrt\Delta g\|_2^2=\int_M|\nabla g|^2\,d\mu.
	$$
	Its $L^2$ nullspace is zero. Indeed, if $g\in\ker\Delta$,
	then $\nabla g=0$ in the weak sense, so $g$ is constant on the
	connected manifold. Since $\mu(M)=\infty$, that constant is zero.
	It follows from the spectral theorem that
	\begin{equation}\label{eq:semigroup-limits}
		P_tf\longrightarrow f\text{ in }L^2\text{ as }t\downarrow0,
		\qquad
		P_tf\longrightarrow0\text{ in }L^2\text{ as }t\to\infty.
	\end{equation}
	Let $\sigma_f$ be the spectral measure of $f$.
	Both limits follow by dominated convergence, since
	$\sigma_f(\{0\})=0$. Also, $\|P_tf\|_2$ decreases with $t$, and
	$$
	\|\nabla_yP_tf\|_2^2
	=\|\partial_tP_tf\|_2^2
	=\int_{[0,\infty)}\lambda e^{-2t\sqrt\lambda}\,d\sigma_f(\lambda).
	$$
	Tonelli's theorem now gives
	\begin{align}
		\int_0^\infty t\big(\|\nabla P_tf\|_2^2
		+\|\partial_tP_tf\|_2^2\big)\,dt
		&=2\int_{(0,\infty)}\int_0^\infty
		t\lambda e^{-2t\sqrt\lambda}\,dt\,d\sigma_f(\lambda)
		=\frac12\|f\|_2^2.
		\label{eq:energy}
	\end{align}
	For every $\lambda>0$, the change of variables
	$q=2t\sqrt\lambda$ gives
	$$
	\int_0^\infty t\lambda e^{-2t\sqrt\lambda}\,dt
	=\frac14\int_0^\infty qe^{-q}\,dq=\frac14.
	$$
	The spatial and time derivative terms each equal $\|f\|_2^2/4$.
	By Tonelli and \eqref{eq:area-definition},
	$$
	\begin{aligned}
		\|S(P_tf)\|_2^2
		&=\int_0^\infty\int_M t|\nabla_{y,t}P_tf(y)|^2
		\frac{\mu\{x:d(x,y)<t\}}{V(y,t)}\,d\mu(y)\,dt\\
		&=\int_0^\infty t\big(\|\nabla P_tf\|_2^2
		+\|\partial_tP_tf\|_2^2\big)\,dt
		=\frac12\|f\|_2^2.
	\end{aligned}
	$$
	Here $\mu\{x:d(x,y)<t\}=V(y,t)$. This proves
	\eqref{eq:area-L2-intro}.
	
	\section{Counterexamples for \texorpdfstring{$1<p<2$}{1 < p < 2}}
	
	\begin{proof}[Proof of Theorem \ref{thm:counterexamples}]
		We use the functions in Lemma \ref{lem:annulus-extension}.
		
		\smallskip
		\noindent\emph{The non-tangential maximal function.}
		Fix $\beta>0$ and choose a fixed $A\geq1$ with $\beta A>D$,
		where $D$ is the constant in \eqref{eq:annuli-distance}.
		For every $x\in X_R$ and $y\in Y_R$,
		$$
		d(x,y)\leq DR<\beta AR.
		$$
		Thus $(y,AR)\in\Gamma^\beta(x)$, and
		\eqref{eq:u-lower} implies
		$$
		\mathcal N_P^\beta f_R(x)\geq P_{AR}f_R(y)\geq c_A,
		\qquad x\in X_R.
		$$
		Taking the level $c_A/2$ in the weak norm gives
		$$
		\|\mathcal N_P^\beta f_R\|_{L^{p,\infty}}
		\geq\frac{c_A}{2}\mu(X_R)^{1/p}\geq c_{p,\beta}R^{m/p}.
		$$
		Divide by the upper bound in \eqref{eq:fR-norm} to obtain
		\eqref{eq:maximal-ratio}.
		
		\smallskip
		\noindent\emph{The area function.}
		First fix $A>\max(1,D)$ and let $c_A$ be as in
		\eqref{eq:u-lower}. Then set
		$$
		B=\max\big\{2A,(2C_1/c_A)^{1/n}\big\}.
		$$
		Both $A$ and $B$ are independent of $R$, and
		$C_1B^{-n}\leq c_A/2$. Thus, for $y\in Y_R$,
		$$
		u_R(y,AR)-u_R(y,BR)\geq c_A/2.
		$$
		The extension is smooth for $t>0$. Therefore
		$$
		(c_A/2)^2
		\leq\left|\int_{AR}^{BR}\partial_tu_R(y,t)\,dt\right|^2
		\leq\left(\int_{AR}^{BR}t|\partial_tu_R(y,t)|^2\,dt\right)
		\left(\int_{AR}^{BR}\frac{dt}{t}\right).
		$$
		Consequently,
		\begin{equation}\label{eq:time-energy-lower}
			\int_{AR}^{BR}t|\partial_tu_R(y,t)|^2\,dt
			\geq\frac{c_A^2}{4\log(B/A)}=:c_2>0,
			\qquad y\in Y_R.
		\end{equation}
		
		For $x\in X_R$, $y\in Y_R$, and $AR\leq t\leq BR$,
		\eqref{eq:annuli-distance} gives $d(x,y)<t$.
		Also, by \eqref{eq:volume-upper},
		$$
		V(y,t)\leq C(BR)^m.
		$$
		We may therefore restrict the area integral to this cylinder and
		retain only the time derivative. Using \eqref{eq:time-energy-lower}, for $x\in X_R$,
		$$
		S(u_R)(x)^2
		\geq\frac{1}{CB^mR^m}\int_{Y_R}\int_{AR}^{BR}
		t|\partial_tu_R(y,t)|^2\,dt\,d\mu(y)\geq\frac{c_2\mu(Y_R)}{CB^mR^m}
		\geq c_3R^{n-m}.
		$$
		Put $\lambda_R=\tfrac12\sqrt{c_3}R^{(n-m)/2}$.
		Then $X_R\subset\{S(u_R)>\lambda_R\}$ and
		$$
		\|S(u_R)\|_{L^{p,\infty}}
		\geq\lambda_R\mu(X_R)^{1/p}
		\geq c_pR^{(n-m)/2+m/p}.
		$$
		After division by $\|f_R\|_p\leq C_pR^{n/p}$, this is
		\eqref{eq:area-ratio}.
		
		For $1<p<2$, both powers in \eqref{eq:maximal-ratio} and
		\eqref{eq:area-ratio} are positive. Letting $R\to\infty$ disproves
		weak type $(p,p)$. Since $\|g\|_{L^{p,\infty}}\leq\|g\|_p$,
		neither operator is bounded on $L^p$.
	\end{proof}
	
	\begin{remark}\label{rem:ranges}
		The maximal-function calculation is valid for every $0<p<\infty$,
		and the area-function calculation is valid for every $0<p<2$.
		At $p=2$, the latter power is zero, in agreement with
		\eqref{eq:area-L2-intro}.
		The Markov property gives
		$\|\mathcal N_P^\beta f\|_\infty\leq\|f\|_\infty$.
		Thus the ordinary maximal function remains bounded on $L^\infty$.
	\end{remark}

	\section{Bounds for the restricted cones}
	\label{sec:restricted-bounds}
	
	Fix
	$
	0<c\leq\frac1{100}.
	$
	The function $\rho$ is 1-Lipschitz and is bounded below by 1.
	Define
	\begin{equation}\label{eq:restricted-cones}
		\begin{split}
			\Gamma_c(x)
			&=\{(y,t):t>0,\ d(x,y)<\min(t,c\rho(x))\},\\
			\Gamma_c^\#(x)
			&=\{(y,t):t>0,\ d(x,y)<\min(2t,16c\rho(x))\}.
		\end{split}
	\end{equation}
	The first region is an ordinary cone for $t<c\rho(x)$ and has
	horizontal radius $c\rho(x)$ for $t\geq c\rho(x)$. There is no
	upper bound on $t$. The second region has both a larger aperture
	and a larger horizontal bound. Put
	\begin{equation}\label{eq:restricted-operators}
		\begin{split}
			\mathcal N_cf(x)&=\sup_{(y,t)\in\Gamma_c(x)}|P_tf(y)|,\qquad
			\mathcal N_c^\#f(x)=\sup_{(y,t)\in\Gamma_c^\#(x)}|P_tf(y)|,\\
			S_c(u)(x)&=\left(\iint_{\Gamma_c(x)}
			|\nabla_{y,t}u(y,t)|^2\frac{t\,d\mu(y)\,dt}{V(y,t)}\right)^{1/2}.
		\end{split}
	\end{equation}
	The denominator is $V(y,t)$ for all $t>0$, including
	$t\geq c\rho(x)$.
	
	\begin{theorem}\label{thm:restricted-bounds}
		For each fixed $0<c\leq1/100$, the following assertions hold.
		For $f\in L^1(M)$ and $\lambda>0$,
		\begin{equation}\label{eq:restricted-weak-intro}
			\mu\{\mathcal N_cf>\lambda\}
			\leq\mu\{\mathcal N_c^\#f>\lambda\}
			\leq\frac{C_c}{\lambda}\|f\|_1.
		\end{equation}
		For $1<p\leq\infty$,
		\begin{equation}\label{eq:restricted-maximal-intro}
			\|\mathcal N_cf\|_p\leq\|\mathcal N_c^\#f\|_p
			\leq C_{p,c}\|f\|_p.
		\end{equation}
		For $1<p\leq2$,
		\begin{equation}\label{eq:restricted-area-intro}
			\|S_c(P_tf)\|_p\leq C_{p,c}\|f\|_p.
		\end{equation}
		At $p=2$, one may take the constant in
		\eqref{eq:restricted-area-intro} to be $1/\sqrt2$.
		The definitions and estimates extend from smooth compactly supported
		functions to the indicated $L^p$ spaces.
	\end{theorem}
	
	The triangle inequality gives
	\begin{equation}\label{eq:rlp-rho-Lipschitz}
		|\rho(x)-\rho(z)|\leq d(x,z).
	\end{equation}
	Both product ends have nonnegative Ricci curvature. Since the
	metric differs from the product metrics only on $K$, the Ricci
	curvature of $M$ has a fixed lower bound.
	
	\subsection{Local estimates for harmonic functions}
	
	On $M\times\mathbb R$ we use the product metric. Thus
	$\partial_t^2v=\Delta v$ is the harmonic equation, and $\nabla v$
	denotes the full product gradient in this subsection.
	
	\begin{lemma}\label{lem:local-harmonic-gradient}
		Let $y\in M$, $t>0$, and
		$$
		0<r\leq\min\{t/4,\rho(y)/16\},\qquad
		Q=B_{M\times\mathbb R}((y,t),2r).
		$$
		There is a constant $C$, depending only on $M$, with the following
		properties. If $v$ is positive and harmonic on $Q$, then
		\begin{equation}\label{eq:rlp-positive-gradient}
			|\nabla v(y,t)|\leq \frac{C}{r}v(y,t).
		\end{equation}
		If $v$ is complex valued and harmonic on $Q$, then
		\begin{equation}\label{eq:rlp-interior-gradient}
			|\nabla v(y,t)|\leq \frac{C}{r}\sup_Q|v|.
		\end{equation}
		Consequently, every positive harmonic function on
		$M\times(0,\infty)$ satisfies
		\begin{equation}\label{eq:rlp-log-gradient}
			|\nabla\log v(y,t)|
			\leq C\left(\frac1t+\frac1{\rho(y)}\right).
		\end{equation}
	\end{lemma}
	
	\begin{proof}
		The Cheng--Yau estimate \cite{CY}, also stated in
		\cite[Theorem A]{WZ}, gives
		$$
		|\nabla\log v(q)|\leq C_D(r^{-1}+\sqrt\kappa)
		$$
		for a positive harmonic function on $B(q,2r)$ in a complete
		$D$-dimensional manifold with $\operatorname{Ric}\geq-(D-1)\kappa$.
		Here $D=m+1$ and $2r\leq t/2$.
		
		If $d(y,K)\geq1$, then
		$$
		2r\leq\rho(y)/8\leq d(y,K)/4.
		$$
		A curve of length less than $2r$ from $y$, in either $M$ or its
		complete model end, cannot reach the removed compact set. Hence $Q$
		is isometric to a ball in
		$\mathbb R^m\times\mathbb R$ or
		$\mathbb R^n\times\mathbb S^{m-n}\times\mathbb R$, where
		$\kappa=0$. If $d(y,K)<1$, then $r<1/8$, and the fixed curvature
		term on $M\times\mathbb R$ is bounded by $C/r$.
		This proves \eqref{eq:rlp-positive-gradient}.
		
		For complex $v$, suppose $A=\sup_Q|v|<\infty$ and apply this
		estimate to $A+\varepsilon+\operatorname{Re}v$ and
		$A+\varepsilon+\operatorname{Im}v$. Their values at $(y,t)$ are
		at most $2A+\varepsilon$. Letting $\varepsilon\downarrow0$ gives
		\eqref{eq:rlp-interior-gradient}; the case $A=\infty$ is immediate.
		Finally,
		$$
		r=\min\{t/4,\rho(y)/16\},\qquad r^{-1}\leq4/t+16/\rho(y)
		$$
		gives \eqref{eq:rlp-log-gradient}.
	\end{proof}
	
	\subsection{The restricted maximal functions}
	
	\begin{prop}\label{prop:restricted-maximal}
		Let
		$$
		P_*g(x)=\sup_{t>0}|P_tg(x)|,\qquad
		H_*g(x)=\sup_{s>0}|e^{-s\Delta}g(x)|.
		$$
		For every $f\in L^p(M)$, $1\leq p<\infty$, one has
		\begin{equation}\label{eq:rlp-maximal-domination}
			\mathcal N_c f(x)\leq\mathcal N_c^\# f(x)
			\leq C P_*|f|(x)\leq C H_*|f|(x)
			\quad\hbox{for almost every }x.
		\end{equation}
		In particular,
		\begin{equation}\label{eq:rlp-weak-maximal}
			\mu\{\mathcal N_c^\# f>\lambda\}
			\leq\frac{C}{\lambda}\|f\|_1,
			\qquad f\in L^1(M),\quad\lambda>0,
		\end{equation}
		and
		\begin{equation}\label{eq:rlp-strong-maximal}
			\|\mathcal N_c^\# f\|_p\leq C_p\|f\|_p,
			\qquad 1<p\leq\infty.
		\end{equation}
		The same bounds hold for $\mathcal N_c$.
	\end{prop}
	
	\begin{proof}
		
		Since $\Gamma_c(x)\subset\Gamma_c^\#(x)$, we have
		$\mathcal N_cf(x)\leq\mathcal N_c^\#f(x)$.
		
		First let $f\in L^2(M)$. By positivity, $|P_tf|\leq P_t|f|$, and
		$v(y,t)=P_t|f|(y)+\varepsilon$ is positive and harmonic for
		$\varepsilon>0$. If $(y,t)\in\Gamma_c^\#(x)$, a minimizing
		geodesic $\gamma$ from $x$ to $y$, parameterized by length, has
		length $\ell=d(x,y)$. Then
		$$
		\ell<\min\{2t,16c\rho(x)\},\qquad
		\rho(\gamma(s))\geq(1-16c)\rho(x).
		$$
		Thus \eqref{eq:rlp-log-gradient} gives
		$$
		\begin{aligned}
			|\log v(y,t)-\log v(x,t)|
			&\leq C\int_0^\ell
			\left(\frac1t+\frac1{\rho(\gamma(s))}\right)ds
			\leq C\left(2+\frac{16c}{1-16c}\right)\leq C.
		\end{aligned}
		$$
		
		Taking exponentials gives the comparison at the same height $t$:
		$$
		P_t|f|(y)+\varepsilon
		\leq C\bigl(P_t|f|(x)+\varepsilon\bigr).
		$$
		The constant is independent of $x,y,t,f$ and $\varepsilon$.
		Letting $\varepsilon\downarrow0$ and using positivity, we obtain
		$$
		|P_tf(y)|
		\leq P_t|f|(y)
		\leq CP_t|f|(x)
		\leq CP_*|f|(x).
		$$
		Taking the supremum over $(y,t)\in\Gamma_c^\#(x)$ proves
		$\mathcal N_c^\#f(x)\leq CP_*|f|(x)$.

		The subordination density in \eqref{eq:subordination} has integral
		one. For $g\geq0$,
		$$
		P_tg(x)=\frac{t}{2\sqrt\pi}\int_0^\infty
		e^{-t^2/(4s)}e^{-s\Delta}g(x)\,\frac{ds}{s^{3/2}}
		\leq H_*g(x).
		$$
		By \cite[Theorem 3]{DLS}, $H_*$ is of weak type $(1,1)$ and
		bounded on $L^p$ for $1<p\leq\infty$. This proves
		\eqref{eq:rlp-weak-maximal} and \eqref{eq:rlp-strong-maximal}
		for $f$ also in $L^2$.
		
		The nonnegative symmetric sub-Markov kernels define consistent
		contractions on $L^q$, $1\leq q\leq\infty$. Indeed, for finite $q$,
		$$
		|P_tf(x)|^q\leq\int_M P_t(x,z)|f(z)|^q\,d\mu(z),\qquad
		\|P_tf\|_q\leq\|f\|_q;
		$$
		the same argument applies to the heat semigroup. For $q=\infty$,
		use the kernel mass bound.
		
		Let $1\leq q<\infty$ and let $F\subset M\times(0,\infty)$ be
		compact. We claim that
		\begin{equation}\label{eq:rlp-compact-evaluation}
			\sup_{(y,t)\in F}|P_th(y)|\leq C_{F,q}\|h\|_q,
			\qquad h\in L^q(M)\cap L^2(M).
		\end{equation}
		Choose $a>0$ with $t\geq a$ on $F$ and put $b=\min\{a,c\}/4$.
		Then
		$$
		x\in B(y,b)\quad\Longrightarrow\quad
		\mathcal N_ch(x)\geq|P_th(y)|,\qquad
		\inf_{(y,t)\in F}\mu(B(y,b))>0.
		$$
		For $q>1$, integrate over $B(y,b)$ and apply
		\eqref{eq:rlp-strong-maximal}. For $q=1$,
		\eqref{eq:rlp-weak-maximal} at level $|P_th(y)|/2$ gives
		$$
		\tfrac12|P_th(y)|\mu(B(y,b))\leq C\|h\|_1.
		$$
		This proves \eqref{eq:rlp-compact-evaluation}.
		
		Choose $f_j\in C_c^\infty(M)$ with $f_j\to f$ in $L^q$.
		By \eqref{eq:rlp-compact-evaluation}, $P_tf_j$ converges uniformly
		on compact subsets. Its limit is distributionally harmonic, hence
		smooth by local elliptic regularity. At each $t>0$, contractivity
		identifies it with $P_tf$ in $L^q$. Applying
		\eqref{eq:rlp-interior-gradient} to differences on finitely many
		slightly larger product balls gives
		\begin{equation}\label{eq:rlp-compact-convergence}
			P_tf_j\longrightarrow P_tf,\qquad
			\nabla_{y,t}P_tf_j\longrightarrow\nabla_{y,t}P_tf
			\quad\hbox{uniformly on compact subsets.}
		\end{equation}
		The doubled balls can be chosen inside $M\times(0,\infty)$ with
		radii satisfying Lemma~\ref{lem:local-harmonic-gradient}.
		
		For every $x$,
		$$
		\mathcal N_c^\#f(x)\leq\liminf_j\mathcal N_c^\#f_j(x).
		$$
		Fatou proves the strong bound for $1<q<\infty$. For $q=1$ and
		$0<\lambda'<\lambda$, it gives
		$$
		\mu\{\mathcal N_c^\#f>\lambda\}
		\leq\liminf_j\mu\{\mathcal N_c^\#f_j>\lambda'\}
		\leq\frac C{\lambda'}\|f\|_1.
		$$
		Let $\lambda'\uparrow\lambda$. The argument with
		$P_t|f|+\varepsilon$ now gives \eqref{eq:rlp-maximal-domination}
		for every finite $q$.
		
		For $f\in L^\infty$, the kernel mass bound gives
		$|P_tf(y)|\leq\|f\|_\infty$ and
		$\mathcal N_c^\#f\leq\|f\|_\infty$.
		To obtain a harmonic representative, restrict each nonnegative
		real or imaginary part of $f$ to increasing compact sets.
		The Poisson extensions increase to the kernel integral and remain
		bounded. The interior gradient estimate gives local equicontinuity;
		the harmonic equation passes to the limit.
	\end{proof}
	
\begin{corollary}\label{cor:restricted-boundary}
Let $f\in L^p(M)$, $1\leq p\leq\infty$. Then
\begin{equation}\label{eq:Hardy-boundary-limit}
 P_tf(y)\longrightarrow f(x)
 \quad\text{as }t\downarrow0,\ d(x,y)<t,
 \quad\text{for almost every }x\in M.
\end{equation}
In particular, this holds for every $f\in H^1_{S_c}(M)$.
\end{corollary}

\begin{proof}
For $h\in C_c^\infty(M)$, the heat equation and contractivity give
$$
 e^{-s\Delta}h-h=-\int_0^s e^{-v\Delta}\Delta h\,dv,\qquad
 \|e^{-s\Delta}h-h\|_\infty
 \leq\min\{s\|\Delta h\|_\infty,2\|h\|_\infty\}.
$$
Subordination gives
\begin{align*}
 \|P_th-h\|_\infty
 &\leq\frac{t}{2\sqrt\pi}\int_0^\infty
 \frac{\min\{s\|\Delta h\|_\infty,2\|h\|_\infty\}}
 {s^{3/2}}\,ds
 \leq Ct\bigl(\|\Delta h\|_\infty+\|h\|_\infty\bigr).
\end{align*}
Also,
$$
 d(x,y)<t\quad\Longrightarrow\quad
 |h(y)-h(x)|\leq t\|\nabla h\|_\infty.
$$
Thus \eqref{eq:Hardy-boundary-limit} holds for $h$ at every $x$.

Let $1\leq p<\infty$ and put
$$
 Df(x)=\limsup_{\substack{t\downarrow0\\d(x,y)<t}}
 |P_tf(y)-f(x)|.
$$
For each fixed $x$,
$$
 0<t<c\rho(x),\quad d(x,y)<t
 \quad\Longrightarrow\quad (y,t)\in\Gamma_c(x).
$$
Hence, for $h\in C_c^\infty(M)$,
$$
 Df(x)\leq\mathcal N_c(f-h)(x)+|f(x)-h(x)|.
$$
Proposition~\ref{prop:restricted-maximal} and Chebyshev's inequality
give, for $\lambda>0$,
$$
 \mu\{Df>2\lambda\}\leq
 \begin{cases}
 C_c\lambda^{-1}\|f-h\|_1,&p=1,\\[1mm]
 C_{p,c}\lambda^{-p}\|f-h\|_p^p,&1<p<\infty.
 \end{cases}
$$
Choose $h$ tending to $f$ in $L^p$. Then $Df=0$ almost everywhere.

Now let $f\in L^\infty(M)$. Choose relatively compact open sets
$U_j$ covering $M$ and $\chi_j\in C_c^\infty(M)$ such that
$0\leq\chi_j\leq1$ and $\chi_j=1$ on a neighbourhood of
$\overline{U_j}$. The result for $p=1$ applies to $\chi_jf$.
If $x\in U_j$, then $\chi_j(y)=1$ whenever $d(x,y)<t$ and $t$
is sufficiently small. Positivity and the kernel mass bound give
\begin{align*}
 |P_t((1-\chi_j)f)(y)|
 &\leq\|f\|_\infty P_t(1-\chi_j)(y)
 \leq\|f\|_\infty\bigl(1-P_t\chi_j(y)\bigr)
 \leq C_{\chi_j}t\|f\|_\infty.
\end{align*}
Therefore \eqref{eq:Hardy-boundary-limit} holds almost everywhere
on each $U_j$, and hence on $M$.
The assertion for $H^1_{S_c}(M)$ follows from its embedding into
$L^1(M)$ in Theorem~\ref{thm:hardy-characterization}(i).
\end{proof}

	\subsection{A nonlinear energy identity}
	
	\begin{lemma}\label{lem:rlp-nonlinear-energy}
		Let $1<p<2$ and let $f\geq0$ belong to $L^p(M)\cap L^2(M)$.
		If $u=P_tf$, then
		\begin{equation}\label{eq:rlp-nonlinear-energy}
			\int_0^\infty\int_M
			t\,u^{p-2}|\nabla_{y,t}u|^2\,d\mu\,dt
			=\frac{\|f\|_p^p}{p(p-1)}.
		\end{equation}
		If $f=0$, the integrand is defined to be zero.
	\end{lemma}
	
	\begin{proof}
		Suppose $f\ne0$. The function $u$ is nonnegative and harmonic,
		and $P_tf\to f$ in $L^2$ at zero. The strong maximum principle
		gives $u>0$ on the connected product $M\times(0,\infty)$.
		
		For $\delta>0$ and $s\geq0$, put
		$$
		\Phi_\delta(s)=(s+\delta)^p-\delta^p-p\delta^{p-1}s,
		\qquad L_\delta=p(p-1)\delta^{p-2}.
		$$
		Then
		$$
		\Phi_\delta(0)=\Phi_\delta'(0)=0,\qquad
		0\leq\Phi_\delta''(s)=p(p-1)(s+\delta)^{p-2}\leq L_\delta,
		$$
		and hence
		\begin{equation}\label{eq:rlp-Phi-bounds}
			0\leq\Phi_\delta(s)\leq\tfrac12L_\delta s^2,
			\qquad 0\leq\Phi_\delta'(s)\leq L_\delta s.
		\end{equation}
		Also
		$$
		\Phi_\delta(s)=p(p-1)\int_0^s(s-a)(a+\delta)^{p-2}\,da
		\uparrow s^p\quad(\delta\downarrow0).
		$$
		
		Put $F_\delta(t)=\int_M\Phi_\delta(u(y,t))\,d\mu(y)$.
		Extend $\Phi_\delta$ to $s<0$ by $L_\delta s^2/2$; then
		$\Phi_\delta\in C^2(\mathbb R)$ and
		$\|\Phi_\delta''\|_\infty\leq L_\delta$.
		The $L^2$ spectral estimates preceding \eqref{eq:semigroup-limits}
		and Taylor's formula give
		$$
		F_\delta'(t)=\int_M\Phi_\delta'(u)u_t\,d\mu
		$$
		and
		\begin{equation}\label{eq:rlp-F-second-first}
			F_\delta''(t)=\int_M\Phi_\delta''(u)|u_t|^2\,d\mu
			+\int_M\Phi_\delta'(u)u_{tt}\,d\mu.
		\end{equation}
		To justify the second derivative, write
		$$
		m_h=\int_0^1\Phi_\delta''\big(u(t)+\theta(u(t+h)-u(t))\big)\,d\theta.
		$$
		Then $|m_h|\leq L_\delta$ and
		$m_h\to\Phi_\delta''(u(t))$ in measure. Since
		$(u(t+h)-u(t))/h\to u_t(t)$ in $L^2$, dominated convergence
		along almost everywhere convergent subsequences gives
		$$
		\frac{\Phi_\delta'(u(t+h))-\Phi_\delta'(u(t))}{h}
		=m_h\frac{u(t+h)-u(t)}h
		\longrightarrow\Phi_\delta''(u(t))u_t(t)\quad\hbox{in }L^2.
		$$
		This also proves continuity of \eqref{eq:rlp-F-second-first}.
		
		Since $\Phi_\delta'$ is Lipschitz and vanishes at zero,
		$\Phi_\delta'(u)$ belongs to the quadratic form domain. The Sobolev
		chain rule and $u_{tt}=\Delta u$ give
		$$
		\int_M\Phi_\delta'(u)u_{tt}\,d\mu
		=\int_M\nabla_y\Phi_\delta'(u)\cdot\nabla_yu\,d\mu
		=\int_M\Phi_\delta''(u)|\nabla_yu|^2\,d\mu.
		$$
		The chain rule follows by approximation in the form norm:
		$\Phi_\delta'$ is Lipschitz, $\Phi_\delta''$ is bounded and continuous,
		and multiplication by the fixed $L^2$ function $\nabla_yu$ preserves
		the convergence of these bounded multipliers. Therefore
		\begin{equation}\label{eq:rlp-F-second}
			F_\delta''(t)=p(p-1)\int_M
			(u+\delta)^{p-2}|\nabla_{y,t}u|^2\,d\mu\geq0.
		\end{equation}
		
		By \eqref{eq:rlp-Phi-bounds},
		$$
		0\leq F_\delta(t)\leq\tfrac12L_\delta\|P_tf\|_2^2,
		\qquad
		|\Phi_\delta(a)-\Phi_\delta(b)|
		\leq L_\delta(|a|+|b|)|a-b|.
		$$
		Thus \eqref{eq:semigroup-limits} and Cauchy--Schwarz give
		\begin{equation}\label{eq:rlp-F-limits}
			F_\delta(t)\longrightarrow
			F_\delta(0):=\int_M\Phi_\delta(f)\,d\mu
			\quad(t\downarrow0),\qquad
			F_\delta(t)\longrightarrow0\quad(t\to\infty).
		\end{equation}
		Boundedness and convexity imply $F_\delta'\leq0$. Also,
		$F_\delta(t/2)\geq F_\delta(t)-\tfrac t2F_\delta'(t)$, so
		\begin{equation}\label{eq:rlp-F-slope}
			0\leq-tF_\delta'(t)
			\leq2\big(F_\delta(t/2)-F_\delta(t)\big).
		\end{equation}
		The right-hand side tends to zero at both endpoints. Hence
		$$
		\int_a^b tF_\delta''(t)\,dt
		=bF_\delta'(b)-aF_\delta'(a)-F_\delta(b)+F_\delta(a)
		\longrightarrow F_\delta(0)
		$$
		as $a\downarrow0$ and $b\to\infty$. Using
		\eqref{eq:rlp-F-second}, we obtain
		$$
		p(p-1)\int_0^\infty\int_M
		t(u+\delta)^{p-2}|\nabla_{y,t}u|^2\,d\mu\,dt
		=\int_M\Phi_\delta(f)\,d\mu.
		$$
		Both integrands increase as $\delta\downarrow0$; monotone
		convergence proves \eqref{eq:rlp-nonlinear-energy}.
	\end{proof}
	
	\subsection{The area function}
	
	In this subsection we use the abbreviation
	$S_cf=S_c(P_tf)$ for the restricted area function of the Poisson
	extension of $f$.
	
	The following estimate applies to every measurable subregion
	of the full cone. It uses the normalization by $V(y,t)$.
	
	\begin{lemma}\label{lem:rlp-subcone-area}
		Let $\Gamma_*(x)$ be a family of measurable sets such that
		$\{(x,y,t):(y,t)\in\Gamma_*(x)\}$ is measurable and
		$\Gamma_*(x)\subset\{(y,t):d(x,y)<t\}$. Define
		$$
		\mathcal N_*f(x)=\sup_{(y,t)\in\Gamma_*(x)}|P_tf(y)|,
		\qquad
		S_*f(x)^2=\int_{\Gamma_*(x)}
		|\nabla_{y,t}P_tf(y)|^2\frac{t\,d\mu(y)\,dt}{V(y,t)}.
		$$
		The supremum and the integral over an empty region are both
		defined to be zero.
		Assume also that $\mathcal N_*f$ is measurable. If $1<p<2$,
		$f\geq0$ belongs to $L^p\cap L^2$, and $\mathcal N_*f\in L^p$,
		then
		\begin{equation}\label{eq:rlp-subcone-estimate}
			\|S_*f\|_p
			\leq\frac1{\sqrt{p(p-1)}}
			\|\mathcal N_*f\|_p^{1-p/2}\|f\|_p^{p/2}.
		\end{equation}
	\end{lemma}
	
	\begin{proof}
		Suppose $f\ne0$ and write $u=P_tf>0$. Put
		$$
		A_p(x)=\int_{\Gamma_*(x)}
		u(y,t)^{p-2}|\nabla_{y,t}u(y,t)|^2
		\frac{t\,d\mu(y)\,dt}{V(y,t)}.
		$$
		Since $u^{2-p}\leq(\mathcal N_*f(x))^{2-p}$ on $\Gamma_*(x)$,
		\begin{equation}\label{eq:rlp-pointwise-area}
			S_*f(x)^2\leq\mathcal N_*f(x)^{2-p}A_p(x).
		\end{equation}
		On nonempty regions, $\mathcal N_*f$ is positive and finite almost
		everywhere; empty regions contribute zero.
		
		For fixed $(y,t)$,
		$\{x:(y,t)\in\Gamma_*(x)\}\subset B(y,t)$.
		Tonelli and \eqref{eq:rlp-nonlinear-energy} give
		$$
		\begin{aligned}
			\int_M A_p(x)\,d\mu(x)
			&=\int_0^\infty\int_M t u^{p-2}|\nabla_{y,t}u|^2
			\frac{\mu\{x:(y,t)\in\Gamma_*(x)\}}{V(y,t)}\,d\mu(y)\,dt
			\leq\frac{\|f\|_p^p}{p(p-1)}.
		\end{aligned}
		$$
		H\"older's inequality with exponents $2/(2-p)$ and $2/p$ yields
		$$
		\begin{aligned}
			\|S_*f\|_p^p
			&\leq\int_M(\mathcal N_*f)^{p(2-p)/2}A_p^{p/2}\,d\mu
			\leq\|\mathcal N_*f\|_p^{p(2-p)/2}
			\left(\int_M A_p\,d\mu\right)^{p/2}.
		\end{aligned}
		$$
		Taking the $p$th root proves \eqref{eq:rlp-subcone-estimate}.
		The case $f=0$ is immediate.
	\end{proof}
	
	For $\Gamma_*=\Gamma_c$, the strict continuous cone inequalities
	give measurability of the family and lower semicontinuity of
	$\mathcal N_cf$: each fixed cone point remains admissible at
	nearby vertices. Proposition~\ref{prop:restricted-maximal} gives
	\begin{equation}\label{eq:rlp-area-positive}
		\|S_cf\|_p\leq C_p\|f\|_p,
		\qquad f\geq0,\quad f\in L^p\cap L^2,\quad1<p<2.
	\end{equation}
	
	For complex data in $L^p\cap L^2$, write
	$$
	f=f_1-f_2+i f_3-i f_4,
	\quad
	f_1=(\operatorname{Re}f)_+,\quad
	f_2=(\operatorname{Re}f)_-,\quad
	f_3=(\operatorname{Im}f)_+,\quad
	f_4=(\operatorname{Im}f)_-.
	$$
	By linearity, the triangle inequality on each cone, and
	\eqref{eq:rlp-area-positive},
	$$
	S_cf\leq\sum_{j=1}^4S_cf_j,\qquad
	\|S_cf\|_p\leq C_p\sum_{j=1}^4\|f_j\|_p\leq4C_p\|f\|_p.
	$$
	For $p=2$, cone inclusion and \eqref{eq:energy} give the more
	precise estimate
	\begin{equation}\label{eq:rlp-area-two}
		\|S_cf\|_2^2\leq\frac12\|f\|_2^2.
	\end{equation}
	
	Finally, let $1<p<2$ and $f\in L^p$. Choose
	$f_j\in C_c^\infty(M)$ converging to $f$ in $L^p$.
	By \eqref{eq:rlp-compact-convergence} and Fatou, first on each
	cone and then on $M$,
	$$
	S_cf(x)^2\leq\liminf_j S_cf_j(x)^2,
	\qquad
	\|S_cf\|_p\leq C_p\|f\|_p.
	$$
	Here $S_cf$ uses the smooth extension in
	Proposition~\ref{prop:restricted-maximal}. Applying this bound
	to $f_j-f$ gives
	$$
	\|S_cf_j-S_cf\|_p\leq\|S_c(f_j-f)\|_p
	\leq C_p\|f_j-f\|_p\longrightarrow0.
	$$
	Thus $S_c$ extends continuously from $C_c^\infty(M)$, proving
	Theorem~\ref{thm:restricted-bounds}.
	
	\begin{remark}\label{rem:rlp-height-truncation}
		The cone truncated in height,
		$$
		\Gamma_c^{\mathrm{tr}}(x)
		=\{(y,t):d(x,y)<t,\ 0<t<c\rho(x)\},
		$$
		is contained in $\Gamma_c(x)$. Its maximal and area functions
		therefore satisfy all the corresponding bounds proved above.
		The region $\Gamma_c(x)$ also includes every height $t>0$;
		only its horizontal radius is capped by $c\rho(x)$.
	\end{remark}

	\section{The Fefferman--Stein inequality}
	\label{sec:full-fs}\label{sec:rfs}

	\subsection{A cutoff calculation}
	\label{subsec:fs-cutoff}

	Let $f\in C_c^\infty(M)$, $u(y,t)=P_tf(y)$, and $\alpha>0$.
	Let $E\subset M$ be closed, with
	\begin{equation}\label{eq:finite-bad-set}
		\mu(E^c)<\infty.
	\end{equation}
	Then $E\ne\varnothing$. Put $\delta(y)=d(y,E)$ and assume that
	\begin{equation}\label{eq:common-cutoff-bound}
		|u(y,t)|\leq\alpha\qquad\text{if }\delta(y)<2t.
	\end{equation}
	The triangle inequality and closedness of $E$ give
	$$
	|\delta(y)-\delta(z)|\leq d(y,z),\qquad
	|\nabla\delta|\leq1\ \text{a.e.},\qquad
	\{\delta=0\}=E.
	$$
	Choose a smooth nonincreasing function $\eta:[0,\infty)\to[0,1]$
	with
	$$
	\eta(s)=1\quad(0\leq s\leq1),\qquad
	\eta(s)=0\quad(s\geq2),
	$$
	and put
	\begin{equation}\label{eq:psi}
		\psi(y,t)=\eta(\delta(y)/t).
	\end{equation}
	The function $\psi$ is locally Lipschitz for $t>0$.
	We use weak derivatives and the chain rule in local charts.
	
	If $x\in E$ and $d(x,y)<t$, then $\delta(y)<t$, so
	\begin{equation}\label{eq:psi-cones}
		\psi(y,t)=1\quad\text{whenever }E\cap B(y,t)\ne\varnothing.
	\end{equation}
	The derivatives of $\psi$ vanish outside
	\begin{equation}\label{eq:transition}
		\mathcal U=\{(y,t):y\in E^c,\ \delta(y)/2<t<\delta(y)\}.
	\end{equation}
	By \eqref{eq:common-cutoff-bound},
	\begin{equation}\label{eq:u-transition}
		|u(y,t)|\leq\alpha\quad\text{on }\mathcal U.
	\end{equation}
	The chain rule gives almost everywhere
	$$
	\nabla_y\psi=\eta'(\delta/t)\frac{\nabla\delta}{t},\qquad
	\partial_t\psi=-\eta'(\delta/t)\frac{\delta}{t^2}.
	$$
	On $\mathcal U$, $\delta/t<2$, and hence
	$$
	|\nabla_{y,t}\psi|^2\leq\frac{5\|\eta'\|_\infty^2}{t^2}
	\mathbf1_{\mathcal U}.
	$$
	In particular,
	\begin{equation}\label{eq:psi-energy}
		\int_0^\infty t|\nabla_{y,t}\psi(y,t)|^2\,dt
		\leq5\|\eta'\|_\infty^2\log2\,\mathbf1_{E^c}(y).
	\end{equation}
	Since $\eta'\leq0$, we have $\partial_t\psi\geq0$. Integrating
	$\tfrac12\partial_t(\psi^2)$ gives
	\begin{equation}\label{eq:psi-time}
		\int_0^\infty\psi(y,t)\partial_t\psi(y,t)\,dt
		=\frac12\mathbf1_{E^c}(y).
	\end{equation}
	Finally,
	\begin{equation}\label{eq:psi-trace}
		\psi(y,t)^2\longrightarrow\mathbf1_E(y)\quad(t\downarrow0).
	\end{equation}
	
	\subsubsection{Integration by parts with spatial and temporal cutoffs}
	Write $\nabla=(\nabla_y,\partial_t)$ and put
	$$
	\mathsf T=\int_0^\infty\int_M t\psi^2|\nabla u|^2\,d\mu\,dt.
	$$
	By \eqref{eq:energy}, $\mathsf T\leq\tfrac12\|f\|_2^2<\infty$.
	We shall prove the identity
	\begin{align}
		\mathsf T
		={}&\frac12\int_E|f|^2\,d\mu
		-2\operatorname{Re}\int_0^\infty\int_M
		t\psi\overline u\,\nabla u\cdot\nabla\psi\,d\mu\,dt
		+\int_0^\infty\int_M|u|^2\psi\partial_t\psi\,d\mu\,dt.
		\label{eq:tested-identity}
	\end{align}
	We allow complex-valued $f$.
	
	Choose smooth functions $\chi,\kappa:[0,\infty)\to[0,1]$ such that
	$$
	\chi=0\text{ on }[0,1],\quad\chi=1\text{ on }[2,\infty),
	\qquad
	\kappa=1\text{ on }[0,1],\quad\kappa=0\text{ on }[2,\infty).
	$$
	For $T_0>4\varepsilon>0$, set
	$$
	h(t)=h_{\varepsilon,T_0}(t)
	=t\chi(t/\varepsilon)\kappa(t/T_0).
	$$
	Then $h$ has compact support in $(0,\infty)$,
	\begin{equation}\label{eq:h-bounds}
		0\leq h(t)\leq t,\qquad |h'(t)|\leq C,
	\end{equation}
	and
	\begin{equation}\label{eq:h-second}
		|h''(t)|\leq\frac{C}{\varepsilon}
		\mathbf1_{[\varepsilon,2\varepsilon]}(t)
		+\frac{C}{T_0}\mathbf1_{[T_0,2T_0]}(t).
	\end{equation}
	Indeed, the transition intervals are disjoint, and
	$$
	h''(t)=
	\begin{cases}
		2\varepsilon^{-1}\chi'(t/\varepsilon)
		+t\varepsilon^{-2}\chi''(t/\varepsilon),&\varepsilon<t<2\varepsilon,\\
		2T_0^{-1}\kappa'(t/T_0)+tT_0^{-2}\kappa''(t/T_0),&T_0<t<2T_0,\\
		0,&\text{otherwise}.
	\end{cases}
	$$
	This proves \eqref{eq:h-bounds}--\eqref{eq:h-second}. Also
	$h(t)\to t$ and $h'(t)\to1$ as $\varepsilon\downarrow0$ and
	$T_0\to\infty$, since $h(t)=t$ on $[2\varepsilon,T_0]$.
	
	Fix $o\in M$ and set $\zeta_j(y)=\kappa(d(o,y)/j)$. Then
	$\zeta_j=1$ on $B(o,j)$, and
	$\operatorname{supp}\zeta_j\subset\overline{B(o,2j)}$ is compact
	by Hopf--Rinow. Moreover,
	$$
	0\leq\zeta_j\leq1,\qquad |\nabla\zeta_j|\leq C/j\ \text{a.e.},
	\qquad \zeta_j\longrightarrow1\ \text{pointwise}.
	$$
	
	The function $w=h\zeta_j^2\psi^2u$ is a compactly supported
	$W^{1,2}$ function in $M\times(0,\infty)$. A finite partition of
	unity and convolution in coordinate charts give
	$w_k\in C_c^\infty$ with $w_k\to w$ in $W^{1,2}$, all supported
	in one relatively compact set $\Omega$. Hence
	$$
	\left|\int\nabla u\cdot\overline{\nabla(w_k-w)}\right|
	\leq\|\nabla u\|_{L^2(\Omega)}
	\|\nabla(w_k-w)\|_2\longrightarrow0.
	$$
	Thus $w$ is an admissible test function. Integrals without a domain
	are over $M\times(0,\infty)$ with measure $d\mu(y)\,dt$.
	
	The equation $\partial_t^2u-\Delta_yu=0$ gives
	$\int\nabla u\cdot\overline{\nabla w}=0$.
	Substitute $w=h\zeta_j^2\psi^2u$, take real parts, and use
	$\operatorname{Re}(\overline u\partial_tu)
	=\tfrac12\partial_t|u|^2$. Integration by parts in time gives
	\begin{align}
		\int h\zeta_j^2\psi^2|\nabla u|^2
		={}&-2\operatorname{Re}\int h\zeta_j^2\psi\overline u\,
		\nabla u\cdot\nabla\psi
		+\frac12\int h''\zeta_j^2|u|^2\psi^2+\int h'\zeta_j^2|u|^2\psi\partial_t\psi
		+\mathcal E_j,
		\label{eq:finite-test}
	\end{align}
	where
	$$
	\mathcal E_j=-2\operatorname{Re}\int
	h\zeta_j\psi^2\overline u\,
	\nabla_yu\cdot\nabla_y\zeta_j.
	$$
	The time calculation has no boundary term, since $h'$ has
	compact support in $(0,\infty)$:
	$$
	\begin{aligned}
		-\frac12\int h'\zeta_j^2\psi^2\partial_t|u|^2
		&=\frac12\int |u|^2\partial_t(h'\zeta_j^2\psi^2)
		=\frac12\int h''\zeta_j^2|u|^2\psi^2
		+\int h'\zeta_j^2|u|^2\psi\partial_t\psi.
	\end{aligned}
	$$
	
	For fixed $\varepsilon,T_0$, Cauchy--Schwarz yields
	$$
	|\mathcal E_j|\leq\frac{C_{T_0}}{j}
	\left(\int_\varepsilon^{2T_0}\|u(\cdot,t)\|_2^2\,dt\right)^{1/2}
	\left(\int_\varepsilon^{2T_0}\|\nabla_yu(\cdot,t)\|_2^2\,dt\right)^{1/2}
	\longrightarrow0.
	$$
	Here $L^2$ contractivity and \eqref{eq:energy} give
	$$
	\int_\varepsilon^{2T_0}\|u(\cdot,t)\|_2^2\,dt
	\leq2T_0\|f\|_2^2,
	\qquad
	\int_\varepsilon^{2T_0}\|\nabla u(\cdot,t)\|_2^2\,dt
	\leq\frac{\|f\|_2^2}{2\varepsilon}.
	$$
	Thus $|\mathcal E_j|\leq C_{\varepsilon,T_0}\|f\|_2^2/j$.
	The other integrands in \eqref{eq:finite-test} are bounded by
	$C_{\varepsilon,T_0}$ times the integrable function
	$$
	\mathbf1_{[\varepsilon,2T_0]}(t)
	\big(|\nabla u|^2+|u|\,|\nabla u|+|u|^2\big).
	$$
	Dominated convergence as $j\to\infty$ therefore gives
	\begin{align}
		\int h\psi^2|\nabla u|^2
		={}&-2\operatorname{Re}\int h\psi\overline u\,
		\nabla u\cdot\nabla\psi
		+\frac12\int h''|u|^2\psi^2
		+\int h'|u|^2\psi\partial_t\psi.
		\label{eq:time-cutoff-identity}
	\end{align}
	
	\subsubsection{The limits at zero and infinity}
	Define
	$$
	H(t)=\int_M|u(y,t)|^2\psi(y,t)^2\,d\mu(y),\qquad
	H_0=\int_E|f|^2\,d\mu.
	$$
	As $t\downarrow0$, we have
	$$
	\begin{aligned}
		\left|H(t)-\int_M|f|^2\psi(y,t)^2\,d\mu\right|
		&\leq\big\||u(\cdot,t)|^2-|f|^2\big\|_1\\
		&\leq\|u(\cdot,t)-f\|_2
		\big(\|u(\cdot,t)\|_2+\|f\|_2\big)\longrightarrow0.
	\end{aligned}
	$$
	By \eqref{eq:psi-trace} and dominated convergence, the remaining
	integral tends to $H_0$. Thus $H(t)\to H_0$.
	
	On the lower interval in \eqref{eq:h-second},
	$$
	\int_\varepsilon^{2\varepsilon}h''(t)\,dt
	=h'(2\varepsilon)-h'(\varepsilon)=1.
	$$
	Since \eqref{eq:h-second} also bounds $\int|h''|$ uniformly,
	$$
	\left|\int_\varepsilon^{2\varepsilon}h''(t)H(t)\,dt-H_0\right|
	\leq C\sup_{\varepsilon\leq t\leq2\varepsilon}|H(t)-H_0|
	\longrightarrow0.
	$$
	At the other endpoint, \eqref{eq:semigroup-limits} gives
	$$
	\left|\int_{T_0}^{2T_0}h''(t)H(t)\,dt\right|
	\leq C\sup_{t\geq T_0}\|P_tf\|_2^2
	=C\|P_{T_0}f\|_2^2\longrightarrow0.
	$$
	Hence
	\begin{equation}\label{eq:h-second-limit}
		\int h''|u|^2\psi^2\longrightarrow\int_E|f|^2\,d\mu
		\quad\text{as }\varepsilon\downarrow0,\ T_0\to\infty.
	\end{equation}
	
	For the other terms of \eqref{eq:time-cutoff-identity}, use
	$t|\nabla u|^2\in L^1$ from \eqref{eq:energy}, and
	\eqref{eq:u-transition}--\eqref{eq:psi-energy} to obtain
	\begin{align}
		\int t\psi|u|\,|\nabla u|\,|\nabla\psi|
		&\leq\alpha
		\left(\int t|\nabla u|^2\right)^{1/2}
		\left(\int t|\nabla\psi|^2\right)^{1/2}
		\leq C\alpha\|f\|_2\mu(E^c)^{1/2}<\infty.
		\label{eq:mixed-integrable}
	\end{align}
	Also, \eqref{eq:u-transition} and \eqref{eq:psi-time} give
	\begin{equation}\label{eq:last-integrable}
		\int |u|^2\psi\partial_t\psi
		\leq\frac{\alpha^2}{2}\mu(E^c)<\infty.
	\end{equation}
	Together with \eqref{eq:h-bounds} and \eqref{eq:h-second-limit},
	these bounds justify dominated convergence and prove
	\eqref{eq:tested-identity}.
	
	\subsubsection{The energy estimate}
	By Young's inequality,
	$$
	2t\psi|u|\,|\nabla u|\,|\nabla\psi|
	\leq\frac12t\psi^2|\nabla u|^2+2t|u|^2|\nabla\psi|^2.
	$$
	By \eqref{eq:u-transition} and \eqref{eq:psi-energy},
	$\int t|u|^2|\nabla\psi|^2\leq C\alpha^2\mu(E^c)$.
	Thus \eqref{eq:last-integrable} and \eqref{eq:tested-identity} give
	$$
	\mathsf T\leq\frac12\int_E|f|^2\,d\mu
	+\frac12\mathsf T+C\alpha^2\mu(E^c).
	$$
	Since $\mathsf T<\infty$, this yields
	\begin{equation}\label{eq:localized-energy}
		\mathsf T\leq\int_E|f|^2\,d\mu+C\alpha^2\mu(E^c).
	\end{equation}

	\subsection{The restricted cones}

	We use the cones and operators in
	\eqref{eq:restricted-cones}--\eqref{eq:restricted-operators}.
	
	\begin{theorem}\label{thm:restricted-fs}
		Fix $0<c\leq1/100$. There is a constant $C_c>0$ such that,
		for every $f\in C_c^\infty(M)$ and $\lambda>0$,
		\begin{equation}\label{eq:restricted-fs-intro}
			\mu\{S_c(P_tf)>\lambda\}
			\leq C_c\mu\{\mathcal N_cf>\lambda\}
			+\frac{C_c}{\lambda^2}\int_0^\lambda
			s\,\mu\{\mathcal N_cf>s\}\,ds.
		\end{equation}
		Consequently, if $0<p<2$ and $\mathcal N_cf\in L^p(M)$, then
		\begin{equation}\label{eq:restricted-comparison-intro}
			\|S_c(P_tf)\|_p\leq C_{p,c}\|\mathcal N_cf\|_p.
		\end{equation}
	\end{theorem}
	
	Fix $0<c\leq1/100$ and $f\in C_c^\infty(M)$. Write
	$u(y,t)=P_tf(y)$, $S_cf=S_c(P_tf)$, and $\nabla=\nabla_{y,t}$.
	Constants may depend on $c$ and $M$. Define
	\begin{align}
		S_{c,0}f(x)^2
		&=\int_0^{c\rho(x)}\int_{B(x,t)}
		|\nabla u(y,t)|^2\frac{t\,d\mu(y)\,dt}{V(y,t)},
		\label{rfs:lower-area}\\
		S_{c,\infty}f(x)^2
		&=\int_{c\rho(x)}^\infty\int_{B(x,c\rho(x))}
		|\nabla u(y,t)|^2\frac{t\,d\mu(y)\,dt}{V(y,t)}.
		\label{rfs:upper-area}
	\end{align}
	Since the joining time slice has measure zero,
	\begin{equation}\label{rfs:area-split}
		S_cf(x)^2=S_{c,0}f(x)^2+S_{c,\infty}f(x)^2.
	\end{equation}
	We estimate the first term by integration by parts and the second
	by $\mathcal N_c^\#f$.
	
	\subsubsection{The part below the variable height}
	
	Set
	\begin{equation}\label{rfs:local-maximal}
		N_0f(x)=\sup_{\substack{d(x,y)<2t\\0<t<8c\rho(x)}}|u(y,t)|.
	\end{equation}
	Since $d(x,y)<2t<16c\rho(x)$ in this supremum,
	\begin{equation}\label{rfs:local-inclusion}
		N_0f\leq \mathcal N_c^\#f.
	\end{equation}
	Theorem \ref{thm:restricted-bounds} gives $N_0f\in L^2(M)$.
	
	\begin{lemma}\label{rfs:lower-distribution}
		For every $\alpha>0$,
		\begin{equation}\label{rfs:lower-distribution-estimate}
			\mu\{S_{c,0}f>\alpha\}
			\leq\frac{C}{\alpha^2}
			\int_M\min\{N_0f(x)^2,\alpha^2\}\,d\mu(x)
			=\frac{2C}{\alpha^2}\int_0^\alpha
			s\,\mu\{N_0f>s\}\,ds.
		\end{equation}
	\end{lemma}
	
	\begin{proof}
		For this proof write $N_0=N_0f$, and put
		$$
		E=\{N_0\leq\alpha\},\qquad
		\delta(y)=d(y,E),\qquad
		Q_\alpha=\int_M\min\{N_0^2,\alpha^2\}\,d\mu.
		$$
		We have $\alpha^2\mu(E^c)\leq Q_\alpha\leq\|N_0\|_2^2<\infty$,
		so $E\ne\varnothing$. Moreover,
		$$
		\{N_0>\alpha\}
		=\bigcup_{\substack{y\in M,\ t>0\\|u(y,t)|>\alpha}}
		\left(B(y,2t)\cap\{x:8c\rho(x)>t\}\right)
		$$
		is open. Thus $E$ is closed, $\delta$ is 1-Lipschitz,
		$|\nabla_y\delta|\leq1$ almost everywhere, and $\{\delta=0\}=E$.
		
		Take $\eta$ from \eqref{eq:psi} and a smooth nonincreasing
		$\vartheta:[0,\infty)\to[0,1]$ equal to one on $[0,1]$ and zero
		on $[2,\infty)$. Set
		\begin{equation}\label{rfs:cutoffs}
			\psi(y,t)=\eta(\delta(y)/t),\qquad
			\xi(y,t)=\vartheta\left(\frac{t}{2c\rho(y)}\right),\qquad
			\Psi(y,t)=\psi(y,t)\xi(y,t).
		\end{equation}
		Since $\rho\geq1$ is 1-Lipschitz, the weak chain rule applies
		to these locally Lipschitz functions.
		
		If $x\in E$ and $d(x,y)<t<c\rho(x)$, then $\delta(y)<t$ and
		$$
		\rho(y)\geq\rho(x)-d(x,y)>(1-c)\rho(x),\qquad
		t<\frac{c}{1-c}\rho(y)<2c\rho(y).
		$$
		Hence both cutoffs equal one:
		\begin{equation}\label{rfs:cutoff-one}
			\Psi(y,t)=1
			\quad\text{if some }x\in E\text{ satisfies }
			d(x,y)<t<c\rho(x).
		\end{equation}
		
		If $\psi(y,t)\xi(y,t)>0$, then $\delta(y)<2t$ and
		$t<4c\rho(y)$. Choose $z\in E$ with $d(y,z)<2t$. Then
		$$
		\rho(z)\geq\rho(y)-d(y,z)
		>\rho(y)-2t>(1-8c)\rho(y).
		$$
		Since $4/(1-8c)<8$, we obtain
		$$
		t<4c\rho(y)<8c\rho(z),\qquad d(z,y)<2t.
		$$
		Thus $|u(y,t)|\leq N_0(z)\leq\alpha$ and, by the vertical
		point in \eqref{rfs:local-maximal}, $|u(y,t)|\leq N_0(y)$. Hence
		\begin{equation}\label{rfs:u-cutoff}
			|u(y,t)|\leq\min\{N_0(y),\alpha\}
			\quad\text{whenever }\psi(y,t)>0,\ \xi(y,t)>0.
		\end{equation}
		The same bound holds on both derivative supports:
		$$
		\begin{aligned}
			\{\xi\nabla\psi\ne0\}&\subset\{t<\delta<2t,\ t<4c\rho(y)\},\\
			\{\psi\nabla\xi\ne0\}&\subset\{\delta<2t,\ 2c\rho(y)<t<4c\rho(y)\}.
		\end{aligned}
		$$
		The derivatives vanish at the transition endpoints.
		
		For the distance cutoff, the calculation in
		\eqref{eq:psi-energy}--\eqref{eq:psi-time} gives
		\begin{equation}\label{rfs:distance-derivatives}
			\int_0^\infty t|\nabla\psi(y,t)|^2\,dt
			\leq C\mathbf1_{E^c}(y),\qquad
			\partial_t\psi\geq0,\qquad
			\int_0^\infty\psi\partial_t\psi\,dt
			=\tfrac12\mathbf1_{E^c}(y).
		\end{equation}
		For the height cutoff, write $s=t/(2c\rho(y))$. Its derivatives
		are
		$$
		\partial_t\xi=\frac{\vartheta'(s)}{2c\rho(y)},\qquad
		\nabla_y\xi=-\vartheta'(s)\frac{s\nabla_y\rho(y)}{\rho(y)}.
		$$
		They vanish unless $2c\rho(y)<t<4c\rho(y)$. On this interval,
		$1<s<2$, so
		$$
		|\partial_t\xi|\leq\frac{2\|\vartheta'\|_\infty}{t},\qquad
		|\nabla_y\xi|\leq\frac{8c\|\vartheta'\|_\infty}{t}.
		$$
		It follows that
		\begin{equation}\label{rfs:height-derivatives}
			\int_0^\infty t|\nabla\xi(y,t)|^2\,dt\leq C,
			\qquad \partial_t\xi\leq0,\qquad
			\int_0^\infty\xi|\partial_t\xi|\,dt=\tfrac12.
		\end{equation}
		Here $\xi(y,0)=1$ and $\xi(y,t)=0$ for $t\geq4c\rho(y)$.
		
		The product rule, \eqref{rfs:u-cutoff}, and the two derivative
		estimates imply
		\begin{align}
			D_\alpha
			&:=\int_0^\infty\int_M t|u|^2|\nabla\Psi|^2\,d\mu\,dt
			\nonumber\\
			&\leq2\int_0^\infty\int_M
			t|u|^2\xi^2|\nabla\psi|^2\,d\mu\,dt
			+2\int_0^\infty\int_M
			t|u|^2\psi^2|\nabla\xi|^2\,d\mu\,dt
			\nonumber\\
			&\leq C\alpha^2\mu(E^c)+C Q_\alpha
			\leq C Q_\alpha<\infty.
			\label{rfs:cutoff-error}
		\end{align}
		The first integral uses $y\in E^c$ and $|u|\leq\alpha$;
		the second uses $|u|^2\leq\min\{N_0(y)^2,\alpha^2\}$.
		
		The time derivative has the decomposition
		$$
		\Psi\partial_t\Psi
		=\psi\xi^2\partial_t\psi+\psi^2\xi\partial_t\xi.
		$$
		The first term is nonnegative and the second is nonpositive. By
		\eqref{rfs:u-cutoff}--\eqref{rfs:height-derivatives},
		\begin{align}
			\int_0^\infty\int_M|u|^2|\Psi\partial_t\Psi|\,d\mu\,dt
			&\leq\tfrac12\alpha^2\mu(E^c)+\tfrac12 Q_\alpha<\infty,
			\label{rfs:time-integrability}\\
			\int_0^\infty\int_M|u|^2\Psi\partial_t\Psi\,d\mu\,dt
			&\leq\tfrac12\alpha^2\mu(E^c).
			\label{rfs:time-upper}
		\end{align}
		
		By \eqref{eq:energy},
		$$
		T_\alpha:=\int_0^\infty\int_M t\Psi^2|\nabla u|^2\,d\mu\,dt
		\leq\tfrac12\|f\|_2^2<\infty.
		$$
		Use $h=h_{\varepsilon,T_0}$ and $\zeta_j$ from
		Subsection \ref{subsec:fs-cutoff}. The compactly supported Sobolev test
		$h\zeta_j^2\Psi^2u$ is justified by the smooth approximation
		used for \eqref{eq:finite-test}. The same calculation gives
		\begin{align}
			\int h\zeta_j^2\Psi^2|\nabla u|^2
			={}&-2\operatorname{Re}\int h\zeta_j^2\Psi\overline u\,
			\nabla u\cdot\nabla\Psi
			+\frac12\int h''\zeta_j^2|u|^2\Psi^2\nonumber\\
			&+\int h'\zeta_j^2|u|^2\Psi\partial_t\Psi
			-2\operatorname{Re}\int h\zeta_j\Psi^2\overline u\,
			\nabla_yu\cdot\nabla_y\zeta_j.
			\label{rfs:finite-identity}
		\end{align}
		Here and below, unspecified integrals are over
		$M\times(0,\infty)$ with measure $d\mu\,dt$.
		For fixed $\varepsilon,T_0$, the last term satisfies
		$$
		\left|2\operatorname{Re}\int h\zeta_j\Psi^2\overline u\,
		\nabla_yu\cdot\nabla_y\zeta_j\right|
		\leq\frac{C_{\varepsilon,T_0}}{j}\|f\|_2^2\longrightarrow0,
		$$
		by $L^2$ contractivity and \eqref{eq:energy}. Also
		$|\nabla\Psi|\leq C/\varepsilon$ on $[\varepsilon,2T_0]$,
		so the other absolute integrands are bounded by an integrable
		multiple of
		$$
		\mathbf1_{[\varepsilon,2T_0]}(t)
		\bigl(|\nabla u|^2+|u|\,|\nabla u|+|u|^2\bigr).
		$$
		Letting $j\to\infty$ gives
		\begin{equation}\label{rfs:time-identity}
			\int h\Psi^2|\nabla u|^2
			=-2\operatorname{Re}\int h\Psi\overline u\,
			\nabla u\cdot\nabla\Psi
			+\frac12\int h''|u|^2\Psi^2
			+\int h'|u|^2\Psi\partial_t\Psi.
		\end{equation}
		
		Put $H(t)=\int_M|u(y,t)|^2\Psi(y,t)^2\,d\mu(y)$.
		Since $\xi(y,t)=1$ for $t\leq2c$, we have
		$\Psi(y,t)^2\to\mathbf1_E(y)$ as $t\downarrow0$. Moreover,
		$$
		\left|H(t)-\int_M|f|^2\Psi(y,t)^2\,d\mu\right|
		\leq\|P_tf-f\|_2(\|P_tf\|_2+\|f\|_2)\longrightarrow0.
		$$
		Thus $H(t)\to H_0:=\int_E|f|^2\,d\mu$ at zero, and
		$0\leq H(t)\leq\|P_tf\|_2^2\to0$ at infinity.
		The transition-interval calculation for \eqref{eq:h-second-limit}
		therefore applies to $H$: by \eqref{eq:h-second},
		$$
		\left|\int_0^\infty h''(t)H(t)\,dt-H_0\right|
		\leq C\sup_{\varepsilon\leq t\leq2\varepsilon}|H(t)-H_0|
		+C\|P_{T_0}f\|_2^2\longrightarrow0.
		$$
		For the other terms, use $0\leq h\leq t$, $|h'|\leq C$,
		\eqref{rfs:time-integrability}, and
		$$
		\int t|\nabla u|^2<\infty,\qquad
		\int t\Psi|u|\,|\nabla u|\,|\nabla\Psi|
		\leq T_\alpha^{1/2}D_\alpha^{1/2}<\infty.
		$$
		Dominated convergence in \eqref{rfs:time-identity} now gives
		\begin{align}
			T_\alpha
			={}&\frac12\int_E|f|^2\,d\mu
			-2\operatorname{Re}\int_0^\infty\int_M
			t\Psi\overline u\,\nabla u\cdot\nabla\Psi\,d\mu\,dt
			+\int_0^\infty\int_M
			|u|^2\Psi\partial_t\Psi\,d\mu\,dt.
			\label{rfs:limiting-identity}
		\end{align}
		
		Young's inequality gives
		$$
		2t\Psi|u|\,|\nabla u|\,|\nabla\Psi|
		\leq\tfrac12t\Psi^2|\nabla u|^2
		+2t|u|^2|\nabla\Psi|^2.
		$$
		Since $T_\alpha<\infty$, \eqref{rfs:cutoff-error},
		\eqref{rfs:time-upper}, and \eqref{rfs:limiting-identity} yield
		\begin{equation}\label{rfs:localized-energy}
			T_\alpha\leq\int_E|f|^2\,d\mu+C Q_\alpha.
		\end{equation}
		Choose $t_j\downarrow0$, $t_j<8c$, with $P_{t_j}f\to f$
		almost everywhere, using $L^2$ convergence. The vertical points
		$(x,t_j)$ are allowed in the supremum \eqref{rfs:local-maximal}, so
		$|f|\leq N_0$ almost everywhere. Therefore
		$$
		\int_E|f|^2\,d\mu\leq\int_E N_0^2\,d\mu\leq Q_\alpha,
		\qquad T_\alpha\leq C Q_\alpha.
		$$
		
		Finally, let
		$A(y,t)=\{x\in E:d(x,y)<t<c\rho(x)\}$. This set is contained
		in $B(y,t)$, so $\mu(A(y,t))/V(y,t)\leq1$.
		If it is nonempty, \eqref{rfs:cutoff-one} gives $\Psi(y,t)=1$.
		Tonelli's theorem therefore gives
		$$
		\begin{aligned}
			\int_E S_{c,0}f(x)^2\,d\mu(x)
			&=\int_0^\infty\int_M
			t|\nabla u(y,t)|^2\frac{\mu(A(y,t))}{V(y,t)}\,d\mu(y)\,dt
			\leq T_\alpha\leq C Q_\alpha.
		\end{aligned}
		$$
		Chebyshev's inequality now implies
		$$
		\mu\{S_{c,0}f>\alpha\}
		\leq\mu(E^c)+\alpha^{-2}\int_E (S_{c,0}f)^2\,d\mu
		\leq C\alpha^{-2}Q_\alpha.
		$$
		The equality in \eqref{rfs:lower-distribution-estimate} follows
		from Tonelli and
		$\min\{a^2,\alpha^2\}=2\int_0^\alpha s\mathbf1_{\{a>s\}}\,ds$.
	\end{proof}
	
	\subsubsection{The part above the variable height}
	
	We use the following bounds for geodesic balls in $M$.
	
	\begin{lemma}\label{rfs:tail-volumes}
		Put $b=c\rho(x)$. If $b\geq1$, then $x$ belongs to one of the
		model ends. Write $q=m$ for the large end and $q=n$ for the
		small end. In this case
		\begin{equation}\label{rfs:large-base-volume}
			V(x,b)\leq C b^q,\qquad
			V(y,t)\geq C^{-1}t^q
			\quad(y\in B(x,b),\ t\geq b).
		\end{equation}
		If $b<1$, then
		\begin{equation}\label{rfs:bounded-base-volume}
			\begin{split}
				V(x,b)&\leq C b^m,\\
				V(y,t)&\geq C_c^{-1}t^m\quad(y\in B(x,b),\ b\leq t\leq1),\\
				V(y,t)&\geq C_c^{-1}t^n\quad(y\in B(x,b),\ t\geq1).
			\end{split}
		\end{equation}
	\end{lemma}
	
	\begin{proof}
		If $b\geq1$, then $\rho(x)\geq100$ and
		$b=c\rho(x)<\rho(x)-1=d(x,K)$. Thus $B(x,b)$ stays in the
		model end containing $x$, whose model ball has volume at most
		$Cb^q$. This proves the upper bound in
		\eqref{rfs:large-base-volume}.
		
		For $y\in B(x,b)$ and $t\geq1$, move from $y$ radially
		outwards by $t/2$ to $z$, keeping the sphere coordinate fixed.
		The model ball $B_{\rm mod}(z,t/4)$ stays in the end, since
		its radial coordinates exceed $r(y)+t/4>r_0$. A radial path
		to $z$ followed by a model geodesic gives
		$$
		B_{\rm mod}(z,t/4)\subset B_M(y,3t/4),\qquad
		\mu(B_{\rm mod}(z,t/4))\geq a t^q.
		$$
		For the product end, the latter bound follows by placing a
		Euclidean ball of radius $t/8$ and a sphere ball of radius
		$1/8$ inside the model ball. This proves the lower bound even
		when $B_M(y,t)$ reaches the other end.
		
		If $b<1$, all $y\in B(x,b)$ lie in a fixed compact
		neighbourhood $F_c$ of $K$. Equation \eqref{eq:volume-upper}
		gives $V(x,b)\leq Cb^m$. Finite coordinate charts and
		monotonicity give $V(y,t)\geq a_c t^m$ for $y\in F_c$ and
		$0<t\leq1$.
		
		Fix $y_*$ in the small end and put
		$L_c=\sup_{y\in F_c}d(y,y_*)$. For $t\geq\max(1,4L_c)$,
		move radially outwards from $y_*$ by $t/2$ to $z$. Then
		$$
		B_{\rm mod}(z,t/8)\subset B_M(y,t),\qquad
		L_c+t/2+t/8<t,\qquad
		V(y,t)\geq a t^n.
		$$
		For $1\leq t\leq\max(1,4L_c)$, use
		$V(y,t)\geq V(y,1)\geq a_c$ and decrease the constant.
		This proves \eqref{rfs:bounded-base-volume}.
	\end{proof}
	
	\begin{lemma}\label{rfs:upper-pointwise}
		There is a constant $B\geq1$ such that
		\begin{equation}\label{rfs:upper-pointwise-estimate}
			S_{c,\infty}f(x)\leq B \mathcal N_c^\#f(x),\qquad x\in M.
		\end{equation}
	\end{lemma}
	
	\begin{proof}
		Fix $x$, write $b=c\rho(x)$, and suppose that
		$d(x,y)<b$ and $t\geq b$. Put $r=b/8$ and consider the product
		ball
		$$
		Q=B_{M\times\mathbb R}((y,t),2r).
		$$
		Every $(z,s)\in Q$ satisfies
		$$
		d(x,z)<b+2r=5b/4,\qquad s>t-2r\geq3b/4.
		$$
		In particular, $d(x,z)<2s$ and
		$d(x,z)<16b=16c\rho(x)$. Thus
		$|u(z,s)|\leq \mathcal N_c^\#f(x)$ throughout $Q$.
		Also $r=b/8\leq t/4$, and
		$$
		\rho(y)\geq\rho(x)-d(x,y)>(1-c)\rho(x),\qquad
		r=\frac{c\rho(x)}8\leq\frac{\rho(y)}{16},
		$$
		where the last inequality uses $2c\leq1-c$.
		Lemma \ref{lem:local-harmonic-gradient}, also for complex $u$, gives
		\begin{equation}\label{rfs:tail-gradient}
			|\nabla u(y,t)|\leq\frac{C}{b}\mathcal N_c^\#f(x).
		\end{equation}
		If $b\geq1$, use \eqref{rfs:large-base-volume}, with $q=m$ or
		$q=n$ according to the end containing $x$. Equations
		\eqref{rfs:upper-area} and \eqref{rfs:tail-gradient} give
		$$
		\begin{aligned}
			S_{c,\infty}f(x)^2
			&\leq C b^{-2}\mathcal N_c^\#f(x)^2
			V(x,b)\int_b^\infty t^{1-q}\,dt
			\leq C b^{q-2}\mathcal N_c^\#f(x)^2
			\frac{b^{2-q}}{q-2}
			\leq C \mathcal N_c^\#f(x)^2.
		\end{aligned}
		$$
		Here $q>2$.
		
		If $b<1$, use \eqref{rfs:bounded-base-volume} and split the time
		integral at one. This gives
		$$
		\begin{aligned}
			S_{c,\infty}f(x)^2
			&\leq C_c b^{m-2}\mathcal N_c^\#f(x)^2
			\left(\int_b^1t^{1-m}\,dt+\int_1^\infty t^{1-n}\,dt\right)\\
			&=C_c \mathcal N_c^\#f(x)^2
			\left(\frac{1-b^{m-2}}{m-2}
			+\frac{b^{m-2}}{n-2}\right)
			\leq C_c \mathcal N_c^\#f(x)^2.
		\end{aligned}
		$$
		Since $m>n\geq3$, taking square roots proves the lemma.
	\end{proof}
	
	\subsubsection{Comparison of the two restricted maximal functions}
	
	\begin{lemma}\label{lem:rfs-aperture}
		Let $0<c\leq1/100$. For every $f$ whose Poisson extension is
		defined and smooth for positive times, and every $\lambda>0$,
		\begin{equation}\label{eq:rfs-aperture-distribution}
			\mu\{\mathcal N_c^\#f>\lambda\}
			\leq C\mu\{\mathcal N_cf>\lambda\}.
		\end{equation}
		The constant depends only on the geometry of $M$ and is uniform
		for $0<c\leq1/100$. In particular, for every $0<p<\infty$,
		$$
		\|\mathcal N_c^\#f\|_p\leq C^{1/p}\|\mathcal N_cf\|_p,
		\qquad
		\|\mathcal N_c^\#f\|_{L^{p,\infty}}
		\leq C^{1/p}\|\mathcal N_cf\|_{L^{p,\infty}}.
		$$
		These estimates apply to real or complex $f$.
	\end{lemma}
	
	\begin{proof}
		We first verify the volume estimate
		\begin{equation}\label{eq:rfs-local-doubling}
			V(y,40r)\leq D_0 V(y,r)
			\quad\text{if }r>0\text{ and }40r\leq\rho(y)/5,
		\end{equation}
		where $D_0$ is independent of $y$ and $r$.
		Suppose first that $\rho(y)\geq2$. Then
		$$
		40r\leq\rho(y)/5<\rho(y)-1=d(y,K).
		$$
		Thus both balls agree with the complete model balls. Indeed,
		a path of length at most $40r$ from $y$, in either metric,
		cannot reach the boundary of the end, since its initial segment
		would have length at least $d(y,K)>40r$.
		In the Euclidean model the volume ratio is $40^m$.
		In the product model, ball volumes satisfy
		$$
		W(s)\approx
		\begin{cases}
			s^m,&0<s\leq1,\\
			s^n,&s\geq1.
		\end{cases}
		$$
		These bounds follow by comparison with products of factor balls
		of radii $s/\sqrt2$ and $s$, since sphere-ball volumes are
		comparable to $\min(s^{m-n},1)$.
		If $r\leq1/40$, or if $r\geq1$, the two displayed powers give
		$W(40r)/W(r)\leq C40^m$. If $1/40<r<1$, they give
		$$
		\frac{W(40r)}{W(r)}
		\leq C40^n r^{n-m}\leq C40^m.
		$$
		This proves \eqref{eq:rfs-local-doubling} when $\rho(y)\geq2$.
		
		If $\rho(y)<2$, then $d(y,K)<1$ and $40r<2/5$.
		Finite coordinate charts near $K$ give, uniformly for these centers,
		$$
		a s^m\leq V(y,s)\leq b s^m,\qquad 0<s\leq2/5.
		$$
		The chart bounds first hold for $s\leq s_*$; monotonicity and
		inclusion in a fixed compact set extend them to $s\leq2/5$.
		Applying them at $r$ and $40r$ proves
		\eqref{eq:rfs-local-doubling}.
		
		The set $F=\{\mathcal N_cf>\lambda\}$ is open, since
		$$
		F=\bigcup_{|P_tf(y)|>\lambda}
		\{x:d(x,y)<t,\ d(x,y)<c\rho(x)\}.
		$$
		If $\mu(F)=\infty$, there is nothing to prove. Suppose
		$\mu(F)<\infty$.
		
		Take $x$ such that $\mathcal N_c^\#f(x)>\lambda$. There is $(y,t)$ with
		\begin{equation}\label{eq:rfs-aperture-witness}
			|P_tf(y)|>\lambda,
			\qquad d(x,y)<\min\{2t,16c\rho(x)\}.
		\end{equation}
		The 1-Lipschitz property of $\rho$ gives
		\begin{equation}\label{eq:rfs-aperture-rho}
			\rho(y)\geq\rho(x)-d(x,y)>(1-16c)\rho(x).
		\end{equation}
		Set
		$$
		r=\frac12\min\{t,c\rho(y)\}.
		$$
		If $z\in B(y,r)$, then
		$$
		\rho(z)\geq\rho(y)-d(y,z)
		>\rho(y)-r\geq(1-c/2)\rho(y).
		$$
		Since $c<1$, we have
		$$
		d(z,y)<r\leq t/2<t,
		\qquad
		d(z,y)<r\leq c\rho(y)/2<c\rho(z).
		$$
		Thus $(y,t)\in\Gamma_c(z)$, and
		$\mathcal N_cf(z)\geq|P_tf(y)|>\lambda$. We have proved
		\begin{equation}\label{eq:rfs-aperture-small-ball}
			B(y,r)\subset F.
		\end{equation}
		
		We next check that $x\in B(y,40r)$. If $r=t/2$, then
		\eqref{eq:rfs-aperture-witness} gives $d(x,y)<4r$.
		If $r=c\rho(y)/2$, then \eqref{eq:rfs-aperture-rho} gives
		$$
		d(x,y)<16c\rho(x)
		<\frac{16c}{1-16c}\rho(y)
		=\frac{32}{1-16c}r
		\leq\frac{800}{21}r<40r.
		$$
		Moreover,
		$$
		40r\leq20c\rho(y)\leq\rho(y)/5.
		$$
		It follows from \eqref{eq:rfs-local-doubling} and
		\eqref{eq:rfs-aperture-small-ball} that
		$$
		\frac{1}{V(y,40r)}\int_{B(y,40r)}\mathbf1_F\,d\mu
		\geq\frac{V(y,r)}{V(y,40r)}\geq D_0^{-1}.
		$$
		
		Let $\mathcal M_{\rm HL}$ be the uncentered Hardy--Littlewood
		maximal operator,
		$$
		\mathcal M_{\rm HL}g(x)
		=\sup_{\substack{z\in M,\ s>0\\x\in B(z,s)}}
		\frac{1}{V(z,s)}\int_{B(z,s)}|g|\,d\mu.
		$$
		The preceding ball contains $x$, so we have the inclusion
		$$
		\{\mathcal N_c^\#f>\lambda\}
		\subset
		\{\mathcal M_{\rm HL}\mathbf1_F>(2D_0)^{-1}\}.
		$$
		The weak type $(1,1)$ estimate in \cite[Theorem 2]{DLS} now gives
		$$
		\mu\{\mathcal N_c^\#f>\lambda\}
		\leq 2C_{\rm HL}D_0\|\mathbf1_F\|_1
		=2C_{\rm HL}D_0\mu\{\mathcal N_cf>\lambda\}.
		$$
		This is \eqref{eq:rfs-aperture-distribution}. Integration against
		$p\lambda^{p-1}d\lambda$ gives the strong comparison; taking
		the supremum of $\lambda\mu\{\mathcal N_c^\#f>\lambda\}^{1/p}$
		gives the weak comparison. The proof uses $|P_tf(y)|$ and
		applies to complex $f$.
	\end{proof}

	\subsubsection{The distribution estimate and the norm comparison}
	
	By \eqref{rfs:area-split},
	$$
	\{S_cf>\lambda\}\subset
	\{S_{c,0}f>\lambda/\sqrt2\}\cup
	\{S_{c,\infty}f>\lambda/\sqrt2\}.
	$$
	Lemmas \ref{rfs:lower-distribution} and \ref{rfs:upper-pointwise},
	together with \eqref{rfs:local-inclusion}, give
	\begin{align}
		\mu\{S_cf>\lambda\}
		&\leq\mu\{\mathcal N_c^\#f>\lambda/(\sqrt2 B)\}
		+\frac{C}{\lambda^2}\int_0^{\lambda/\sqrt2}
		s\,\mu\{N_0f>s\}\,ds\nonumber\\
		&\leq\mu\{\mathcal N_c^\#f>\lambda/A\}
		+\frac{C}{\lambda^2}\int_0^\lambda
		s\,\mu\{\mathcal N_c^\#f>s\}\,ds,
		\label{rfs:enlarged-distribution}
	\end{align}
	where $A=\sqrt2 B\geq1$. Lemma \ref{lem:rfs-aperture} gives
	\begin{equation}\label{rfs:distrib-before-absorption}
		\mu\{S_cf>\lambda\}
		\leq C\mu\{\mathcal N_cf>\lambda/A\}
		+\frac{C}{\lambda^2}\int_0^\lambda
		s\,\mu\{\mathcal N_cf>s\}\,ds.
	\end{equation}
	To write the first distribution term at level $\lambda$, split
	$$
	\mu\{\mathcal N_cf>\lambda/A\}
	=\mu\{\mathcal N_cf>\lambda\}
	+\mu\{\lambda/A<\mathcal N_cf\leq\lambda\}.
	$$
	For $0<s<\lambda/A$, the set in the second term is contained
	in $\{\mathcal N_cf>s\}$. Hence
	$$
	\int_0^\lambda s\,\mu\{\mathcal N_cf>s\}\,ds
	\geq\int_0^{\lambda/A}s\,\mu\{\mathcal N_cf>s\}\,ds
	\geq\frac{\lambda^2}{2A^2}
	\mu\{\lambda/A<\mathcal N_cf\leq\lambda\}.
	$$
	Substituting these two relations into
	\eqref{rfs:distrib-before-absorption} gives
	\begin{equation}\label{rfs:final-distribution}
		\mu\{S_cf>\lambda\}
		\leq C_c\mu\{\mathcal N_cf>\lambda\}
		+\frac{C_c}{\lambda^2}\int_0^\lambda
		s\,\mu\{\mathcal N_cf>s\}\,ds.
	\end{equation}
	This proves \eqref{eq:restricted-fs-intro} with the same
	restricted cone in both operators.
	
	For $0<p<2$, assume that $\mathcal N_cf\in L^p(M)$, multiply
	\eqref{rfs:final-distribution} by $p\lambda^{p-1}$, and
	integrate over $\lambda>0$. The first term gives
	$$
	p\int_0^\infty\lambda^{p-1}
	\mu\{\mathcal N_cf>\lambda\}\,d\lambda
	=\|\mathcal N_cf\|_p^p.
	$$
	For the second term, Tonelli and $p<2$ give
	$$
	\begin{aligned}
		p\int_0^\infty\lambda^{p-3}
		\int_0^\lambda s\,\mu\{\mathcal N_cf>s\}\,ds\,d\lambda
		&=p\int_0^\infty s\,\mu\{\mathcal N_cf>s\}
		\int_s^\infty\lambda^{p-3}\,d\lambda\,ds\\
		&=\frac{p}{2-p}\int_0^\infty
		s^{p-1}\mu\{\mathcal N_cf>s\}\,ds
		=\frac{1}{2-p}\|\mathcal N_cf\|_p^p.
	\end{aligned}
	$$
	Tonelli applies to nonnegative extended integrals, so this also
	proves $S_cf\in L^p$ and
	\begin{equation}\label{rfs:norm-comparison}
		\|S_cf\|_p^p
		\leq C_c\left(1+\frac1{2-p}\right)\|\mathcal N_cf\|_p^p.
	\end{equation}
	This proves Theorem \ref{thm:restricted-fs}. For $1<p<2$,
	Theorem \ref{thm:restricted-bounds} also gives
	$\|S_cf\|_p\leq C_{p,c}\|f\|_p$. At $p=2$, use the energy
	identity \eqref{eq:energy}.

	\section{The two-sided \texorpdfstring{$L^2$}{L2} estimate}
	\label{sec:L2-statement}

	\subsection{The theorem and the geometric estimates}

	Let $M$ be complete and connected. Outside a smooth compact set
	$K$, its two ends are isometric to the complements of bounded sets
	in $\mathbb R^m$ and $\mathbb R^n\times\mathbb S^{m-n}$.
	The sphere has its usual metric. Write
	$$
	\Delta=-\operatorname{div}\nabla,\qquad P_t=e^{-t\sqrt \Delta},\qquad
	\rho(x)=1+d(x,K),\qquad V(x,r)=\mu(B(x,r)).
	$$
	Here $\Delta$ is the nonnegative self-adjoint Laplace--Beltrami operator
	on $L^2(M)$. For $f\in L^2(M)$, put $u(x,t)=P_tf(x)$ and
	\begin{equation}\label{eq:def}
		\Sc f(x)^2
		=\int_0^\infty\int_{B(x,\min\{t,c\rho(x)\})}
		|\nabla_{y,t}u(y,t)|^2\frac{t\dd\mu(y)\dd t}{V(y,t)},
		\qquad 0<c\leq\frac1{100}.
	\end{equation}
	Thus $|\nabla_{y,t}u|^2=|\nabla_yu|^2+|\partial_tu|^2$.
	The denominator in \eqref{eq:def} is the volume of the full ball
	$B(y,t)$.
	
	\begin{theorem}\label{thm:main}
		There is $C=C(M,c)<\infty$ such that, for every $f\in L^2(M)$,
		\begin{equation}\label{eq:main}
			C^{-1}\|f\|_2\leq\|\Sc f\|_2\leq 2^{-1/2}\|f\|_2.
		\end{equation}
		The lower estimate holds already for a part with
		$0<t<\varepsilon R(y)$, where $R$ is a smooth function comparable
		to $\rho$. Here $\varepsilon=10\varepsilon_0$, with
		$\varepsilon_0=\varepsilon_0(M,c)>0$ fixed in \eqref{eq:height}.
	\end{theorem}
	
	All constants may depend on $M$ and $c$. We shall use the
	following geometric facts. There is a smooth proper function $R$
	such that
	\begin{equation}\label{eq:radius}
		1\leq R,\qquad C_R^{-1}\rho\leq R\leq C_R\rho,
		\qquad |\nabla R|\leq C_R.
	\end{equation}
	Such a function is obtained by joining the radial coordinate on
	each end to a positive smooth function on a compact set and then
	multiplying by a constant.
	
	For sufficiently small $a_0>0$, balls with $r\leq a_0\rho(x)$
	satisfy uniform volume doubling and a uniform Poincar\'e
	inequality. Far from $K$ they lie in one of the two product
	manifolds. The remaining balls have bounded radii and centers in
	a fixed compact set. On these balls, the local estimates for
	harmonic functions and their gradients have uniform constants.
	
	We also use the diagonal Poisson estimate
	\begin{equation}\label{eq:poisson-diagonal}
		P_s(x,x)\leq C s^{-n}\quad(s\geq1),\qquad
		\sup_x P_s(x,x)<\infty\quad(s>0\text{ fixed}).
	\end{equation}
	It follows from the heat kernel bounds in \cite{GS} by
	subordination; see also \cite[Theorem 2.2]{BDLW}.
	In particular, symmetry and the semigroup identity give
	\begin{equation}\label{eq:decay}
		|P_tf(x)|\leq P_{2t}(x,x)^{1/2}\|f\|_2
		\leq Ct^{-n/2}\|f\|_2\qquad(t\geq1).
	\end{equation}
	
	\subsection{The upper estimate and the part at small heights}
	
	\begin{lemma}\label{lem:energy}
		For $f\in L^2(M)$,
		\begin{equation}\label{eq:L2-energy}
			\int_0^\infty t\|\nabla_{x,t}P_tf\|_2^2\dd t
			=\frac12\|f\|_2^2.
		\end{equation}
		If $A=\Sc f$, then
		\begin{equation}\label{eq:lower-part}
			\int_M\int_0^{c\rho(y)/(1+c)}
			t|\nabla_{y,t}u(y,t)|^2\dd t\dd\mu(y)
			\leq\|A\|_2^2\leq\frac12\|f\|_2^2.
		\end{equation}
	\end{lemma}
	
	\begin{proof}
		The kernel of $\Delta$ in $L^2(M)$ is zero. Indeed,
		$\Delta h=0$ implies $\|\nabla h\|_2^2=0$; connectedness and
		$\mu(M)=\infty$ then give $h=0$. The spectral theorem and
		the quadratic form of $\Delta$ give
		$$
		\|\nabla_x P_tf\|_2^2
		=\|\sqrt \Delta P_tf\|_2^2
		=\|\partial_tP_tf\|_2^2,
		\qquad
		\int_0^\infty t\lambda e^{-2t\sqrt\lambda}\dd t=\frac14
		\quad(\lambda>0).
		$$
		This proves \eqref{eq:L2-energy}.
		
		Set
		$$
		W_c(y,t)=\{x:d(x,y)<t,\ d(x,y)<c\rho(x)\},\qquad
		\omega_c(y,t)=\frac{\mu(W_c(y,t))}{V(y,t)}.
		$$
		Tonelli's theorem gives the exact identity
		\begin{equation}\label{eq:shadow}
			\|A\|_2^2
			=\int_M\int_0^\infty
			t|\nabla_{y,t}u(y,t)|^2\omega_c(y,t)\dd t\dd\mu(y).
		\end{equation}
		Since $W_c(y,t)\subset B(y,t)$, we have $0\leq\omega_c\leq1$.
		If $t\leq c\rho(y)/(1+c)$ and $d(x,y)<t$, then
		$$
		c\rho(x)\geq c\rho(y)-c\,d(x,y)
		>c\rho(y)-ct\geq t>d(x,y).
		$$
		Hence $W_c(y,t)=B(y,t)$ and $\omega_c(y,t)=1$ in this range.
		Equations \eqref{eq:L2-energy} and \eqref{eq:shadow} prove
		\eqref{eq:lower-part}.
	\end{proof}
	
	Fix $\varepsilon_0>0$ sufficiently small that
	\begin{equation}\label{eq:height}
		H(x):=\varepsilon_0R(x),\quad C_R\varepsilon_0\leq1,\quad
		10H(x)\leq\min\{a_0,c/(1+c)\}\rho(x).
	\end{equation}
	Then $H\geq\varepsilon_0$ and $|\nabla H|\leq C_R\varepsilon_0$.
	
	\begin{lemma}\label{lem:pointwise}
		For almost every $x\in M$ and every $0<t\leq2H(x)$,
		\begin{equation}\label{eq:gradient}
			t|\nabla_{x,t}u(x,t)|\leq C A(x).
		\end{equation}
		For $g(x)=u(x,H(x))$, we have
		\begin{equation}\label{eq:trace-pointwise}
			\rho(x)|\nabla g(x)|\leq C A(x),\qquad
			|u(x,t)-g(x)|\leq C A(x)|\log(t/H(x))|.
		\end{equation}
		The second inequality holds for $0<t\leq2H(x)$.
	\end{lemma}
	
	\begin{proof}
		The function $u$ is harmonic on $M\times(0,\infty)$ for
		$\D=\operatorname{div}_x\nabla_x+\partial_t^2$.
		Fix $0<t\leq2H(x)$ and put
		$q(y,s)=|\nabla_{y,s}u(y,s)|^2$ and $r=t/8$.
		By \eqref{eq:height}, $t\leq a_0\rho(x)/5$.
		If $d(x,K)\geq1$, then
		$$
		5r=\frac{5t}{8}\leq\frac{a_0}{8}(1+d(x,K))
		\leq\frac{a_0}{4}d(x,K).
		$$
		Thus $B(x,5r)$ lies in one product end, where
		$\operatorname{Ric}_M\geq0$. In this case put $\kappa=0$.
		If $d(x,K)<1$, then $t<2a_0/5$, and $B(x,5r)$ lies
		in a fixed compact neighbourhood of $K$. In this case choose a
		fixed $\kappa\geq0$ such that $\operatorname{Ric}_M\geq-\kappa$
		there. In both cases $\kappa r^2\leq C$.
		
		Bochner's formula on $M\times(0,\infty)$ gives
		$$
		\D q=2|\nabla^2_{y,s}u|^2
		+2\operatorname{Ric}_M(\nabla_yu,\nabla_yu)
		\geq-2\kappa q
		$$
		on $B_{M\times\mathbb R}((x,t),5r)$.
		On $M\times\mathbb R^2$ set
		$$
		Q(y,s,v)=q(y,s)\cosh(\sqrt{2\kappa}\,v).
		$$
		Since $5r<t$, this is defined on
		$B_{M\times\mathbb R^2}((x,t,0),5r)$, and
		$$
		(\D+\partial_v^2)Q
		=\cosh(\sqrt{2\kappa}\,v)(\D q+2\kappa q)\geq0.
		$$
		Apply the local $L^1$ mean value inequality for nonnegative
		subharmonic functions \cite[Theorem 2.1, $p=1$]{LiSchoen}
		on this ball. Its constant is uniform since $\kappa r^2\leq C$.
		With $B_r=B_{M\times\mathbb R^2}((x,t,0),r)$, it gives
		$$
		q(x,t)\leq\frac{C}{|B_r|}\int_{B_r}Q(y,s,v)
		\dd\mu(y)\dd s\dd v.
		$$
		The product metric and the uniform volume doubling for
		$t\leq a_0\rho(x)/5$ give
		$$
		\begin{aligned}
			|B_r|&\geq r^2V(x,r/2)\geq c\,t^2V(x,t),\\
			\int_{B_r}Q\dd\mu\dd s\dd v
			&\leq C r\int_{t-r}^{t+r}\int_{B(x,r)}q(y,s)
			\dd\mu(y)\dd s.
		\end{aligned}
		$$
		Here $\cosh(\sqrt{2\kappa}\,v)\leq C$ for $|v|<r$.
		Since $r=t/8$, we obtain, for $0<t\leq2H(x)$,
		\begin{equation}\label{eq:meanvalue}
			|\nabla u(x,t)|^2
			\leq\frac{C}{tV(x,t)}
			\int_{t/2}^{3t/2}\int_{B(x,t/4)}|\nabla u(y,s)|^2
			\dd\mu(y)\dd s.
		\end{equation}
		
		For $y\in B(x,t/4)$ and $t/2<s<3t/2$, \eqref{eq:height} gives
		$$
		d(x,y)<\min\{s,c\rho(x)\},\qquad
		V(y,s)\leq C V(x,t).
		$$
		Since $s\geq t/2$, equations \eqref{eq:def} and
		\eqref{eq:meanvalue} imply
		$$
		A(x)^2\geq\frac{t}{C V(x,t)}
		\int_{t/2}^{3t/2}\int_{B(x,t/4)}|\nabla u(y,s)|^2
		\dd\mu(y)\dd s
		\geq C^{-1}t^2|\nabla u(x,t)|^2.
		$$
		The chain rule now gives
		$$
		|\nabla g|
		\leq |\nabla_xu(x,H)|+|\partial_tu(x,H)||\nabla H|
		\leq C H^{-1}A.
		$$
		Since $\rho/H\leq C/\varepsilon_0$, the first estimate in
		\eqref{eq:trace-pointwise} follows. Integrating
		$|\partial_su(x,s)|\leq CA(x)/s$ from $H(x)$ to $t$ gives the
		second estimate.
	\end{proof}
	
	\subsection{The trace on the two ends}
	
	\begin{lemma}\label{lem:hardy}
		Suppose that $g$ is smooth, tends to zero on both ends, and
		$\rho\nabla g\in L^2(M)$. Then
		\begin{equation}\label{eq:hardy}
			\|g\|_2\leq C\|\rho\nabla g\|_2.
		\end{equation}
		Consequently the trace in Lemma \ref{lem:pointwise} satisfies
		\begin{equation}\label{eq:trace-L2}
			\|g\|_2+\|\rho\nabla g\|_2\leq C\|A\|_2.
		\end{equation}
	\end{lemma}
	
	\begin{proof}
		Enlarge a compact set so that each remaining end has coordinates
		$(r,\omega)$, $r>r_0$, and measure $r^{q-1}\dd r\dd\omega$,
		where $q=m$ or $q=n$. In the second case $\omega$ also includes
		the compact sphere variable. On each end $r\simeq\rho$.
		
		For almost every $\omega$, Cauchy--Schwarz gives
		$$
		\int_r^\infty|\partial_sg(s,\omega)|\dd s
		\leq
		\left(\int_r^\infty|\partial_sg|^2s^{q+1}\dd s\right)^{1/2}
		\left(\int_r^\infty s^{-q-1}\dd s\right)^{1/2}<\infty.
		$$
		The limit of $g$ at infinity is zero. Hence
		$$
		g(r,\omega)=-\int_1^\infty sr\,\partial_rg(sr,\omega)
		\frac{\dd s}{s}.
		$$
		For $E_q=\{r>r_0\}$, a change of variable gives
		$$
		\|sr\,\partial_rg(sr,\omega)\|_{L^2(E_q)}
		\leq s^{-q/2}\|r\partial_rg\|_{L^2(E_q)}.
		$$
		Minkowski's inequality, first with a finite upper limit and then
		with Fatou's lemma, yields
		\begin{equation}\label{eq:radial}
			\|g\|_{L^2(E_q)}
			\leq\int_1^\infty s^{-1-q/2}\dd s\,
			\|r\partial_rg\|_{L^2(E_q)}
			=\frac2q\|r\partial_rg\|_{L^2(E_q)}.
		\end{equation}
		
		Choose a fixed bounded connected smooth set $U$ containing the
		remaining compact part and a fixed annulus $D\subset U\cap E_m$
		of positive measure. Put $g_U=\mu(U)^{-1}\int_Ug\dd\mu$.
		Poincar\'e's inequality and Cauchy--Schwarz give
		$$
		\|g-g_U\|_{L^2(U)}\leq C\|\nabla g\|_{L^2(U)},\qquad
		|g_U|\leq\mu(D)^{-1/2}
		\bigl(\|g\|_{L^2(D)}+\|g-g_U\|_{L^2(U)}\bigr).
		$$
		Consequently,
		$$
		\|g\|_{L^2(U)}
		\leq C\bigl(\|\nabla g\|_{L^2(U)}+\|g\|_{L^2(D)}\bigr).
		$$
		Together with \eqref{eq:radial}, this proves \eqref{eq:hardy}.
		
		For $g(x)=P_{H(x)}f(x)$, \eqref{eq:decay} and the properness of
		$H$ give $g(x)\to0$ on both ends. The function $g$ is smooth,
		since $u$ and $H$ are smooth and $H\geq\varepsilon_0>0$.
		Lemma \ref{lem:pointwise} and \eqref{eq:lower-part} imply
		$\rho\nabla g\in L^2$. Applying \eqref{eq:hardy} proves
		\eqref{eq:trace-L2}.
	\end{proof}
	
	\subsection{Green's formula and the lower estimate}
	
	Choose $\theta\in C^\infty([0,\infty))$ with
	$$
	0\leq\theta\leq1,\qquad
	\theta=1\text{ on }[0,1],\qquad
	\theta=0\text{ on }[2,\infty),
	\qquad w(x,t)=t\theta(t/H(x)).
	$$
	With $s=t/H(x)$,
	\begin{equation}\label{eq:w-derivatives}
		w_{tt}=H^{-1}\{2\theta'(s)+s\theta''(s)\},\qquad
		\nabla_xw=-s^2\theta'(s)\nabla H.
	\end{equation}
	Both derivatives vanish unless $H<t<2H$. In particular,
	\begin{equation}\label{eq:w-bounds}
		0\leq w\leq t\mathbf1_{\{t<2H\}},\qquad
		|w_{tt}|\leq C H^{-1}\mathbf1_{\{H<t<2H\}},\qquad
		|\nabla_xw|\leq C\varepsilon_0\mathbf1_{\{H<t<2H\}}.
	\end{equation}
	
	\begin{lemma}\label{lem:green}
		For every $f\in L^2(M)$,
		\begin{equation}\label{eq:green}
			\|f\|_2^2
			=2\int_M\int_0^\infty w|\nabla_{x,t}u|^2\dd t\dd\mu
			-\int_M\int_0^\infty w_{tt}|u|^2\dd t\dd\mu
			+\int_M\int_0^\infty
			\nabla_xw\cdot\nabla_x|u|^2\dd t\dd\mu.
		\end{equation}
		All three integrals are absolutely convergent.
	\end{lemma}
	
	\begin{proof}
		For complex $u$, all gradient products below are understood in
		their real form. Harmonicity gives
		\begin{equation}\label{eq:lap-square}
			\D|u|^2=2|\nabla_{x,t}u|^2.
		\end{equation}
		We first insert cutoffs in $x$ and $t$.
		Let $\chi\in C^\infty([0,\infty))$ equal one on $[0,1]$ and
		zero on $[2,\infty)$, with $0\leq\chi\leq1$, and set
		$\chi_N(x)=\chi(R(x)/N)$. Then
		$$
		|\nabla\chi_N|\leq C/N,\qquad
		\supp\nabla\chi_N\subset F_N:=\{N\leq R\leq2N\}.
		$$
		Let $\zeta\in C^\infty([0,\infty))$ equal zero on $[0,1]$ and
		one on $[2,\infty)$, with $0\leq\zeta\leq1$. Put
		$$
		h_\tau(t)=t\zeta(t/\tau),\qquad
		w_\tau(x,t)=h_\tau(t)\theta(t/H(x)),\qquad
		0<\tau<\varepsilon_0/4.
		$$
		The function $\chi_Nw_\tau$ has compact support in
		$M\times(0,\infty)$. Integrating \eqref{eq:lap-square} twice in
		$t$ and once in $x$ gives
		\begin{equation}\label{eq:cutoff-green}
			2\iint\chi_Nw_\tau|\nabla u|^2
			=\iint\chi_N(w_\tau)_{tt}|u|^2
			-\iint\chi_N\nabla_xw_\tau\cdot\nabla_x|u|^2
			-\iint w_\tau\nabla\chi_N\cdot\nabla_x|u|^2.
		\end{equation}
		Here and below $\iint$ denotes integration over
		$M\times(0,\infty)$ with measure $\dd\mu\dd t$.
		
		The two transition intervals are disjoint. Thus
		$$
		(w_\tau)_{tt}=h_\tau''+w_{tt},\qquad
		\supp h_\tau''\subset[\tau,2\tau],\qquad
		\int_0^\infty h_\tau''\dd t=1,\qquad
		\int_0^\infty|h_\tau''|\dd t\leq C.
		$$
		Since $P_tf\to f$ in $L^2$,
		$$
		\lim_{\tau\downarrow0}
		\int_0^\infty h_\tau''(t)
		\int_M\chi_N|u(x,t)|^2\dd\mu\dd t
		=\int_M\chi_N|f|^2\dd\mu.
		$$
		Moreover, $\nabla_xw_\tau=\nabla_xw$, and the first integral in
		\eqref{eq:cutoff-green} converges by \eqref{eq:L2-energy}.
		For the last integral use Lemma \ref{lem:pointwise}:
		\begin{equation}\label{eq:cutoff-error}
			\begin{split}
				\iint|w\nabla\chi_N\cdot\nabla_x|u|^2|
				&\leq\frac C N\int_{F_N}A(x)
				\int_0^{2H(x)}\bigl(|g(x)|+A(x)|\log(t/H(x))|\bigr)
				\dd t\dd\mu(x)\\
				&\leq C\int_{F_N}(A^2+|g|^2)\dd\mu.
			\end{split}
		\end{equation}
		Here $H/N\leq2\varepsilon_0$ on $F_N$, and
		$\int_0^2|\log s|\dd s<\infty$. The same bound applies to
		$w_\tau$, since $0\leq w_\tau\leq w$.
		Dominated convergence therefore applies as $\tau\downarrow0$.
		After rearranging \eqref{eq:cutoff-green}, we obtain
		\begin{equation}\label{eq:spatial-cutoff}
			\int_M\chi_N|f|^2\dd\mu
			=2\iint\chi_Nw|\nabla u|^2
			-\iint\chi_Nw_{tt}|u|^2
			+\iint\chi_N\nabla_xw\cdot\nabla_x|u|^2
			+\iint w\nabla\chi_N\cdot\nabla_x|u|^2.
		\end{equation}
		
		For $H<t<2H$, Lemma \ref{lem:pointwise} gives
		\begin{equation}\label{eq:band}
			|u(x,t)|\leq|g(x)|+CA(x),\qquad
			|\nabla u(x,t)|\leq CA(x)/t.
		\end{equation}
		Equations \eqref{eq:w-bounds} and \eqref{eq:band} imply
		\begin{align}
			\iint|w_{tt}|\,|u|^2
			&\leq C\int_M(|g|^2+A^2)\dd\mu,\label{eq:error-t}\\
			\iint|\nabla_xw|\,|\nabla_x|u|^2|
			&\leq C\int_M A(|g|+A)\dd\mu
			\leq C(\|g\|_2^2+\|A\|_2^2).\label{eq:error-x}
		\end{align}
		These are finite by Lemma \ref{lem:hardy}.
		The first integral in \eqref{eq:green} is finite by
		\eqref{eq:L2-energy}. Finally, the right side of
		\eqref{eq:cutoff-error} tends to zero as $N\to\infty$.
		Dominated convergence in \eqref{eq:spatial-cutoff} proves
		\eqref{eq:green}.
	\end{proof}
	
	\begin{proof}[Proof of Theorem \ref{thm:main}]
		By \eqref{eq:height}, \eqref{eq:lower-part}, and
		\eqref{eq:w-bounds},
		$$
		\iint w|\nabla u|^2
		\leq\int_M\int_0^{2H(x)}t|\nabla u(x,t)|^2\dd t\dd\mu(x)
		\leq\|A\|_2^2.
		$$
		Lemmas \ref{lem:hardy} and \ref{lem:green}, together with
		\eqref{eq:error-t}--\eqref{eq:error-x}, give
		$$
		\|f\|_2^2
		\leq2\|A\|_2^2+C(\|g\|_2^2+\|A\|_2^2)
		\leq C\|A\|_2^2.
		$$
		The upper estimate is Lemma \ref{lem:energy}.
		
		To check the last assertion of the theorem, replace $A$ throughout
		by the function obtained from \eqref{eq:def} by also requiring
		$t<10H(y)$. If $t\leq2H(x)$, $y\in B(x,t/4)$, and
		$t/2<s<3t/2$, then \eqref{eq:height} gives
		$$
		H(y)\geq H(x)-C_R\varepsilon_0t/4\geq H(x)/2,
		\qquad s\leq3H(x)\leq6H(y)<10H(y).
		$$
		Thus the cylinders in Lemma \ref{lem:pointwise} are retained.
		The energy over $0<t<2H(y)$ is also retained. All estimates used
		in the lower bound therefore hold for this smaller area function.
		Since $10H(y)=\varepsilon R(y)$ with
		$\varepsilon=10\varepsilon_0$, this proves the last assertion.
	\end{proof}
	
	\begin{corollary}\label{cor:H2}
		The completion of $L^2(M)$ in the norm $\|\Sc f\|_2$ is
		$L^2(M)$, with equivalent norms. In particular,
		$$
		H^2_{S_c}(M)=L^2(M).
		$$
	\end{corollary}

	\section{The Hardy space: atoms and maximal functions}
	\label{sec:interp-statement}
	
	We prove the atomic and maximal assertions in Theorem
	\ref{thm:hardy-characterization}.
	
	\subsection{The atoms and the area comparison}
	
	We use the global space in \eqref{eq:H1}. An atom has the form
	$a=\Delta b$, with
	\begin{equation}\label{eq:atoms}
		\begin{gathered}
			b\in\operatorname{Dom}(\Delta),\quad
			\supp b\cup\supp \Delta b\subset B=B(x,r),\quad r\leq\delta\rho(x),\\
			\|b\|_2\leq r^2\mu(B)^{-1/2},\qquad
			\|\Delta b\|_2\leq\mu(B)^{-1/2}.
		\end{gathered}
	\end{equation}
	The number $\delta>0$ is fixed and sufficiently small.
	Let $H^1_{\mathrm{at},\Delta,\delta}(M)$ be the space of $L^1$ sums
	$\sum_\nu\lambda_\nu a_\nu$ of these atoms with
	$\sum_\nu|\lambda_\nu|<\infty$. Its norm is the infimum of this
	coefficient sum. Cauchy--Schwarz gives $\|a\|_1\leq1$.
	
	Every atom cancels against bounded global harmonic functions.
	If $\Delta h=0$ and $h$ is bounded, choose
	$\chi\in C_c^\infty(M)$ equal to one near $\overline B$. Then
	\begin{equation}\label{eq:atoms-harmonic-cancellation}
		\int_M a\overline h\dd\mu
		=\langle \Delta b,\chi h\rangle
		=\langle b,\Delta(\chi h)\rangle=0.
	\end{equation}
	Here $\chi h\in\operatorname{Dom}(\Delta)$, and $\Delta(\chi h)=0$ near
	$\supp b$. The same cancellation holds for every $L^1$ atomic
	sum by absolute convergence.
	
	Put
	$$
	\mathcal N_c^\#h(x)
	=\sup_{d(x,y)<\min\{2t,16c\rho(x)\}}|P_th(y)|.
	$$
	For $1<q<2$ and $h\in C_c^\infty(M)$, the restricted maximal
	estimate and Fefferman--Stein inequality give
	\begin{equation}\label{eq:prior-N}
		\|\mathcal N_c^\#h\|_q\leq C_q\|h\|_q,
		\qquad \|\Sc h\|_1\leq C\|\mathcal N_c^\#h\|_1.
	\end{equation}
	These follow from Sections \ref{sec:restricted-bounds} and
	\ref{sec:rfs}.
	We first extend the area comparison to $L^2$ functions.
	
	\begin{lemma}\label{lem:Hardy-FS-L2}
		If $f\in L^2(M)$ and $\mathcal N_cf\in L^1(M)$, then
		\begin{equation}\label{eq:Hardy-FS-L2}
			\|\Sc f\|_1\leq C_c\|\mathcal N_cf\|_1.
		\end{equation}
	\end{lemma}
	\begin{proof}
		Choose $f_j\in C_c^\infty(M)$ with $f_j\to f$ in $L^2$.
		The $L^2$ bounds in Theorem \ref{thm:restricted-bounds} give
		$$
		\|\Sc f_j-\Sc f\|_2
		+\|\mathcal N_cf_j-\mathcal N_cf\|_2
		\leq C_c\|f_j-f\|_2\longrightarrow0.
		$$
		For the first term use Minkowski on each cone; for the second,
		use the inequality between the difference of two suprema and
		the supremum of their difference.
		After passage to a subsequence, both functions converge almost
		everywhere. Their level-set measures converge at every common
		continuity point $\lambda>0$ of the two limiting distributions.
		Indeed, convergence in $L^2$ gives convergence in measure, and
		for $0<\eta<\lambda$ the level sets are bounded between those
		at $\lambda-\eta$ and $\lambda+\eta$, apart from a set whose
		measure tends to zero.
		Also,
		$$
		\int_0^\lambda s\,\mu\{\mathcal N_cf_j>s\}\,ds
		=\tfrac12\|\min\{\mathcal N_cf_j,\lambda\}\|_2^2
		\longrightarrow
		\tfrac12\|\min\{\mathcal N_cf,\lambda\}\|_2^2.
		$$
		Here truncation is a contraction in $L^2$.
		Pass to the limit in \eqref{eq:restricted-fs-intro}.
		The resulting distribution inequality holds for almost every
		$\lambda>0$. Its integral term satisfies, by Tonelli,
		$$
		\int_0^\infty\lambda^{-2}
		\int_0^\lambda s\,\mu\{\mathcal N_cf>s\}\,ds\,d\lambda
		=\int_0^\infty\mu\{\mathcal N_cf>s\}\,ds
		=\|\mathcal N_cf\|_1.
		$$
		Integration in $\lambda$ proves \eqref{eq:Hardy-FS-L2}.
	\end{proof}
	
	\begin{lemma}\label{lem:Hardy-kernel}
		For every ball $B=B(x_B,r)$ with $r\leq\delta\rho(x_B)$,
		\begin{equation}\label{eq:prior-kernel}
			\sup_{z\in B}\int_{M\setminus2B}
			\sup_{d(x,y)<\min\{2t,16c\rho(x)\}}
			|(\Delta P_t)(y,z)|\dd\mu(x)\leq Cr^{-2}.
		\end{equation}
	\end{lemma}
	
	\begin{proof}
		In \cite[Theorem~2.2]{BDLW}, $P_{t,k}$ is the kernel of
		$(t\sqrt \Delta)^k e^{-t\sqrt \Delta}$. In particular,
		$$
		P_{t,2}=t^2\Delta P_t,\qquad
		t^{-2}t^{-d}\left(\frac{t}{t+d(y,z)}\right)^{d+2}
		=(t+d(y,z))^{-d-2}.
		$$
		Thus $k=2$ gives a bound uniform as $t\downarrow0$ away
		from the diagonal. The terms of dimension $m$ are bounded by
		$C(t+d(y,z))^{-m-2}$; the extra cross-end coefficient is at
		most a constant. The remaining terms are
		$$
		\begin{cases}
			C\rho(y)^{2-m}\rho(z)^{2-m}
			(t+\rho(y)+\rho(z))^{-n-2},
			&y,z\in E_m,\\
			C\rho(y)^{2-m}(t+d(y,z))^{-n-2},
			&y\in E_m,\ z\in E_n\cup K,\\
			C\rho(z)^{2-m}(t+d(y,z))^{-n-2},
			&z\in E_m,\ y\in E_n\cup K,\\
			C(t+d(y,z))^{-n-2},
			&y,z\in E_n\cup K.
		\end{cases}
		$$
		The radii in that theorem are comparable to $\rho$.
		Every coefficient displayed here is at most a constant.
		Also $d(y,z)\leq C(\rho(y)+\rho(z))$, by joining the two
		points through the fixed compact part. Consequently, in all
		cases,
		$$
		|(\Delta P_t)(y,z)|
		\leq C\bigl((t+d(y,z))^{-m-2}
		+(t+d(y,z))^{-n-2}\bigr).
		$$
		
		Fix $z\in B$ and write
		$$
		D=d(x,z),\qquad R=\rho(z),\qquad
		H_z(x)=\sup_{d(x,y)<\min\{2t,16c\rho(x)\}}
		|(\Delta P_t)(y,z)|.
		$$
		For every point in this supremum,
		$$
		t+d(y,z)\geq D/2,\qquad
		(1-16c)\rho(x)\leq\rho(y)\leq(1+16c)\rho(x).
		$$
		Moreover,
		$$
		D\geq r\quad(x\notin2B),\qquad
		R\geq(1-\delta)\rho(x_B),\qquad
		r\leq\frac{\delta}{1-\delta}R.
		$$
		Choose a fixed sufficiently large $R_*$, and put
		$$
		F=E_m\cap\{\rho\geq R_*\},\qquad Q=M\setminus F.
		$$
		The set $Q$ is the smaller end together with a fixed bounded
		set. For every $z\in M$ and $u>0$,
		$$
		\mu(Q\cap B(z,u))\leq Cu^n.
		$$
		For $u\leq1$, this follows from the volume upper bound
		$Cu^m\leq Cu^n$. For $u\geq1$, the bounded part of $Q$
		has fixed volume. If $z\in E_n$, a path of length less than
		$u$ either stays in that end or passes through its fixed
		inner boundary. In its Euclidean coordinates, its endpoint
		therefore satisfies $|X-Z|\leq u+C$. The compact factor has
		fixed volume, so the end contribution is at most $Cu^n$.
		If $z\notin E_n$, every such path reaching $E_n$ enters
		through that inner boundary, and its endpoint has radial
		coordinate at most $u+C$. This gives the same volume bound.
		The general kernel bound above yields
		$$
		H_z(x)\leq C(D^{-m-2}+D^{-n-2})\qquad(x\in Q).
		$$
		Dyadic annuli, the global bound $\mu(B(z,u))\leq Cu^m$,
		and the bound for $Q$ give
		\begin{align*}
			\int_{\{D\geq r\}}D^{-m-2}\dd\mu(x)
			&\leq C\sum_{j\geq0}(2^jr)^{-m-2}(2^{j+1}r)^m
			\leq Cr^{-2},\\
			\int_{Q\cap\{D\geq r\}}D^{-n-2}\dd\mu(x)
			&\leq C\sum_{j\geq0}(2^jr)^{-n-2}(2^{j+1}r)^n
			\leq Cr^{-2}.
		\end{align*}
		
		For $x\in F$, put $s=\rho(x)$. The horizontal restriction
		forces $y\in E_m$ and $\rho(y)\approx s$. If $z\in E_m$,
		its term of dimension $n$ is bounded by
		$$
		Cs^{2-m}R^{2-m}(s+R)^{-n-2}
		\leq Cs^{2-m}(s+R)^{-n-2}.
		$$
		If $z\in E_n$, every path from $y$ to $z$ meets the compact
		part, and
		$$
		d(y,z)\geq\rho(y)+R-2\geq c(s+R).
		$$
		If $z\in K$, then $R=1$ and
		$d(y,z)\geq\rho(y)-1\geq c(s+R)$.
		These cases give the bound
		$$
		H_z(x)\leq CD^{-m-2}
		+Cs^{2-m}(s+R)^{-n-2}\qquad(x\in F).
		$$
		The radial volume element on the large end is at most
		$Cs^{m-1}\,ds$ after integration over the angular variables.
		Hence
		\begin{align*}
			\int_Fs^{2-m}(s+R)^{-n-2}\dd\mu(x)
			&\leq C\int_0^\infty s(s+R)^{-n-2}\,ds\\
			&=\frac{C}{n(n+1)}R^{-n}
			\leq CR^{-2}\leq Cr^{-2}.
		\end{align*}
		Combining the estimates on $Q$ and $F$ proves
		\eqref{eq:prior-kernel}.
	\end{proof}
	
	\begin{lemma}\label{lem:note-atom-forward}
		For sufficiently small fixed $\delta$, every atom in
		\eqref{eq:atoms} satisfies
		$$
		\|\mathcal N_c^\#a\|_1+\|\mathcal S_ca\|_1\leq C.
		$$
	\end{lemma}
	
	\begin{proof}
		Write $a=\Delta b$ and $B=B(x_B,r)$. If $\rho(x_B)$ is large,
		$2B$ stays in one exact product end, where the ball volumes
		have a uniform doubling bound. If $\rho(x_B)$ is bounded,
		the radii are bounded, and the volume upper and lower bounds
		for bounded radii give the same conclusion. Thus
		$\mu(2B)\leq C\mu(B)$. The $L^2$ bound for
		$\mathcal N_c^\#$ therefore gives
		$$
		\int_{2B}\mathcal N_c^\#a\dd\mu
		\leq\mu(2B)^{1/2}\|\mathcal N_c^\#a\|_2
		\leq C\mu(B)^{1/2}\|a\|_2\leq C.
		$$
		The spectral theorem gives $P_ta=\Delta P_tb$. Since $b$ is
		supported in $B$, \eqref{eq:prior-kernel} gives
		\begin{align*}
			\int_{M\setminus2B}\mathcal N_c^\#a\dd\mu
			&\leq\int_B|b(z)|
			\int_{M\setminus2B}
			\sup_{d(x,y)<\min\{2t,16c\rho(x)\}}
			|(\Delta P_t)(y,z)|\dd\mu(x)\dd\mu(z)\\
			&\leq Cr^{-2}\|b\|_1
			\leq Cr^{-2}\mu(B)^{1/2}\|b\|_2\leq C.
		\end{align*}
		Thus $a\in L^2$ and $\mathcal N_ca\leq\mathcal N_c^\#a\in L^1$.
		Lemma~\ref{lem:Hardy-FS-L2} now gives
		$$
		\|\mathcal S_ca\|_1\leq C\|\mathcal N_ca\|_1
		\leq C\|\mathcal N_c^\#a\|_1\leq C.
		$$
				The proof is complete.
	\end{proof}
	
	\subsection{The trace and the Green potential}
	\label{sec:reverse-green}
	
	For the reverse atomic estimate, assume
	$$
	f\in L^2(M),\qquad u(x,t)=P_tf(x),\qquad
	A(x)=\Sc f(x),\qquad A\in L^1(M).
	$$
	We choose the smooth function $R$ in \eqref{eq:radius} with two
	additional properties.
	
	\begin{lemma}\label{lem:reverse-radius}
		The function $R$ may be chosen so that, for some $a>0$,
		\begin{equation}\label{eq:reverse-radius}
			-\Delta\log R\geq aR^{-2},\qquad -\Delta R\geq aR^{-1}.
		\end{equation}
		All its positive-order covariant derivatives are bounded.
	\end{lemma}
	
	\begin{proof}
		Choose a smooth function $v_0$ equal to $\log r$ outside a
		fixed compact set in each end. On the end of dimension $q=m,n$,
		$$
		-\Delta v_0=\frac{q-2}{r^2}.
		$$
		Since $n\geq3$, there are $a_1>0$ and
		$g_0\in C_c^\infty(M)$, $g_0\geq0$, such that
		$$
		-\Delta v_0+g_0\geq a_1\rho^{-2}.
		$$
		Let $G$ be the minimal positive Green operator of $\Delta$ and put
		$v=v_0-Gg_0$. Then
		$$
		-\Delta v=-\Delta v_0+g_0\geq a_1\rho^{-2}.
		$$
		Integration of the heat kernel estimates in \cite{GS} gives,
		on the $q$-dimensional end,
		$$
		0\leq Gg_0(x)\leq Cr^{2-q}.
		$$
		Outside a larger compact set, $Gg_0$ is positive and harmonic.
		The harmonic gradient estimate on a model ball of radius
		comparable to $r$ gives
		$$
		|\nabla Gg_0(x)|\leq Cr^{1-q}.
		$$
		Interior elliptic estimates on fixed-radius model balls give
		$|\nabla^kGg_0(x)|\leq C_kr^{2-q}$ for $k\geq2$.
		On the compact remainder, $Gg_0$ is smooth.
		
		Take $R=C_1e^v$, with $C_1$ large enough that $R\geq1$.
		Then $R\approx\rho$, $|\nabla R|\leq C$, and all its
		positive-order derivatives are bounded. Also,
		$$
		-\Delta\log R=-\Delta v\geq aR^{-2},\qquad
		-\Delta R=R\bigl(-\Delta v+|\nabla v|^2\bigr)\geq aR^{-1}.
		$$
		These properties include \eqref{eq:radius}.
	\end{proof}
	
	Take $\varepsilon>0$ small enough for the proof of
	Lemma~\ref{lem:pointwise} to apply with
	$$
	H(x)=\varepsilon R(x),\qquad g(x)=u(x,H(x)).
	$$
	Require also $4H(x)\leq c\rho(x)/10$. The estimates
	\eqref{eq:gradient} and \eqref{eq:trace-pointwise} then give
	\begin{equation}\label{eq:reverse-local-gradient}
		t|\nabla_{x,t}u(x,t)|\leq CA(x),\qquad 0<t\leq2H(x),
	\end{equation}
	and
	\begin{equation}\label{eq:reverse-trace-difference}
		\rho(x)|\nabla g(x)|\leq C_\varepsilon A(x),\qquad
		|u(x,t)-g(x)|\leq CA(x)\left|\log\frac{t}{H(x)}\right|.
	\end{equation}
	The second inequality holds for $0<t\leq2H(x)$.
	
	\begin{lemma}\label{lem:reverse-trace}
		The trace satisfies
		\begin{equation}\label{eq:reverse-trace}
			\|g\|_1+\|\rho\nabla g\|_1\leq C_\varepsilon\|A\|_1.
		\end{equation}
	\end{lemma}
	
	\begin{proof}
		The gradient estimate follows from
		\eqref{eq:reverse-trace-difference}. By \eqref{eq:decay},
		$g$ tends uniformly to zero along both ends.
		Choose $r_1$ beyond the compact part. On the end of dimension
		$q=m,n$, write
		$$
		\dd\mu=r^{q-1}\dd r\dd\nu_q(\omega),
		$$
		where $\dd\nu_n$ includes the compact sphere variable.
		For almost every $\omega$,
		$$
		g(r,\omega)=-\int_r^\infty\partial_sg(s,\omega)\dd s.
		$$
		The integral is absolutely convergent, since
		$\rho\nabla g\in L^1$ and $r\approx\rho$ on the ends.
		Tonelli gives
		$$
		\begin{aligned}
			\int_{r_1}^\infty|g(r,\omega)|r^{q-1}\dd r
			&\leq\frac1q\int_{r_1}^\infty
			|\partial_sg(s,\omega)|(s^q-r_1^q)\dd s\\
			&\leq C\int_{r_1}^\infty
			\rho(s)|\nabla g(s,\omega)|s^{q-1}\dd s.
		\end{aligned}
		$$
		Let $K_1$ be the compact remainder. Choose a connected bounded
		smooth domain $U$ containing $K_1$ and a fixed end annulus
		$D\subset U$ beyond $r_1$. The $L^1$ Poincar\'e inequality gives
		$$
		\int_{K_1}|g|\dd\mu
		\leq C\left(\int_U|\nabla g|\dd\mu+
		\int_D|g|\dd\mu\right).
		$$
		Together with the two end estimates, this proves
		$\|g\|_1\leq C\|\rho\nabla g\|_1$ and hence
		\eqref{eq:reverse-trace}.
	\end{proof}
	
	\begin{lemma}\label{lem:reverse-harmonic-extension}
		Let
		$$
		\Omega=\{(y,t):t>\varepsilon R(y)\},\qquad
		\Sigma=\partial\Omega,\qquad W(y,t)=tR(y)^{-2}.
		$$
		If $h\geq0$ belongs to $L^1(\Sigma)$ and
		$\mathcal H_\Omega h$ is its minimal harmonic extension, then
		\begin{equation}\label{eq:reverse-weighted-extension}
			\int_\Omega W\mathcal H_\Omega h\dd\mu\dd t
			\leq C_\varepsilon\|h\|_{L^1(\Sigma)}.
		\end{equation}
	\end{lemma}
	
	\begin{proof}
		Recall $\D=-\Delta+\partial_t^2$ on $M\times\mathbb R$. Put
		$$
		\Psi(y,t)=t\log\frac{t}{H(y)}+C_2(t-H(y)),\qquad
		C_2\geq\frac1{a\varepsilon^2}.
		$$
		Then $\Psi\geq0$ on $\Omega$, $\Psi=0$ on $\Sigma$, and
		$$
		\begin{aligned}
			-\D\Psi
			&=t(-\Delta\log R)-\frac1t+C_2\varepsilon(-\Delta R)
			\geq\frac{at}{R^2}-\frac1t+
			\frac{C_2a\varepsilon}{R}
			\geq aW.
		\end{aligned}
		$$
		The outward normal on $\Sigma$ is
		$\nu=(\nabla H,-1)/(1+|\nabla H|^2)^{1/2}$, so
		\begin{equation}\label{eq:reverse-barrier-flux}
			-\partial_\nu\Psi
			=(1+C_2)(1+\varepsilon^2|\nabla R|^2)^{1/2}
			\leq C_\varepsilon.
		\end{equation}
		
		The product heat kernel satisfies
		$$
		\int_1^\infty p_s^{M\times\mathbb R}(z,z)\dd s
		\leq C\int_1^\infty s^{-(n+1)/2}\dd s<\infty.
		$$
		Thus the product has a minimal positive Green kernel.
		Let $G_\Omega$ and $P_\Omega$ be its minimal Dirichlet Green
		and Poisson kernels, obtained by smooth bounded exhaustion.
		Choose $W_j\in C_c^\infty(\Omega)$ with $0\leq W_j\uparrow W$,
		and put
		$$
		V_j(z)=\int_\Omega G_\Omega(z,w)W_j(w)\dd w.
		$$
		On each exhaustion domain, comparison with $\Psi/a$ gives
		$0\leq V_j\leq\Psi/a$ after increasing the domain.
		Both functions vanish on $\Sigma$. Their inward normal
		difference quotients therefore give
		$$
		-\partial_\nu V_j(\xi)
		\leq-a^{-1}\partial_\nu\Psi(\xi)\leq C_\varepsilon.
		$$
		Since $W_j$ is supported inside $\Omega$,
		$$
		-\partial_\nu V_j(\xi)
		=\int_\Omega P_\Omega(z,\xi)W_j(z)\dd z.
		$$
		Monotone convergence gives
		$\int_\Omega P_\Omega(z,\xi)W(z)\dd z\leq C_\varepsilon$.
		Tonelli now yields
		$$
		\begin{aligned}
			\int_\Omega W(z)\mathcal H_\Omega h(z)\dd z
			&=\int_\Sigma h(\xi)
			\left(\int_\Omega W(z)P_\Omega(z,\xi)\dd z\right)
			\dd\Sigma(\xi)
			\leq C_\varepsilon\|h\|_{L^1(\Sigma)}.
		\end{aligned}
		$$
				The proof is complete.
	\end{proof}
	
	\begin{proposition}\label{prop:reverse-green}
		The integral
		$$
		b=\int_0^\infty tP_tf\dd t
		$$
		converges absolutely in $L^1(M,\rho^{-2}\dd\mu)$, satisfies
		$\Delta b=f$ in distributions, and obeys
		\begin{equation}\label{eq:reverse-green-bound}
			\|\rho^{-2}b\|_1+\|\rho^{-1}\nabla b\|_1\leq C\|\Sc f\|_1.
		\end{equation}
		More precisely,
		\begin{equation}\label{eq:reverse-absolute-potentials}
			\int_M\int_0^\infty
			t\bigl(\rho(x)^{-2}|u(x,t)|+
			\rho(x)^{-1}|\nabla_xu(x,t)|\bigr)\dd t\dd\mu(x)
			\leq C\|A\|_1.
		\end{equation}
		Also $b\in W^{2,2}_{\mathrm{loc}}(M)$.
	\end{proposition}
	
	\begin{proof}
		The graph measure satisfies
		$\dd\Sigma=(1+\varepsilon^2|\nabla R|^2)^{1/2}\dd\mu
		\approx\dd\mu$. By Lemma~\ref{lem:reverse-trace},
		$g=u(\cdot,H)$ is bounded, continuous and integrable.
		Put $v=\mathcal H_\Omega|g|$. The Dirichlet Poisson kernel has
		mass at most one, so $v$ is bounded and nonnegative.
		Every sequence escaping compact subsets of $\overline\Omega$
		has $t\to\infty$, since $R$ is proper and $t\geq\varepsilon R$.
		By \eqref{eq:decay}, $u$ tends uniformly to zero there.
		The maximum principle on bounded exhaustion domains gives
		$|u|\leq v+\eta_j$, where $\eta_j$ is the supremum of $|u|$
		on their artificial boundaries and $\eta_j\to0$.
		For complex $u$, apply this argument to
		$\operatorname{Re}(e^{i\theta}u)$ for every $\theta$.
		Consequently $|u|\leq v$ on $\Omega$, and
		Lemma~\ref{lem:reverse-harmonic-extension} gives
		\begin{equation}\label{eq:reverse-upper-integral}
			\int_M R(x)^{-2}\int_{H(x)}^\infty t|u(x,t)|\dd t\dd\mu(x)
			\leq C\|g\|_1\leq C\|A\|_1.
		\end{equation}
		
		For $0<t<H(x)$, use \eqref{eq:reverse-trace-difference} and
		$$
		\int_0^H t\dd t=\frac{H^2}{2},\qquad
		\int_0^H t\log(H/t)\dd t=\frac{H^2}{4}.
		$$
		They give
		$$
		\int_M R^{-2}\int_0^H t|u|\dd t\dd\mu
		\leq C_\varepsilon(\|g\|_1+\|A\|_1)\leq C\|A\|_1.
		$$
		This proves the first term of \eqref{eq:reverse-absolute-potentials}.
		
		For the spatial gradient, \eqref{eq:reverse-local-gradient} gives
		$$
		R(x)^{-1}\int_0^{H(x)}t|\nabla_xu(x,t)|\dd t\leq CA(x).
		$$
		For $t\geq H(x)$, take a product ball $Q(x,t)$ of radius
		$\kappa R(x)$ centered at $(x,t)$, where $\kappa>0$ is small
		compared with $\varepsilon$. Lemma~\ref{lem:local-harmonic-gradient}
		and the local $L^1$ mean value inequality for $|u|$ give,
		after a fixed enlargement of this ball,
		$$
		|\nabla_xu(x,t)|
		\leq\frac{C}{R(x)|Q(x,t)|}
		\int_{Q(x,t)}|u(y,s)|\dd\mu(y)\dd s.
		$$
		The balls used here lie in one complete model or a fixed bounded
		part of the product. On them,
		$$
		R(y)\approx R(x),\qquad s\approx t,\qquad
		s\geq\tfrac12\varepsilon R(y).
		$$
		For each $(y,s)$ in such a ball, local doubling gives
		$|Q(x,t)|\approx|Q(y,s)|$. Its possible centers $(x,t)$ lie
		in a fixed enlargement of $Q(y,s)$. Hence
		$$
		\int_{\{(x,t):(y,s)\in Q(x,t)\}}
		\frac{\dd\mu(x)\dd t}{|Q(x,t)|}\leq C.
		$$
		Tonelli therefore gives
		$$
		\int_{t\geq H(x)}\frac{t}{R(x)}|\nabla_xu(x,t)|\dd\mu(x)\dd t
		\leq C\int_{s\geq\varepsilon R(y)/2}
		\frac{s}{R(y)^2}|u(y,s)|\dd\mu(y)\dd s.
		$$
		Apply \eqref{eq:reverse-upper-integral} with $\varepsilon/2$.
		This proves \eqref{eq:reverse-absolute-potentials}.
		Local absolute convergence gives, in distributions,
		$$
		\nabla b=\int_0^\infty t\nabla_xP_tf\dd t.
		$$
		It also proves \eqref{eq:reverse-green-bound}.
		
		For $0<\alpha<T$, integration by parts in $t$ gives
		$$
		\Delta\int_\alpha^T tP_tf\dd t
		=T\partial_TP_Tf-\alpha\partial_\alpha P_\alpha f-P_Tf+P_\alpha f.
		$$
		The spectral theorem implies
		$$
		\begin{aligned}
			P_\alpha f\to f,\quad
			\alpha\partial_\alpha P_\alpha f\to0
			&\quad(\alpha\downarrow0),\\
			P_Tf\to0,\quad T\partial_TP_Tf\to0
			&\quad(T\to\infty)
		\end{aligned}
		$$
		in $L^2$. Here $\ker_{L^2}\Delta=\{0\}$: a vector in the kernel
		has zero gradient and is constant, while $M$ is connected and
		has infinite volume. Pair the preceding identity with a compactly
		supported smooth function. Weighted $L^1$ convergence on the left
		then gives $\Delta b=f$ in distributions. Interior elliptic regularity
		gives $b\in W^{2,2}_{\mathrm{loc}}(M)$.
	\end{proof}
	
	\subsection{Estimates on the complete models}
	\label{sec:complete-model-estimates}
	
	We use auxiliary atomic norms on the complete models to estimate
	$\chi_jf$ and $\Delta(\chi_jb)$. Summing these estimates over $j$
	gives a bound by $\|\Sc f\|_1$.
	
	Write $A(x)=\Sc f(x)$ and $u(x,t)=P_tf(x)$. We use the bounds
	\begin{equation}\label{eq:reverse-local-known}
		\|g\|_1+\|\rho^{-2}b\|_1
		+\|\rho^{-1}\nabla b\|_1\leq C\|A\|_1,
		\qquad \Delta b=f,
	\end{equation}
	where $g(x)=u(x,H(x))$, $H(x)=\varepsilon R(x)$ for the smooth
	radius $R\approx\rho$ chosen above, and
	\begin{equation}\label{eq:reverse-local-pointwise}
		t|\nabla_{x,t}u(x,t)|\leq CA(x)
		\quad\text{for }0<t\leq2H(x).
	\end{equation}
	
	Choose a nonnegative smooth partition of unity $\{\chi_j\}$ with
	$$
	\operatorname{supp}\chi_j\subset B_j=B(x_j,r_j),\qquad
	r_j=\nu_0\rho(x_j),\qquad
	|\nabla^k\chi_j|\leq C_k r_j^{-k}\quad(k\leq3).
	$$
	Take $\nu_0$ small relative to $\varepsilon$ so that
	$2r_j\leq H(x)$ on every fixed enlargement used below. Require
	also that all these enlargements
	lie in one model end, except for finitely many balls in a fixed
	compact set. The enlarged balls have bounded overlap. Indeed,
	$\rho$ is Lipschitz, so intersecting enlarged balls have comparable
	radii, and the model volume bounds give the packing estimate.
	On the product end, if $r_j$ is larger than a fixed constant,
	choose the cutoff independent of the sphere variable. Its support
	may be placed in a product cylinder contained in $B_j$, with
	Euclidean radius comparable to $r_j$. The volume of this cylinder
	is comparable to $\mu(B_j)$. Thus $B_j$ always denotes a ball.
	
	For an end ball, let $X_j$ be its complete model and let
	$\Delta_j$ be the model Laplacian. First choose a finite
	coordinate cover of the fixed compact region. On each chart,
	extend the metric from a larger coordinate neighbourhood to a
	smooth metric on $\mathbb R^m$ which is Euclidean outside a
	compact set. Its density and metric are uniformly comparable
	to the Euclidean ones. These finitely many complete models are
	fixed before $\nu_0$ is chosen. Assign each central ball to a
	chart containing all its required enlargements and use the
	corresponding complete model as $X_j$. Thus $\Delta_j=\Delta$
	on every enlarged ball where it is used. These complete models are doubling. The constants
	on the two end models are uniform for every $r_j>0$; the finitely
	many central choices have a common constant.
	
	The common region is chosen larger than every fixed enlargement
	used below. If a curve of length at most $C_1r_j$ starts in
	$C_1B_j$, it cannot leave that region: the distance to its boundary
	is greater than $2C_1r_j$. This applies to curves on $M$ and on
	$X_j$. Hence, for the short distances and radii used below,
	$$
	d_M(x,y)=d_j(x,y),\qquad
	B_M(y,t)=B_{X_j}(y,t),\qquad V_M(y,t)=V_j(y,t).
	$$
	On an end this follows by taking $\nu_0$ small. On a central patch
	it follows by taking the coordinate neighbourhood larger first
	and then decreasing $\nu_0$.
	
	The operator domains also agree for the compactly supported
	functions used below. Suppose that
	$\beta\in\operatorname{Dom}(\Delta_j)$ and its support
	is a compact subset of the common region. Its extension by zero
	belongs to the form domain of $\Delta$. Choose a smooth cutoff
	$\zeta=1$ near its support, contained in that region. For
	$\varphi\in C_c^\infty(M)$, the model form identity gives
	$$
	\int_M\nabla\beta\cdot\nabla\overline\varphi\dd\mu
	=\int_{X_j}\nabla\beta\cdot
	\nabla\overline{\zeta\varphi}\dd\mu_j
	=\int_{X_j}\Delta_j\beta\,
	\overline{\zeta\varphi}\dd\mu_j
	=\int_M\Delta_j\beta\,\overline\varphi\dd\mu.
	$$
	Locality puts $\supp\Delta_j\beta$ in $\supp\beta$.
	The last function belongs to $L^2$. Density in the form domain
	therefore gives $\beta\in\operatorname{Dom}(\Delta)$ and
	$\Delta\beta=\Delta_j\beta$.

	\subsubsection{Auxiliary atomic spaces and their duals}
	
	Let $X$ be one of these models, with measure $\nu$ and
	nonnegative Laplacian $\Delta_X$, and let $r>0$.
	For a central model, $r$ is one of the finitely many radii just
	chosen. An auxiliary atom is a function $a$ supported in a ball
	$B$ of radius $s\leq r$ such that
	$$
	\|a\|_2\leq\nu(B)^{-1/2},\qquad
	\int a\,d\nu=0\quad\text{if }s<r.
	$$
	Let $\mathfrak A_r(X)$ consist of sums $F=\sum_k\lambda_k a_k$
	in $L^1(X)$ with $\sum_k|\lambda_k|<\infty$. Its norm is the
	infimum of this coefficient sum. We omit $X$ when it is fixed.
	Its dual norm is equivalent to
	\begin{equation}\label{eq:reverse-local-bmo}
		\|\varphi\|_{\operatorname{bmo}_r}
		=\sup_{s(B)<r}\left(\frac1{\nu(B)}
		\int_B|\varphi-\varphi_B|^2\,d\nu\right)^{1/2}
		+\sup_{s(B)\geq r}\left(\frac1{\nu(B)}
		\int_B|\varphi|^2\,d\nu\right)^{1/2}.
	\end{equation}
	Here and below $s(B)$ denotes the radius of $B$.
	Testing on atoms proves one direction of duality. Conversely,
	a bounded functional restricts to an $L^2$ functional on every
	ball of radius $r$, since all normalized $L^2$ functions on such
	a ball are atoms. The resulting local $L^2$ functions agree on
	overlaps. Testing on mean-zero $L^2$ functions gives the first
	term of \eqref{eq:reverse-local-bmo}. Covering larger balls by
	balls of radius $r$ gives the second term. Hahn--Banach therefore
	gives, for $F\in \mathfrak A_r$,
	\begin{equation}\label{eq:reverse-local-duality}
		\|F\|_{\mathfrak A_r}\approx
		\sup_{\|\varphi\|_{\operatorname{bmo}_r}\leq1}
		\left|\int_XF\overline\varphi\,d\nu\right|.
	\end{equation}
	All constants depend only on the doubling constants. In particular,
	an $L^2$ function supported in a ball of radius $r$ already belongs
	to $\mathfrak A_r$.
	
	\begin{lemma}\label{lem:reverse-model-tent}
		Let $X$ be one of the complete models above. Put
		$$
		\mathcal AF(x)^2=\int_0^\infty\int_{B(x,t)}
		|F(y,t)|^2\frac{d\nu(y)\,dt}{t\nu(B(y,t))},
		\qquad \|F\|_{T^p}=\|\mathcal AF\|_p.
		$$
		Every $F\in T^1$ has a representation
		$$
		F=\sum_\ell\lambda_\ell A_\ell,\qquad
		\sum_\ell|\lambda_\ell|\leq C\|F\|_{T^1},
		$$
		where, for a ball $B_\ell$,
		$$
		\operatorname{supp}A_\ell\subset B_\ell\times(0,r_{B_\ell}),
		\qquad \int_0^{r_{B_\ell}}\|A_\ell(\cdot,t)\|_2^2\,
		\frac{dt}{t}\leq\nu(B_\ell)^{-1}.
		$$
		The series converges in $T^1$. If $F\in T^2$, it also converges
		in $T^2$. If
		$$
		\operatorname{supp}F\subset 2B_0\times(0,r_0),
		\qquad B_0=B(x_0,r_0),
		$$
		all the balls may be chosen with
		$B_\ell\subset CB_0$ and $r_{B_\ell}\leq Cr_0$.
		Moreover, if
		$$
		\|G\|_{T^\infty}
		=\sup_B\left(\frac1{\nu(B)}\int_0^{r_B}\int_B
		|G(y,t)|^2\,d\nu(y)\,\frac{dt}{t}\right)^{1/2},
		$$
		then
		\begin{equation}\label{eq:reverse-model-tent-duality}
			\int_0^\infty\int_X|F(y,t)|\,|G(y,t)|\,d\nu(y)\,
			\frac{dt}{t}
			\leq C\|F\|_{T^1}\|G\|_{T^\infty}.
		\end{equation}
		The same assertions hold for vector fields.
	\end{lemma}
	
	\begin{proof}
		We give the construction. It is the usual tent decomposition
		on a doubling space; see \cite{Russ}. Let $\mathcal M_X$ be the
		centred Hardy--Littlewood maximal operator. For $k\in\mathbb Z$,
		put
		$$
		O_k=\{\mathcal AF>2^k\},\qquad
		O_k^*=\{\mathcal M_X\mathbf1_{O_k}>\alpha\},\qquad
		T_k=\{(y,t):B(y,t)\subset O_k^*\}.
		$$
		Choose $\alpha>0$ so small that
		$$
		\nu(B(x,2t))\leq C_D\nu(B(y,t)),\ d(x,y)<t
		\quad\Longrightarrow\quad \alpha C_D\leq\tfrac12.
		$$
		Doubling gives such a constant $C_D$. The maximal theorem gives
		\begin{equation}\label{eq:reverse-tent-level-volumes}
			\nu(O_k^*)\leq C_\alpha\nu(O_k).
		\end{equation}
		Choose Whitney balls $D_{k\ell}$ and a nonnegative partition
		$\vartheta_{k\ell}$ on $O_k^*$ such that
		$$
		\sum_\ell\vartheta_{k\ell}=1,\qquad
		\operatorname{supp}\vartheta_{k\ell}\subset2D_{k\ell},\qquad
		\sum_\ell\mathbf1_{2D_{k\ell}}\leq C,
		$$
		and
		$$
		r_{D_{k\ell}}\approx
		d(y,(O_k^*)^c)\quad(y\in2D_{k\ell}),\qquad
		\sum_\ell\nu(D_{k\ell})\leq C\nu(O_k^*).
		$$
		For example, select disjoint balls whose radii are a fixed small
		multiple of their distances to $(O_k^*)^c$, enlarge them to cover
		$O_k^*$, and divide the corresponding nonnegative cutoff functions
		by their sum. Intersecting enlarged balls have comparable radii.
		Doubling gives the overlap and volume bounds.
		
		Set
		$$
		F_{k\ell}(y,t)
		=F(y,t)\mathbf1_{T_k\setminus T_{k+1}}(y,t)
		\vartheta_{k\ell}(y).
		$$
		Its support is contained in
		$2D_{k\ell}\times(0,Cr_{D_{k\ell}})$. If
		$(y,t)\notin T_{k+1}$, choose
		$x_1\in B(y,t)\setminus O_{k+1}^*$. Then
		$$
		\nu(O_{k+1}\cap B(y,t))
		\leq\alpha\nu(B(x_1,2t))
		\leq\tfrac12\nu(B(y,t)).
		$$
		Thus at least half of $B(y,t)$ lies outside $O_{k+1}$.
		Choose a fixed enlargement $B_{k\ell}=C_1D_{k\ell}$ which
		contains every $B(y,t)$ occurring in $F_{k\ell}$. Tonelli gives
		$$
		\begin{aligned}
			\int_0^\infty\|F_{k\ell}(\cdot,t)\|_2^2\,\frac{dt}{t}
			&\leq2\int_{B_{k\ell}\setminus O_{k+1}}
			\int_0^\infty\int_{B(x,t)}
			|F_{k\ell}(y,t)|^2
			\frac{d\nu(y)\,dt}{t\nu(B(y,t))}\,d\nu(x)\\
			&\leq2\int_{B_{k\ell}\setminus O_{k+1}}
			\mathcal AF(x)^2\,d\nu(x)\\
			&\leq C2^{2k}\nu(B_{k\ell}).
		\end{aligned}
		$$
		For a fixed sufficiently large $C_2$, take
		$$
		\lambda_{k\ell}=C_2 2^k\nu(B_{k\ell}),\qquad
		A_{k\ell}=\lambda_{k\ell}^{-1}F_{k\ell},
		$$
		omitting zero pieces. These are the asserted atoms. Their
		coefficients satisfy
		$$
		\begin{aligned}
			\sum_{k,\ell}|\lambda_{k\ell}|
			&\leq C\sum_k2^k\nu(O_k^*)
			\leq C\sum_k2^k\nu(O_k)
			\leq C\int_X\mathcal AF\,d\nu.
		\end{aligned}
		$$
		
		For almost every $(y,t)$ with $F(y,t)\ne0$, the area function
		is positive on $B(y,t)$, up to a null set. This follows by taking
		$(y,t)$ to be a density point of a set where $|F|$ is bounded
		below. Consequently,
		$$
		\nu(O_k\cap B(y,t))\longrightarrow\nu(B(y,t))
		\quad(k\to-\infty).
		$$
		The preceding doubling comparison, with the chosen $\alpha$,
		then gives $B(y,t)\subset O_k^*$ for all sufficiently negative
		$k$. In the other direction, \eqref{eq:reverse-tent-level-volumes}
		and $\nu(O_k)\to0$ show that no fixed ball $B(y,t)$ lies in
		$O_k^*$ for all sufficiently positive $k$. Hence
		$$
		F=\sum_{k,\ell}F_{k\ell}\quad\text{almost everywhere}.
		$$
		Each piece is $F$ multiplied by a nonnegative function, and these
		functions sum to one. Every tent atom has $T^1$ norm at most $C$:
		its area function is supported in $2B$, and
		$$
		\|A\|_{T^1}
		\leq\nu(2B)^{1/2}\|A\|_{T^2}
		\leq C.
		$$
		The coefficient bound therefore gives convergence in $T^1$.
		If $F\in T^2$, dominated convergence gives convergence in
		$L^2(d\nu\,dt/t)$ and hence in $T^2$, since Tonelli gives
		$$
		\|F\|_{T^2}^2
		=\int_0^\infty\|F(\cdot,t)\|_2^2\,\frac{dt}{t}.
		$$
		
		For the additional support assertion,
		$\operatorname{supp}\mathcal AF\subset3B_0$. The volume bounds
		on our complete models imply, for $d(x,x_0)\geq Cr_0$,
		$$
		\mathcal M_X\mathbf1_{3B_0}(x)
		\leq C\left(\frac{r_0}{d(x,x_0)}\right)^n.
		$$
		Indeed, every ball centred at $x$ which meets $3B_0$ has radius
		at least $d(x,x_0)-3r_0$. On $\mathbb R^m$ and on each central
		model its volume is comparable to the $m$th power of its radius.
		On the product model the volume is comparable to
		$\min(s^m,s^n)$, whose growth from $r_0$ to a larger radius is
		at least the $n$th power of the ratio. Thus $O_k^*\subset CB_0$.
		The Whitney radii and their fixed enlargements are at most $Cr_0$.
		
		Finally, Cauchy--Schwarz on the box of each tent atom gives
		$$
		\int_0^\infty\int_X|A_\ell|\,|G|\,d\nu\,\frac{dt}{t}
		\leq\nu(B_\ell)^{-1/2}
		\left(\int_0^{r_{B_\ell}}\int_{B_\ell}|G|^2
		\,d\nu\,\frac{dt}{t}\right)^{1/2}
		\leq\|G\|_{T^\infty}.
		$$
		Summing the nonnegative pieces proves
		\eqref{eq:reverse-model-tent-duality}.
	\end{proof}

	We record three estimates for these auxiliary spaces on the
	complete models.
	
	\begin{lemma}\label{lem:reverse-model-estimates}
		The following constants are uniform for the model choices above.
		\begin{enumerate}[(i)]
			\item If $v=e^{-t\sqrt{\Delta_X}}\varphi$ and
			$\|\varphi\|_{\operatorname{bmo}_r}\leq1$, then, for $0<t\leq Cr$,
			\begin{equation}\label{eq:reverse-bmo-poisson}
				|v(x,t)|\leq C\bigl(1+\log^+(r/t)\bigr),\qquad
				t|\nabla_{x,t}v(x,t)|\leq C.
			\end{equation}
			For balls $B$ with $s(B)\leq Cr$,
			\begin{equation}\label{eq:reverse-bmo-carleson}
				\int_0^{s(B)}\int_B t|\nabla_{x,t}v|^2\,d\nu\,dt
				\leq C\nu(B).
			\end{equation}
			\item If $\|\beta\|_\infty+r\operatorname{Lip}(\beta)\leq C_0$,
			then multiplication by $\beta$ is bounded on $\mathfrak A_r$, with norm
			at most $CC_0$.
			\item Put $J_r=(I+r^2\Delta_X)^{-1}$. For $F\in L^1(X)$,
			\begin{equation}\label{eq:reverse-spectral-smoothing}
				\|J_rF\|_{\mathfrak A_r}\leq C\|F\|_1,\qquad
				\|D J_rF\|_{\mathfrak A_r}\leq Cr^{-1}\|F\|_1.
			\end{equation}
			Here $D$ is a Euclidean derivative on an end, a rotation field on
			the sphere, or one of the fixed smooth coordinate fields used on
			a central patch. For central fields one may insert a fixed cutoff
			supported in the chosen coordinate neighbourhood.
		\end{enumerate}
	\end{lemma}
	
	\begin{proof}
		The end models have Gaussian heat kernel bounds and first
		derivative bounds for every positive heat time. Their ball
		volumes are $s^m$ on $\mathbb R^m$ and comparable to
		$\min(s^m,s^n)$ on the product. On each central model the
		coefficients and density are smooth, uniformly elliptic, and
		constant outside a compact set. Its heat kernel has Gaussian
		bounds. For heat times $0<s\leq1$, interior parabolic derivative
		estimates give the corresponding first derivative bounds.
		For $s\geq1$, apply those estimates on fixed balls to obtain
		$$
		|\nabla_x e^{-s\Delta_X}(x,z)|
		\leq Cs^{-m/2}\exp\bigl(-d(x,z)^2/(Cs)\bigr).
		$$
		In the subordination formula, the part with $s\geq1$ is therefore
		a smooth kernel $R_t$ satisfying, for $t$ in any fixed bounded
		interval,
		$$
		|R_t(x,z)|+t|\partial_tR_t(x,z)|+|\nabla_xR_t(x,z)|
		\leq Ct(1+d(x,z))^{-m-1}.
		$$
		Its integral against $|\varphi|$ is bounded by $Ct$, since the
		large-ball averages in \eqref{eq:reverse-local-bmo} are bounded.
		The part with $s\leq1$ has the Poisson kernel and derivative
		bounds used below. All the heat kernels preserve constants.
		
		For (i), averages of $\varphi$ over nested balls satisfy
		$$
		|\varphi_{B(x,t)}|\leq C\bigl(1+\log^+(r/t)\bigr),\qquad
		\left(\frac1{\nu(2^kB)}\int_{2^kB}|\varphi-\varphi_B|^2\,d\nu\right)^{1/2}
		\leq C(k+1).
		$$
		The first assertion follows by telescoping up to radius $r$.
		The second follows from the same telescoping and
		\eqref{eq:reverse-local-bmo}. Integrating over the annuli
		$2^{k+1}B\setminus2^kB$ proves the first inequality in
		\eqref{eq:reverse-bmo-poisson}. For the gradient, subtract
		$\varphi_{B(x,t)}$ and use that the derivative of the Poisson
		kernel has integral zero. The resulting series is bounded by
		$Ct^{-1}\sum_{k\geq0}(k+1)2^{-k}$. For a central model, the
		smooth term just estimated is bounded by the same expression.
		
		For \eqref{eq:reverse-bmo-carleson}, subtract $\varphi_{2B}$ and
		split the difference into its restrictions to $2B$ and $(2B)^c$.
		The $L^2$ energy identity bounds the first contribution by
		$C\nu(B)$. The kernel derivative bounds and the preceding annular
		estimate bound the gradient of the second contribution by
		$C/s(B)$ on $B\times(0,s(B))$. Its weighted square integral is
		also at most $C\nu(B)$. Constant terms have zero gradient.
		
		For (ii), let $a$ be an atom in a ball of radius $s<r$. Write
		$$
		\beta a=\bigl(\beta a-c_a\nu(B)^{-1}\mathbf1_B\bigr)
		+c_a\nu(B)^{-1}\mathbf1_B,
		\qquad c_a=\int_B\beta a\,d\nu.
		$$
		The first term is a fixed multiple of an atom, and
		$|c_a|\leq C_0s/r$. Telescoping normalized indicators over doubled
		balls up to radius $r$ shows that
		\begin{equation}\label{eq:reverse-small-support}
			\|F\|_{\mathfrak A_r}\leq C\nu(B)^{1/2}\|F\|_2
			+C\bigl(1+\log(r/s)\bigr)\|F\|_1
			\quad\text{if }\operatorname{supp}F\subset B,\ s\leq r.
		\end{equation}
		Since $(s/r)(1+\log(r/s))\leq C$, (ii) follows. At radius $r$
		there is no cancellation condition.
		
		For (iii), write the Bessel kernel as
		$$
		J_r(x,z)=r^{-2}\int_0^\infty e^{-s/r^2}
		e^{-s\Delta_X}(x,z)\,ds.
		$$
		For $a=d(x,z)\leq r$, its kernel and first derivatives satisfy
		\begin{equation}\label{eq:reverse-bessel-kernels}
			|J_r(x,z)|\leq\frac{Ca^2}{r^2\nu(B(z,a))},\qquad
			|D_xJ_r(x,z)|\leq\frac{Ca}{r^2\nu(B(z,a))}.
		\end{equation}
		For $a>r$ they are bounded, respectively, by
		$$
		\frac{C}{\nu(B(z,r))}e^{-ca/r},\qquad
		\frac{C}{r\nu(B(z,r))}e^{-ca/r}.
		$$
		To obtain \eqref{eq:reverse-bessel-kernels}, split the heat
		integral at $s=a^2$ and use the Gaussian bounds. The smallest
		volume exponent on either complete end model is $n>2$, so the
		part with $s>a^2$ is bounded by the displayed expression.
		The central models have local dimension $m>2$; their remaining
		integrals are bounded smooth terms at the fixed radii in question.
		
		Apply \eqref{eq:reverse-small-support} on
		$2^{-k-1}r<d(x,z)\leq2^{-k}r$. The contributions are at most
		$C(k+1)2^{-2k}$ for $J_r(\cdot,z)$ and
		$Cr^{-1}(k+1)2^{-k}$ for $D_xJ_r(\cdot,z)$.
		For the outer annuli, cover by balls of radius $r$ and use the
		atoms without cancellation; polynomial growth is absorbed by
		$e^{-ca/r}$. Summing gives
		$$
		\sup_z\|J_r(\cdot,z)\|_{\mathfrak A_r}\leq C,\qquad
		\sup_z\|D_xJ_r(\cdot,z)\|_{\mathfrak A_r}\leq Cr^{-1}.
		$$
		Integration against $F(z)$ proves (iii).
	\end{proof}
	
	\subsubsection{Estimates for the boundary pieces}
	
	\begin{lemma}\label{lem:reverse-boundary-atoms}
		For the partition above,
		\begin{equation}\label{eq:reverse-boundary-atoms}
			\sum_j\|\chi_j f\|_{\mathfrak A_{r_j}(X_j)}
			\leq C\bigl(\|A\|_1+\|g\|_1\bigr)
			\leq C\|A\|_1.
		\end{equation}
	\end{lemma}
	
	\begin{proof}
		Fix $j$, write $r=r_j$, and take
		$\|\varphi\|_{\operatorname{bmo}_r(X_j)}\leq1$.
		Put $v=e^{-t\sqrt{\Delta_j}}\varphi$. On the enlarged
		patch, $u$ and $v$ solve the same harmonic equation. Let
		$\theta=1$ on $[0,1]$, $\theta=0$ on $[2,\infty)$, and set
		$$
		w(x,t)=t\theta(t/r)\chi_j(x),\qquad
		g_r(x)=u(x,r).
		$$
		The radius choices and \eqref{eq:reverse-local-pointwise} imply,
		for $x$ on the support of $\chi_j$ or its derivatives,
		\begin{equation}\label{eq:reverse-local-log-u}
			|g_r(x)|\leq|g(x)|+CA(x),\qquad
			|u(x,t)|\leq|g_r(x)|+CA(x)|\log(t/r)|\quad(0<t<2r).
		\end{equation}
		The ratio $H(x)/r$ lies between two fixed positive constants.
		
		On the common patch, put $\mathcal D=-\Delta+\partial_t^2$.
		Then
		$$
		\mathcal Du=\mathcal Dv=0,\qquad
		\mathcal D(u\overline v)=2\nabla u\cdot\nabla\overline v.
		$$
		Integration by parts gives
		\begin{equation}\label{eq:reverse-local-green}
			\int_M\chi_jf\overline\varphi\,d\mu
			=2\iint w\nabla u\cdot\nabla\overline v
			-\iint w_{tt}u\overline v
			+\iint\nabla_xw\cdot\nabla_x(u\overline v).
		\end{equation}
		All integrals on the right use $d\mu\,dt$.
		For the first term, put
		$$
		F=t\nabla u\,\mathbf1_{B_j\times(0,2r)},\qquad
		G=t\nabla v\,\mathbf1_{B_j\times(0,2r)}.
		$$
		The area function of $F$ is supported in $3B_j$ and satisfies
		$$
		\mathcal AF(x)\leq CA(x),\qquad
		\|F\|_{T^1}\leq C\int_{3B_j}A\,d\mu.
		$$
		Indeed, $t<2r\ll c\rho(x)$ on every contributing cone, so these
		balls lie in the original restricted cones. The distances and
		volumes agree by the choice of the common region.
		Equation \eqref{eq:reverse-bmo-carleson} gives
		$\|G\|_{T^\infty}\leq C$. For a ball of radius greater than $Cr$,
		its box integral is bounded by the integral over
		$B_j\times(0,2r)$; if the ball meets $B_j$, doubling compares
		its volume with that of $B_j$. Smaller balls are covered by
		\eqref{eq:reverse-bmo-carleson}. Since $|w|\leq Ct$, Lemma
		\ref{lem:reverse-model-tent} gives
		$$
		\left|\iint w\nabla u\cdot\nabla\overline v\right|
		\leq C\int_0^\infty\int_{X_j}|F|\,|G|\,d\mu_j\,\frac{dt}{t}
		\leq C\int_{3B_j}A\,d\mu.
		$$
		
		The temporal term is supported on $r<t<2r$, where
		$|w_{tt}|\leq C/r$ and $|v|\leq C$. Its absolute value is at most
		$C\int_{B_j}(|g|+A)\,d\mu$. For the spatial term, use
		$|\nabla_xw|\leq Ct/r$ and
		\eqref{eq:reverse-bmo-poisson}--\eqref{eq:reverse-local-log-u}:
		$$
		\begin{aligned}
			\iint|\nabla_xw|\,|\nabla_x(u\overline v)|
			&\leq\frac Cr\int_{B_j}\int_0^{2r}
			\bigl[|u|+A(x)(1+\log^+(r/t))\bigr]\,dt\,d\mu(x)
			\leq C\int_{B_j}(|g|+A)\,d\mu.
		\end{aligned}
		$$
		Here $\int_0^{2r}(1+|\log(t/r)|)\,dt\leq Cr$.
		
		To justify \eqref{eq:reverse-local-green}, choose a smooth
		$\zeta$ with $\zeta=0$ on $[0,1]$ and $\zeta=1$ on $[2,\infty)$,
		and put
		$$
		h_\tau(t)=t\zeta(t/\tau),\qquad
		w_\tau(x,t)=h_\tau(t)\theta(t/r)\chi_j(x),\qquad 0<\tau<r/4.
		$$
		Then
		$$
		\operatorname{supp}h_\tau''\subset[\tau,2\tau],\qquad
		\int_0^\infty h_\tau''(t)\,dt=1,\qquad
		\int_0^\infty|h_\tau''(t)|\,dt\leq C.
		$$
		Integration by parts with $w_\tau$ has no boundary term. Move
		the term containing $h_\tau''$ to the other side. Since
		$u(\cdot,t)\to f$ and $v(\cdot,t)\to\varphi$ in $L^2$ on the
		fixed compact patch,
		$$
		\int_0^\infty h_\tau''(t)
		\int_M\chi_j u(x,t)\overline{v(x,t)}\,d\mu(x)\,dt
		\longrightarrow\int_M\chi_j f\overline\varphi\,d\mu.
		$$
		The displayed absolute bounds justify all the other limits.
		The convergence of $v$ in $L^2$ on the fixed patch follows
		by multiplying $\varphi$ by a compactly supported cutoff equal
		to one on a larger patch. For that part, use strong continuity
		of the model Poisson semigroup on $L^2$. The other part is a
		positive distance away. Its Poisson integral tends uniformly
		to zero on the fixed patch, by the kernel bounds and the bounded
		large-ball averages of $|\varphi|$. All integrations by parts
		are made at positive heights, where $v$ is smooth.
		
		Since $\chi_j f$ is compactly supported in $B_j$ and belongs to
		$L^2$, it already belongs to $\mathfrak A_r$. Taking the supremum in
		\eqref{eq:reverse-local-duality} gives
		$$
		\|\chi_jf\|_{\mathfrak A_{r_j}(X_j)}
		\leq C\int_{3B_j}(A+|g|)\,d\mu.
		$$
		Bounded overlap proves \eqref{eq:reverse-boundary-atoms}.
	\end{proof}
	
	Every auxiliary atom defined above has $L^1$ norm at most one.
	Since $\sum_j\chi_j=1$, Lemma~\ref{lem:reverse-boundary-atoms}
	also gives
	\begin{equation}\label{eq:reverse-local-L1}
		\|f\|_1\leq\sum_j\|\chi_j f\|_1
		\leq\sum_j\|\chi_j f\|_{\mathfrak A_{r_j}(X_j)}
		\leq C\|\Sc f\|_1.
	\end{equation}
	
	\subsubsection{Commutators with the spatial cutoffs}
	
	Set $h_j=\chi_jb$ and $a_j=\Delta h_j$. These functions are regarded
	as compactly supported functions on $X_j$. Since $\Delta$ is the
	nonnegative Laplacian,
	\begin{equation}\label{eq:reverse-local-commutator}
		a_j=\chi_jf-2\langle\nabla\chi_j,\nabla b\rangle+(\Delta\chi_j)b.
	\end{equation}
	
	\begin{lemma}\label{lem:reverse-potential-atoms}
		One has
		\begin{equation}\label{eq:reverse-potential-atoms}
			\sum_j\left(\|a_j\|_{\mathfrak A_{r_j}(X_j)}
			+r_j^{-2}\|h_j\|_1\right)
			\leq C\|\Sc f\|_1.
		\end{equation}
	\end{lemma}
	
	\begin{proof}
		Choose $\eta_j=1$ on a neighbourhood of $\operatorname{supp}\chi_j$,
		with support in $2B_j$ and
		$|\nabla^k\eta_j|\leq C_kr_j^{-k}$ for $k\leq2$.
		Put $B=\eta_jb$, $F=\Delta_jB$, and $r=r_j$. Then
		\begin{equation}\label{eq:reverse-local-F-bound}
			\|F\|_1\leq C\left(\|f\|_{L^1(2B_j)}
			+r^{-1}\|\nabla b\|_{L^1(2B_j)}
			+r^{-2}\|b\|_{L^1(2B_j)}\right).
		\end{equation}
		The identity
		$$
		B=J_rB+r^2J_rF
		$$
		and \eqref{eq:reverse-spectral-smoothing} give
		\begin{equation}\label{eq:reverse-B-atomic-norm}
			\|B\|_{\mathfrak A_r}\leq C(\|B\|_1+r^2\|F\|_1),\qquad
			\|DB\|_{\mathfrak A_r}\leq C(r^{-1}\|B\|_1+r\|F\|_1).
		\end{equation}
		On Euclidean end patches use the coordinate derivatives. On
		large product-end patches the cutoffs are independent of the
		sphere variable. On the remaining product patches, the finitely
		many rotation fields on the sphere express the scalar product
		of gradients. On central patches use the fixed coordinate
		fields. In each case the coefficients in
		$r\langle\nabla\chi_j,\nabla B\rangle$ and $r^2(\Delta\chi_j)B$
		have bounded supremum and $r$ times their Lipschitz constant
		bounded uniformly. Part (ii) of Lemma~\ref{lem:reverse-model-estimates}
		therefore gives
		$$
		\|-2\langle\nabla\chi_j,\nabla b\rangle+(\Delta\chi_j)b\|_{\mathfrak A_r}
		\leq C\bigl(\|F\|_1+r^{-2}\|B\|_1\bigr).
		$$
		Combine this with \eqref{eq:reverse-local-commutator},
		\eqref{eq:reverse-local-F-bound}, and
		\eqref{eq:reverse-boundary-atoms}. Since $r_j\approx\rho$
		on $2B_j$, bounded overlap yields
		$$
		\begin{aligned}
			\sum_j\|a_j\|_{\mathfrak A_{r_j}}
			&\leq C\bigl(\|A\|_1+\|g\|_1+\|f\|_1
			+\|\rho^{-1}\nabla b\|_1+\|\rho^{-2}b\|_1\bigr),\\
			\sum_j r_j^{-2}\|h_j\|_1&\leq C\|\rho^{-2}b\|_1.
		\end{aligned}
		$$
		Now apply \eqref{eq:reverse-local-known} and
		\eqref{eq:reverse-local-L1}.
	\end{proof}
	
	\subsubsection{A truncated square function on the complete models}
	
	The finite-propagation construction below uses one further model
	estimate. Let $\Phi$ be an even Schwartz function with compactly
	supported Fourier transform, let $J\geq1$, and set
	$Q(s)=s^{2J}\Phi(s)$. Then
	\begin{equation}\label{eq:reverse-truncated-model-square}
		\left\|\left(\int_0^r\int_{d(x,y)<t}
		|Q(t\sqrt{\Delta_X})F(y)|^2
		\frac{d\nu(y)\,dt}{t\nu(B(y,t))}\right)^{1/2}\right\|_1
		\leq C\|F\|_{\mathfrak A_r}.
	\end{equation}
	The constant is uniform for the models and radii used above.
	Here $J$ and $\Phi$ are fixed.
	
	We prove this estimate on the auxiliary atoms defined above. Finite propagation places the kernel of
	$Q(t\sqrt{\Delta_X})$ in $d(x,y)\leq C_\Phi t$.
	The model heat kernels have Gaussian H\"older bounds for every
	positive time on the complete ends and for times at most one
	on the finitely many central models. These bounds follow from
	the product heat kernels on the ends and uniform ellipticity
	on the central models. Take $\nu_0$ small enough that $r_j\leq1$
	for every central ball.
	Since $Q$ is even and Schwartz,
	\cite[Theorem~3.4, equations~(3.1)--(3.2)]{KP}
	gives both the size and H\"older estimates. On an end model
	apply that theorem to $t^2\Delta_X$ with distance $d/t$ and
	parameter one. Its heat kernel is $p_{t^2s}$, so its constants
	are independent of $t$. On a central model apply the theorem
	for $0<t\leq r_j$. Finite propagation and model doubling give
	$$
	|K_t(y,z)|\leq\frac{C}{\nu(B(y,t))}.
	$$
	For some $\gamma>0$, we also have
	$$
	|K_t(y,z)-K_t(y,z')|
	\leq\frac{C}{\nu(B(y,t))}
	\left(\frac{d(z,z')}{t}\right)^\gamma
	\quad\text{if }d(z,z')\leq t,
	$$
	with the same finite support restriction up to a fixed enlargement.
	The $L^2$ square estimate follows from
	$\int_0^\infty|Q(t\lambda)|^2\,dt/t<\infty$.
	Let $a$ be an atom in $B=B(x_B,s)$, with $s<r$.
	The $L^2$ estimate and doubling give
	$$
	\int_{CB}\mathcal A_r(Qa)\,d\nu
	\leq\nu(CB)^{1/2}\|\mathcal A_r(Qa)\|_2\leq C,
	$$
	where $\mathcal A_r$ denotes the square function in
	\eqref{eq:reverse-truncated-model-square}.
	For $R=d(x,x_B)\geq C_1s$, finite propagation shows that only
	$t\geq\kappa R$ can contribute, where $\kappa>0$ is fixed.
	Choose $C_1$ so that these times also satisfy $t\geq2s$.
	Cancellation gives
	$$
	Q(t\sqrt{\Delta_X})a(y)
	=\int_B[K_t(y,z)-K_t(y,x_B)]a(z)\,d\nu(z).
	$$
	Since $\|a\|_1\leq1$ and $d(x,y)<t$, the kernel estimate and
	doubling imply
	$$
	|Q(t\sqrt{\Delta_X})a(y)|
	\leq\frac{C}{\nu(B(x,t))}\left(\frac{s}{t}\right)^\gamma.
	$$
	Consequently,
	$$
	\begin{aligned}
		\mathcal A_r(Qa)(x)^2
		&\leq C\int_{\kappa R}^r
		\frac{(s/t)^{2\gamma}}{\nu(B(x,t))^2}\,\frac{dt}{t}
		\leq\frac{C}{\nu(B(x,R))^2}\left(\frac{s}{R}\right)^{2\gamma}.
	\end{aligned}
	$$
	The integral is zero if $\kappa R\geq r$.
	On $2^ks\leq R<2^{k+1}s$, doubling compares
	$\nu(B(x,R))$ with $\nu(B(x_B,2^ks))$. Therefore
	$$
	\int_{2^{k+1}B\setminus2^kB}\mathcal A_r(Qa)\,d\nu
	\leq C2^{-k\gamma}.
	$$
	The sum converges. If the atom has radius $r$, finite propagation
	places its truncated square function in a fixed enlargement
	of its ball; the $L^2$ estimate alone then suffices. Summing
	an atomic representation gives the assertion for finite sums.
	For an $L^1$ atomic sum, the model multiplier kernels give
	locally uniform convergence at every fixed $t>0$. Fatou gives the same
	bound for that sum and proves
	\eqref{eq:reverse-truncated-model-square}.
	
	\subsection{The completion in \texorpdfstring{$L^1$}{L1}}
	
	We use the global Hardy space defined in \eqref{eq:H1}. Its area function is
	$\Sc f=S_c(P_tf)$ of \eqref{eq:restricted-operators}. The time interval in its
	definition is $(0,\infty)$. By Theorem \ref{thm:main},
	$$
	\|\Sc f\|_1=0\quad\Longrightarrow\quad
	\Sc f=0\text{ almost everywhere}\quad\Longrightarrow\quad f=0
	\qquad(f\in L^2).
	$$
	Homogeneity and the triangle inequality follow from the
	linearity of $P_t$ and Minkowski's inequality on each cone.
	Thus $\|\Sc(\cdot)\|_1$ is a norm on the space being completed.
	
	Equation \eqref{eq:reverse-local-L1}, proved above, gives
	$$
	\|f\|_1\leq C\|\Sc f\|_1
	\qquad(f\in L^2,\ \Sc f\in L^1).
	$$
	We justify its use on the completion. Suppose that $(f_j)$ is
	Cauchy in this norm. The last inequality gives $f_j\to f$ in
	$L^1$ for some $f$. For every compact
	$E\subset M\times(0,\infty)$, the Poisson kernel bounds in
	\cite{BDLW} and the local harmonic gradient estimate give
	$$
	\sup_{(y,t)\in E}
	\bigl(|P_th(y)|+|\nabla_{y,t}P_th(y)|\bigr)
	\leq C_E\|h\|_1\qquad(h\in L^1).
	$$
	Indeed, the kernel bound proves the estimate for $P_th$ on a
	slightly larger compact set. The interior gradient estimate on
	finitely many product balls proves the derivative estimate.
	The $L^1$ and $L^2$ definitions of $P_t$ agree on their intersection.
	Consequently, the gradients of $P_tf_k$ converge locally to
	those of $P_tf$. Fatou's lemma, first on the cones and then on
	$M$, gives
	$$
	\|\Sc(f-f_j)\|_1
	\leq\liminf_{k\to\infty}\|\Sc(f_k-f_j)\|_1
	\longrightarrow0\qquad(j\to\infty).
	$$
	If $f=0$, then $\|\Sc f_j\|_1\to0$. The map from the completion
	to $L^1$ is therefore injective. We identify the completion with
	its image; its norm is $\|\Sc f\|_1$.
	\subsection{Atomic decomposition in the global Hardy space}
	\label{sec:reverse-atoms}
	
	\begin{proposition}\label{thm:restricted-atomic-reverse}
		If $f\in L^2(M)$ and $\Sc f\in L^1(M)$, then
		\begin{equation}\label{eq:note-main-series}
			f=\sum_\nu\lambda_\nu a_\nu,\qquad
			\sum_\nu|\lambda_\nu|\leq C_{c,\delta}\|\Sc f\|_1,
		\end{equation}
		where the $a_\nu$ satisfy \eqref{eq:atoms}. The series converges
		in $L^1(M)$ and in $ H^1_{S_c}(M) $. Conversely, every $L^1$ atomic
		sum satisfies
		\begin{equation}\label{eq:note-main-forward}
			\|\Sc f\|_1\leq C_{c,\delta}\sum_\nu|\lambda_\nu|.
		\end{equation}
	\end{proposition}
	
	\begin{proof}[Proof of Proposition~\ref{thm:restricted-atomic-reverse}]
		Let $b$ be the distributional solution $\Delta b=f$ constructed above.
		Choose a locally finite partition of unity $\{\chi_j\}$ on
		admissible balls $B_j=B(x_j,r_j)$, with $r_j=\nu_0\rho(x_j)$,
		where $\nu_0>0$ is the small constant chosen above. The preceding weighted
		Green estimate and the estimates on the complete models give
		\begin{equation}\label{eq:atomic-local-sum}
			\|f\|_1+\|\rho^{-2}b\|_1
			+\sum_j\|f_j\|_{\mathfrak A_{r_j}(X_j)}
			\leq C_{c,\delta}\|\Sc f\|_1,
			\qquad b_j=\chi_jb,\quad f_j=\Delta b_j.
		\end{equation}
		Here $\Delta_j$ is the Laplacian on a complete model agreeing with
		$\Delta$ on a fixed sufficiently large dilation of $B_j$.
		On the ends these are $\mathbb R^m$ and
		$\mathbb R^n\times\mathbb S^{m-n}$. For the finitely many balls
		near $K$, use the complete models from the local construction.
		The notation $\mathfrak A_r(X_j)$ denotes the atomic space defined in
		Section~\ref{sec:complete-model-estimates}. All constants in its
		truncated tent estimate are uniform in $j$.
		
		$b_j\in\operatorname{Dom}(\Delta_j)$ and $f_j\in L^2$.
		On the common region the operators agree, so $f_j=\Delta_jb_j$.
		Moreover,
		\begin{equation}\label{eq:atomic-potential-sum}
			\sum_j r_j^{-2}\|b_j\|_1\leq C\|\rho^{-2}b\|_1,
			\qquad \sum_j f_j=\Delta\Big(\sum_j\chi_jb\Big)=f
			\quad\hbox{in distributions}.
		\end{equation}
		We now decompose each $f_j$ on its complete model. All supports
		obtained below remain in the region where that model agrees
		with $M$.
		
		Choose a real even Schwartz function $\Phi$ with
		$\operatorname{supp}\widehat\Phi\subset[-1,1]$, and put
		$$
		\psi(s)=s^2\Phi(s),\qquad
		\int_0^\infty\psi(s)^2\,\frac{ds}{s}=1.
		$$
		A constant multiple of any nonzero such $\Phi$ gives this
		normalization. Define
		\begin{equation}\label{eq:atomic-low-multiplier}
			\eta(s)=1-\int_0^s\psi(v)^2\,\frac{dv}{v}.
		\end{equation}
		Since $\psi$ is even and vanishes to order at least two at zero,
		the integral defines an even entire function. On the real axis,
		$$
		\eta(s)=\int_{|s|}^\infty\psi(v)^2\,\frac{dv}{v}.
		$$
		Thus $\eta$ is Schwartz. Moreover,
		$$
		\eta'(s)=-s^3\Phi(s)^2,\qquad
		\operatorname{supp}\widehat{\eta'}\subset[-2,2].
		$$
		The identity $i\xi\widehat\eta(\xi)=\widehat{\eta'}(\xi)$
		therefore gives $\operatorname{supp}\widehat\eta\subset[-2,2]$.
		Finite propagation therefore gives
		\begin{equation}\label{eq:atomic-propagation}
			\begin{split}
				\operatorname{supp}K_{\Phi(t\sqrt{\Delta_j})}\cup
				\operatorname{supp}K_{\psi(t\sqrt{\Delta_j})}
				&\subset\{(x,y):d_j(x,y)\leq t\},\\
				\operatorname{supp}K_{\eta(r\sqrt{\Delta_j})}
				&\subset\{(x,y):d_j(x,y)\leq2r\}.
			\end{split}
		\end{equation}
		The same support assertion holds for $r^2\Delta_j\eta(r\sqrt{\Delta_j})$.
		
		The spectral identity
		$1=\int_0^r\psi(t\lambda)^2\,dt/t+\eta(r\lambda)$ gives
		\begin{equation}\label{eq:atomic-inhomogeneous-formula}
			f_j=\int_0^{r_j}\psi(t\sqrt{\Delta_j})
			\bigl[\psi(t\sqrt{\Delta_j})f_j\bigr]\,\frac{dt}{t}
			+\Delta_j\eta(r_j\sqrt{\Delta_j})b_j
			\quad\hbox{in }L^2.
		\end{equation}
		The last term is $\Delta_j\beta$ with
		$\beta=\eta(r_j\sqrt{\Delta_j})b_j$. It will also give an
		$\Delta$-atom satisfying \eqref{eq:atoms-harmonic-cancellation}.
		
		\medskip
		\noindent\emph{The integral over $0<t<r_j$.}
		Set
		$$
		F_j(y,t)=\mathbf1_{(0,r_j)}(t)\psi(t\sqrt{\Delta_j})f_j(y).
		$$
		Let $T_j^p$ be the usual tent space on the complete model, with
		$$
		\mathcal A_jF(x)^2
		=\int_0^\infty\int_{B_j^{\rm mod}(x,t)}
		|F(y,t)|^2\frac{d\mu_j(y)\,dt}{tV_j(y,t)},
		\qquad \|F\|_{T_j^p}=\|\mathcal A_jF\|_p.
		$$
		Estimate \eqref{eq:reverse-truncated-model-square} with $J=1$
		and the spectral theorem give
		\begin{equation}\label{eq:atomic-local-tent-estimate}
			\|F_j\|_{T_j^1}\leq C\|f_j\|_{\mathfrak A_{r_j}(X_j)},
			\qquad
			\int_0^{r_j}\|F_j(\cdot,t)\|_2^2\,\frac{dt}{t}
			\leq\|f_j\|_2^2.
		\end{equation}
		By \eqref{eq:atomic-propagation},
		$\operatorname{supp}F_j\subset 2B_j\times(0,r_j)$.
		Lemma \ref{lem:reverse-model-tent} gives
		\begin{equation}\label{eq:atomic-tent-decomposition}
			F_j=\sum_k\lambda_{jk}A_{jk},\qquad
			\sum_k|\lambda_{jk}|\leq C\|F_j\|_{T_j^1}.
		\end{equation}
		Here each $A=A_{jk}$ has a ball $B=B_{jk}$ such that
		\begin{equation}\label{eq:atomic-tent-size}
			\operatorname{supp}A\subset B\times(0,r_B),\qquad
			\int_0^{r_B}\|A(\cdot,t)\|_2^2\,\frac{dt}{t}
			\leq\mu_j(B)^{-1}.
		\end{equation}
		The support assertion in that lemma gives
		$r_B\leq C_0r_j$ and $B\subset C_0B_j$ for a fixed $C_0$.
		Its proof also gives convergence in $T_j^1$ and $T_j^2$.
		Thus every ball and every finite-propagation kernel below remains
		inside the chosen common region.
		
		For an atom $A$ in \eqref{eq:atomic-tent-size}, put
		\begin{equation}\label{eq:atomic-tent-synthesis}
			a_A=\int_0^{r_B}\psi(t\sqrt{\Delta_j})A(\cdot,t)\,\frac{dt}{t},
			\qquad
			\beta_A=\int_0^{r_B}t^2\Phi(t\sqrt{\Delta_j})A(\cdot,t)\,\frac{dt}{t}.
		\end{equation}
		For $h\in L^2$, Cauchy--Schwarz and the spectral theorem give
		$$
		|\langle a_A,h\rangle|
		\leq\left(\int_0^{r_B}\|A(\cdot,t)\|_2^2\,\frac{dt}{t}\right)^{1/2}
		\left(\int_0^\infty
		\|\psi(t\sqrt{\Delta_j})h\|_2^2\,\frac{dt}{t}\right)^{1/2}
		\leq\mu_j(B)^{-1/2}\|h\|_2.
		$$
		Also,
		$$
		\|\beta_A\|_2
		\leq\|\Phi\|_\infty
		\left(\int_0^{r_B}t^4\,\frac{dt}{t}\right)^{1/2}
		\left(\int_0^{r_B}\|A(\cdot,t)\|_2^2\,\frac{dt}{t}\right)^{1/2}
		\leq\tfrac12\|\Phi\|_\infty r_B^2\mu_j(B)^{-1/2}.
		$$
		Thus
		\begin{equation}\label{eq:atomic-synthesis-bounds}
			\|a_A\|_2\leq\mu_j(B)^{-1/2},\qquad
			\|\beta_A\|_2\leq C r_B^2\mu_j(B)^{-1/2}.
		\end{equation}
		First integrate over $\varepsilon<t<r_B$. Then
		$\Delta_j\beta_{A,\varepsilon}=a_{A,\varepsilon}$.
		The two integrals converge in $L^2$ as $\varepsilon\downarrow0$:
		for $\beta_A$ use the preceding estimate, and for $a_A$ use the
		boundedness of the synthesis map from $T_j^2$ to $L^2$.
		Closedness of $\Delta_j$ gives
		$\beta_A\in\operatorname{Dom}(\Delta_j)$ and $a_A=\Delta_j\beta_A$.
		Finite propagation gives
		$$
		\operatorname{supp}\beta_A\cup\operatorname{supp}a_A
		\subset\overline{2B}.
		$$
		Model doubling gives $\mu_j(2B)\leq C\mu_j(B)$.
		After division by one fixed constant,
		\eqref{eq:atomic-synthesis-bounds} therefore has exactly the
		normalization \eqref{eq:atoms} on a slightly larger ball.
		
		Choose $\nu_0$ so small that every slightly enlarged ball $2B$
		obtained here has radius at most $\delta\rho$ at its centre and lies in the common region of the model and $M$.
		Then $\Delta_j\beta_A=\Delta\beta_A$, and the normalized $a_A$ is an
		$\Delta$-atom satisfying \eqref{eq:atoms}. Equations
		\eqref{eq:atomic-local-tent-estimate}--\eqref{eq:atomic-tent-decomposition}
		and the $T_j^2$ synthesis estimate show that the first term in
		\eqref{eq:atomic-inhomogeneous-formula} is a sum of such atoms,
		with coefficient sum at most
		$C\|f_j\|_{\mathfrak A_{r_j}(X_j)}$.
		
		\medskip
		\noindent\emph{The remaining term.}
		Put $\beta_j^0=\eta(r_j\sqrt{\Delta_j})b_j$.
		The functions $\eta(s)$ and $s^2\eta(s)$ are even and Schwartz.
		The size estimate in
		\cite[Theorem~3.4, equation~(3.1)]{KP}, applied on the models
		as above, gives
		\begin{equation}\label{eq:atomic-low-kernel-bound}
			\|\beta_j^0\|_2+r_j^2\|\Delta_j\beta_j^0\|_2
			\leq C\mu_j(B_j)^{-1/2}\|b_j\|_1.
		\end{equation}
		Finite propagation makes both kernels vanish when
		$d_j(x,y)>2r_j$. The size estimate and model doubling give
		$|K(x,y)|\leq C/V_j(y,r_j)$. For $y\in B_j$, model doubling
		then gives
		$$
		\|K(\cdot,y)\|_2\leq C\mu_j(B_j)^{-1/2},
		$$
		which proves \eqref{eq:atomic-low-kernel-bound} by Minkowski.
		Furthermore,
		$$
		\operatorname{supp}\beta_j^0\cup
		\operatorname{supp}\Delta_j\beta_j^0\subset\overline{3B_j}.
		$$
		Thus $\Delta_j\beta_j^0=\Delta\beta_j^0$ is a multiple of one normalized
		$\Delta$-atom satisfying \eqref{eq:atoms}, with coefficient at most
		$C r_j^{-2}\|b_j\|_1$. Use a slightly enlarged $3B_j$ to
		contain the closed support in its interior. Terms with zero
		coefficient are omitted. By \eqref{eq:atomic-potential-sum} and
		\eqref{eq:reverse-green-bound}, the total coefficient of these
		remaining terms is at most $C\|\rho^{-2}b\|_1\leq C\|\Sc f\|_1$.
		
		Combining the two parts of
		\eqref{eq:atomic-inhomogeneous-formula}, and then summing in $j$, gives
		$$
		\sum_\nu|\lambda_\nu|
		\leq C\sum_j\left(
		\|f_j\|_{\mathfrak A_{r_j}(X_j)}
		+r_j^{-2}\|b_j\|_1\right)
		\leq C_{c,\delta}\|\Sc f\|_1.
		$$
		For each $j$, the $L^2$ identity
		\eqref{eq:atomic-inhomogeneous-formula} identifies its atomic sum
		with $f_j$. The normalized atoms have $L^1$ norm at most one,
		so the complete series converges absolutely in $L^1(M)$.
		Its distributional sum is $f$ by
		\eqref{eq:atomic-potential-sum}. This proves the $L^1$ assertion
		in \eqref{eq:note-main-series}.
		
		Finally, write $g_N=\sum_{\nu\leq N}\lambda_\nu a_\nu$.
		For $N'>N$, the forward atomic estimate gives
		$$
		\|\Sc(g_{N'}-g_N)\|_1
		\leq C_{c,\delta}\sum_{N<\nu\leq N'}|\lambda_\nu|.
		$$
		Since $g_{N'}\to f$ in $L^1$, their Poisson extensions and gradients
		converge locally for positive times, by the kernel bounds and
		interior harmonic estimates. Fatou on the cones and on $M$ gives
		$$
		\|\Sc(f-g_N)\|_1
		\leq C_{c,\delta}\sum_{\nu>N}|\lambda_\nu|\longrightarrow0.
		$$
		This proves the reverse implication and its asserted convergence.
		The same argument with finite partial sums, the forward atomic
		estimate, and Fatou proves \eqref{eq:note-main-forward} for every
		$L^1$ atomic sum. Finite sums of atoms belong to $L^2$.
		Their density in the atomic norm, and the two estimates just
		proved, identify $H^1_{S_c}(M)$ with
		$H^1_{\mathrm{at},\Delta,\delta}(M)$ and give
		\eqref{eq:note-global-norm}.
	\end{proof}
	
	\subsection{The atomic description of the completion}
	
	We now prove Theorem \ref{thm:hardy-characterization}(i) on the
	whole completion. For $f\in H^1_{S_c}(M) $, $f\ne0$, choose functions
	$f_j\in L^2$ with $\Sc f_j\in L^1$ and
	$$
	\|f-f_j\|_{ H^1_{S_c}(M) }\leq2^{-j}\|f\|_{ H^1_{S_c}(M) }.
	$$
	Then
	$$
	f=f_1+\sum_{j\geq1}(f_{j+1}-f_j),\qquad
	\|f_1\|_{ H^1_{S_c}(M) }
	+\sum_{j\geq1}\|f_{j+1}-f_j\|_{ H^1_{S_c}(M) }
	\leq3\|f\|_{ H^1_{S_c}(M) }.
	$$
	Apply Proposition \ref{thm:restricted-atomic-reverse} to $f_1$
	and each difference. The combined coefficient sum is at most
	$C\|f\|_{ H^1_{S_c}(M) }$. The resulting series converges in
	$ H^1_{S_c}(M) $ and $L^1$ and has sum $f$ by the injective
	embedding proved above. Conversely, \eqref{eq:note-main-forward}
	shows that the finite partial sums of any $L^1$ atomic series
	are Cauchy in $ H^1_{S_c}(M) $. Their limit has the given $L^1$ sum.
	Thus
	\begin{equation}\label{eq:prior-H1}
		H^1_{S_c}(M) \hookrightarrow L^1(M),\qquad
		H^1_{S_c}(M) =H^1_{\mathrm{at},\Delta,\delta}(M)
		\quad\text{with equivalent norms},
	\end{equation}
	and
	\begin{equation}\label{eq:note-global-norm}
		\|f\|_{ H^1_{S_c}(M) }\approx
		\inf_{f=\sum_\nu\lambda_\nu a_\nu}\sum_\nu|\lambda_\nu|.
	\end{equation}
	This gives the series, both convergence assertions, and the norm
	equivalence in \eqref{eq:intro-atomic}. It also gives
	$$
	\int_M f\overline h\dd\mu=0
	\qquad(f\in H^1_{S_c}(M) , h\text{ bounded},\ \Delta h=0)
	$$
	by \eqref{eq:atoms-harmonic-cancellation}.
	Both $ H^1_{S_c}(M) $ and $L^2$ are continuously and injectively
	embedded in $L^1+L^2$. Hence $( H^1_{S_c}(M) ,L^2)$ is a compatible
	pair of Banach spaces.
	
	\subsection{The maximal characterization}
	\label{subsec:Hardy-maximal}
	
	\begin{proof}[Proof of Theorem \ref{thm:hardy-characterization}(ii)]
		Let $a=\Delta b$ be an atom in \eqref{eq:atoms}. Local doubling and
		the $L^2$ maximal bound give
		$$
		\int_{2B}\mathcal N_ca\dd\mu
		\leq\mu(2B)^{1/2}\|\mathcal N_ca\|_2
		\leq C_c\mu(B)^{1/2}\|a\|_2\leq C_c.
		$$
		Since $b\in\operatorname{Dom}(\Delta)$, the spectral theorem gives
		$P_ta=\Delta P_tb$. The kernel estimate \eqref{eq:prior-kernel} and
		$\mathcal N_c\leq\mathcal N_c^\#$ therefore give
		\begin{align*}
			\int_{M\setminus2B}\mathcal N_ca\dd\mu
			&\leq\int_B|b(z)|
			\int_{M\setminus2B}
			\sup_{(y,t)\in\Gamma_c(x)}|(\Delta P_t)(y,z)|
			\dd\mu(x)\dd\mu(z)
			\leq Cr^{-2}\|b\|_1\leq C.
		\end{align*}
		Thus every atom satisfies
		\begin{equation}\label{eq:Hardy-atom-maximal}
			\|\mathcal N_ca\|_1\leq C_{c,\delta}.
		\end{equation}
		If $f=\sum_\nu\lambda_\nu a_\nu$ is its atomic representation,
		the partial sums converge in $L^1$. The local Poisson estimate
		used above for the completion gives convergence of their
		extensions at every fixed $(y,t)$. Hence
		$$
		\mathcal N_cf(x)
		\leq\sum_\nu|\lambda_\nu|\mathcal N_ca_\nu(x),
		\qquad
		\|\mathcal N_cf\|_1
		\leq C_{c,\delta}\sum_\nu|\lambda_\nu|.
		$$
		Taking the infimum proves
		$$
		\|\mathcal N_cf\|_1\leq C_{c,\delta}\|f\|_{ H^1_{S_c}(M) },
		\qquad f\in H^1_{S_c}(M) .
		$$
		
		Conversely, let $f\in L^1(M)$ and $\mathcal N_cf\in L^1(M)$.
		We prove convergence of its Poisson approximations in the
		maximal norm. The full Poisson kernel estimates in
		\cite[Theorem~2.2, with $k=0$]{BDLW} give
		$$
		0\leq P_t(y,z)\leq C(t^{-m}+t^{-n}),\qquad y,z\in M.
		$$
		Together with $L^1$ contractivity, this gives
		\begin{equation}\label{eq:Hardy-Poisson-regularization}
			\|P_\varepsilon f\|_2^2
			\leq\|P_\varepsilon f\|_\infty\|P_\varepsilon f\|_1
			\leq C(\varepsilon^{-m}+\varepsilon^{-n})\|f\|_1^2,
		\end{equation}
		and
		$$
		\sup_y|P_tf(y)|
		\leq C(t^{-m}+t^{-n})\|f\|_1\longrightarrow0
		\qquad(t\to\infty).
		$$
		Corollary~\ref{cor:restricted-boundary} gives
		\eqref{eq:Hardy-boundary-limit}.

		Put $f_\varepsilon=P_\varepsilon f$. If $(y,t)\in\Gamma_c(x)$,
		then $(y,t+\varepsilon)\in\Gamma_c(x)$. The semigroup identity
		therefore gives
		$$
		\mathcal N_cf_\varepsilon(x)\leq\mathcal N_cf(x),\qquad
		\mathcal N_c(f_\varepsilon-f)(x)\leq2\mathcal N_cf(x).
		$$
		Fix an $x$ for which \eqref{eq:Hardy-boundary-limit} holds.
		For $0<t<\eta$ and $0<\varepsilon<\eta$, that boundary limit
		bounds both $P_tf(y)-f(x)$ and $P_{t+\varepsilon}f(y)-f(x)$
		uniformly on the cone, by quantities tending to zero as
		$\eta\downarrow0$. On $\eta\leq t\leq T$, local smoothness
		gives uniform continuity on
		$$
		\overline{B(x,c\rho(x))}\times[\eta,T+1].
		$$
		This set is compact because $M$ is complete. For $t>T$, the
		global decay estimate above bounds both terms by a quantity
		tending to zero as $T\to\infty$. First choose $\eta$ and $T$,
		and then let $\varepsilon\downarrow0$. It follows that
		$$
		\mathcal N_c(f_\varepsilon-f)(x)\longrightarrow0
		\quad\text{for almost every }x.
		$$
		Dominated convergence proves
		\begin{equation}\label{eq:Hardy-maximal-approximation}
			\|\mathcal N_c(f_\varepsilon-f)\|_1\longrightarrow0.
		\end{equation}
		The boundary limit gives $|f|\leq\mathcal N_cf$ almost
		everywhere. The same domination also proves
		$f_\varepsilon\to f$ in $L^1$.
		
		By \eqref{eq:Hardy-Poisson-regularization}, each
		$f_\varepsilon$ belongs to $L^2$. Lemma
		\ref{lem:Hardy-FS-L2} gives
		$$
		\|\Sc f_\varepsilon\|_1\leq C_c\|\mathcal N_cf\|_1,
		$$
		and, using \eqref{eq:Hardy-maximal-approximation},
		\begin{align*}
			\|\Sc(f_\varepsilon-f_{\varepsilon'})\|_1
			&\leq C_c\|\mathcal N_c(f_\varepsilon-f_{\varepsilon'})\|_1
			\leq C_c\bigl(\|\mathcal N_c(f_\varepsilon-f)\|_1
			+\|\mathcal N_c(f_{\varepsilon'}-f)\|_1\bigr)
			\longrightarrow0.
		\end{align*}
		Thus these approximations converge in $ H^1_{S_c}(M) $.
		Its injective $L^1$ embedding identifies their limit with $f$.
		In particular,
		$$
		f\in H^1_{S_c}(M) ,\qquad
		\|f\|_{ H^1_{S_c}(M) }\leq C_c\|\mathcal N_cf\|_1.
		$$
		This proves both assertions in \eqref{eq:intro-maximal-Hardy}.
	\end{proof}
	
\subsection{The reverse \texorpdfstring{$L^p$}{Lp} estimate}
\label{subsec:reverse-Lp}

\begin{theorem}\label{thm:reverse-Lp}
For $1<p\leq2$ and $0<c\leq1/100$,
\begin{equation}\label{eq:reverse-Lp}
 C_{p,c}^{-1}\|f\|_p\leq\|\Sc f\|_p
 \leq C_{p,c}\|f\|_p,\qquad f\in L^p(M).
\end{equation}
\end{theorem}

Choose a partition $\{\chi_j\}$ and complete models $X_j$ as in
Section~\ref{sec:complete-model-estimates}, with $\nu_0$ fixed
using only $M$ and $c$. Their tent norms are
those of Lemma~\ref{lem:reverse-model-tent}. Let $\mathcal M_j$
denote the Hardy--Littlewood maximal operator on $X_j$.

\begin{lemma}\label{lem:lp-model-test}
For $2<q<\infty$, $r=r_j$, and
$v=e^{-t\sqrt{\Delta_j}}\varphi$, one has
\begin{equation}\label{eq:lp-model-test}
 \left\|t\nabla_{x,t}v\,
 \mathbf1_{B_j\times(0,2r)}\right\|_{T^q}
 \leq C_q\|\varphi\|_{L^q(X_j)}.
\end{equation}
Moreover,
\begin{equation}\label{eq:lp-model-maximal}
 \sup_{0<t<2r}\bigl(|v(x,t)|+t|\nabla_{x,t}v(x,t)|\bigr)
 \leq C\mathcal M_j\varphi(x),
\end{equation}
where $\mathcal M_j\varphi=\mathcal M_j(|\varphi|)$.
The constants are uniform in $j$.
\end{lemma}

\begin{proof}
 Put
 $$
 T_j\varphi=t\nabla_{x,t}e^{-t\sqrt{\Delta_j}}\varphi
 \,\mathbf1_{B_j\times(0,2r)}.
 $$
 Tonelli's theorem and the spectral theorem give
 \begin{equation}\label{eq:lp-model-L2}
 \begin{aligned}
 \|T_j\varphi\|_{T^2}^2
 &=\int_0^{2r}\int_{B_j}
 t|\nabla_{x,t}v|^2\,d\mu_j\,dt
 \leq\int_0^\infty t
 \left(\|\Delta_j^{1/2}v\|_2^2
 +\|\partial_tv\|_2^2\right)dt
 \leq\tfrac12\|\varphi\|_2^2.
 \end{aligned}
 \end{equation}
 If $\varphi\in L^\infty(X_j)$, then
 $\|\varphi\|_{\operatorname{bmo}_r}\leq3\|\varphi\|_\infty$.
 Equation~\eqref{eq:reverse-bmo-carleson} therefore gives
 $$
 \int_0^{s(B)}\int_B|T_j\varphi|^2\,d\mu_j\,\frac{dt}{t}
 \leq C\mu_j(B)\|\varphi\|_\infty^2
 \qquad(s(B)\leq4r).
 $$
 If $s(B)>4r$ and $B\cap B_j\ne\varnothing$, then
 $\mu_j(B_j)\leq C\mu_j(B)$, and
 $$
 \begin{aligned}
 \int_0^{s(B)}\int_B|T_j\varphi|^2\,d\mu_j\,\frac{dt}{t}
 &\leq\int_0^{2r}\int_{B_j}t|\nabla_{x,t}v|^2\,d\mu_j\,dt
 \leq C\mu_j(B_j)\|\varphi\|_\infty^2
 \leq C\mu_j(B)\|\varphi\|_\infty^2.
 \end{aligned}
 $$
 If $B\cap B_j=\varnothing$, the integral is zero. Thus
 \begin{equation}\label{eq:lp-model-Linfty}
 \|T_j\varphi\|_{T^\infty}\leq C\|\varphi\|_\infty.
 \end{equation}
 On a doubling model the box norm in
 Lemma~\ref{lem:reverse-model-tent} is equivalent to the Carleson
 norm defined using tents. The complex interpolation identity
  gives
 $$
 [T^2,T^\infty]_{1-2/q}=T^q,
 \qquad [L^2,L^\infty]_{1-2/q}=L^q.
 $$
 Gradients have finitely many Euclidean components on each model;
 on the sphere, use its embedding in Euclidean space. Apply the
 scalar interpolation identity to each component.
 Interpolating \eqref{eq:lp-model-L2} and
 \eqref{eq:lp-model-Linfty} proves \eqref{eq:lp-model-test}.

 Write $v_0=e^{-t\sqrt{\Delta_j}}|\varphi|$. On either end
 model, subordination of the Gaussian heat kernel bound gives
 $$
 e^{-t\sqrt{\Delta_j}}(x,z)
 \leq\frac{Ct}{(t+d_j(x,z))\mu_j(B(x,t+d_j(x,z)))}.
 $$
 Integration over $B(x,t)$ and the annuli
 $2^kt\leq d_j(x,z)<2^{k+1}t$ yields
 \begin{equation}\label{eq:lp-model-positive-maximal}
 v_0(x,t)\leq C\left(1+\sum_{k\geq0}2^{-k}\right)
 \mathcal M_j\varphi(x)
 \leq C\mathcal M_j\varphi(x).
 \end{equation}
 For a central model, split the heat integral at heat time one,
 as in Lemma~\ref{lem:reverse-model-estimates}. The first part has
 the same bound. The remaining Poisson kernel is at most
 $Ct(1+d_j(x,z))^{-m-1}$. Since $0<t<2r_j$ and the central
 radii are bounded,
 $$
 t\int_{X_j}(1+d_j(x,z))^{-m-1}|\varphi(z)|\,d\mu_j(z)
 \leq Ct\left(1+\sum_{k\geq0}2^{-k}\right)
 \mathcal M_j\varphi(x)
 \leq C\mathcal M_j\varphi(x).
 $$
 This proves \eqref{eq:lp-model-positive-maximal} on every model.

 The end models have nonnegative Ricci curvature. Each central
 model has a fixed lower Ricci bound, and $t<2r_j$ is bounded.
 The Cheng--Yau estimate \cite{CY}, applied on the product ball
 of radius $t/2$ about $(x,t)$, gives
 $$
 t|\nabla_{x,t}h(x,t)|\leq Ch(x,t)
 $$
 for every positive harmonic $h$ on the model product.
 Apply this to
 $$
 v_0+\varepsilon+\operatorname{Re}v,\quad
 v_0+\varepsilon-\operatorname{Re}v,\quad
 v_0+\varepsilon+\operatorname{Im}v,\quad
 v_0+\varepsilon-\operatorname{Im}v.
 $$
 These functions are positive because $|v|\leq v_0$.
 Subtracting the two real parts and the two imaginary parts gives
 $$
 t|\nabla_{x,t}v(x,t)|\leq C(v_0(x,t)+\varepsilon).
 $$
 Let $\varepsilon\downarrow0$ and use
 \eqref{eq:lp-model-positive-maximal}. This proves
 \eqref{eq:lp-model-maximal}.
\end{proof}

\begin{proof}[Proof of Theorem~\ref{thm:reverse-Lp}]
 The upper estimate is already proved. The lower estimate for
 $p=2$ is Theorem~\ref{thm:main}. Let $1<p<2$, put
 $q=p/(p-1)>2$, and first take $f\in L^p(M)\cap L^2(M)$.
 Write
 $$
 u(x,t)=P_tf(x),\qquad A(x)=\Sc f(x),\qquad
 g(x)=u(x,H(x)),\quad H(x)=\varepsilon R(x).
 $$
 The pointwise estimates \eqref{eq:reverse-local-gradient} and
 \eqref{eq:reverse-trace-difference} give
 \begin{equation}\label{eq:lp-trace-pointwise}
 \rho(x)|\nabla g(x)|\leq CA(x),\qquad
 |u(x,t)-g(x)|\leq CA(x)|\log(t/H(x))|\quad(0<t\leq2H(x)).
 \end{equation}
 The decay estimate \eqref{eq:decay} gives $g(x)\to0$ uniformly
 along both ends. Choose $r_1$ beyond the compact part, as in
 Lemma~\ref{lem:reverse-trace}, and let $E_d$ denote the part
 $r>r_1$ of the end of dimension $d=m$ or $n$. Then
 $$
 g(r,\omega)=-\int_1^\infty sr\,\partial_rg(sr,\omega)\,
 \frac{ds}{s}.
 $$
 For almost every $\omega$, H\"older's inequality gives
 $$
 \int_r^\infty|\partial_sg(s,\omega)|\,ds
 \leq
 \left(\int_r^\infty|s\partial_sg(s,\omega)|^p s^{d-1}\,ds\right)^{1/p}
 \left(\int_r^\infty s^{-(d+p-1)/(p-1)}\,ds\right)^{1/q}
 <\infty.
 $$
 Apply Minkowski's inequality first on $1<s<N$, and let
 $N\to\infty$ by Fatou's lemma. The change of variable $r\mapsto sr$ gives
 $$
 \begin{aligned}
 \|g\|_{L^p(E_d)}
 &\leq\int_1^\infty
 \|sr\,\partial_rg(sr,\omega)\|_{L^p(r>r_1)}\,\frac{ds}{s}
 \leq\int_1^\infty s^{-1-d/p}\,ds\,
 \|r\partial_rg\|_{L^p(E_d)}
 =\frac pd\|r\partial_rg\|_{L^p(E_d)}.
 \end{aligned}
 $$
 For the compact remainder $K_1$, choose the bounded connected
 domain $U$ and end annulus $D$ used in Lemma~\ref{lem:reverse-trace}.
 The $L^p$ Poincar\'e inequality gives
 $$
 \|g\|_{L^p(K_1)}
 \leq C_p\left(\|\nabla g\|_{L^p(U)}+\|g\|_{L^p(D)}\right).
 $$
 Consequently,
 \begin{equation}\label{eq:lp-trace}
 \|g\|_p\leq C_p\|\rho\nabla g\|_p
 \leq C_{p,c}\|A\|_p.
 \end{equation}

 Fix $j$, write $r=r_j$, and take
 $\varphi\in C_c^\infty(X_j)$. Put
 $$
 v=e^{-t\sqrt{\Delta_j}}\varphi,\qquad
 w=t\theta(t/r)\chi_j(x),\qquad
 V_\varphi(x)=\sup_{0<t<2r}
 \bigl(|v(x,t)|+t|\nabla_{x,t}v(x,t)|\bigr).
 $$
 By \eqref{eq:lp-model-maximal} and the model maximal theorem,
 \begin{equation}\label{eq:lp-test-maximal-norm}
 \|V_\varphi\|_{L^q(X_j)}\leq C_q\|\varphi\|_{L^q(X_j)}.
 \end{equation}
 Use \eqref{eq:reverse-local-green} with this $v$ and $w$.
 For its gradient term, put
 $$
 F=t\nabla u\,\mathbf1_{B_j\times(0,2r)},\qquad
 G=t\nabla v\,\mathbf1_{B_j\times(0,2r)}.
 $$
 The agreement of the model balls with the balls of $M$ gives
 $$
 \mathcal AF(x)\leq CA(x)\mathbf1_{3B_j}(x),\qquad
 \|G\|_{T^q}\leq C_q\|\varphi\|_{L^q(X_j)}.
 $$
 Tonelli's theorem, Cauchy--Schwarz on each cone, and H\"older's
 inequality give
 $$
 \begin{aligned}
 \int_0^\infty\int_{X_j}|F|\,|G|\,d\mu_j\,\frac{dt}{t}
 &=\int_{X_j}\int_0^\infty\int_{B(x,t)}
 \frac{|F(y,t)|\,|G(y,t)|}{\mu_j(B(y,t))}
 \,d\mu_j(y)\,\frac{dt}{t}\,d\mu_j(x)\\
 &\leq\int_{X_j}\mathcal AF(x)\mathcal AG(x)\,d\mu_j(x)\\
 &\leq C_q\|A\|_{L^p(3B_j)}\|\varphi\|_{L^q(X_j)}.
 \end{aligned}
 $$
 Since $|w|\leq Ct$, it follows that
 \begin{equation}\label{eq:lp-green-gradient}
 \left|\iint w\nabla u\cdot\nabla\overline v\right|
 \leq C_q\|A\|_{L^p(3B_j)}\|\varphi\|_{L^q(X_j)}.
 \end{equation}

 Equation~\eqref{eq:reverse-local-log-u} gives, on $B_j$,
 $$
 |u(x,t)|\leq|g(x)|+CA(x)(1+|\log(t/r)|),
 \qquad t|\nabla u(x,t)|\leq CA(x).
 $$
 Also,
 $$
 |w_{tt}|\leq Cr^{-1}\mathbf1_{\{r<t<2r\}},\qquad
 |\nabla_xw|\leq Ct/r,\qquad
 \int_0^{2r}(1+|\log(t/r)|)\,dt\leq Cr.
 $$
 Hence the two cutoff terms satisfy
 \begin{equation}\label{eq:lp-green-cutoffs}
 \begin{aligned}
 \iint |w_{tt}|\,|u|\,|v|
 +\iint|\nabla_xw|\,|\nabla_x(u\overline v)|
 &\leq C\int_{B_j}(|g|+A)V_\varphi\,d\mu\\
 &\leq C_q\bigl(\|g\|_{L^p(B_j)}+\|A\|_{L^p(B_j)}\bigr)
 \|\varphi\|_{L^q(X_j)}.
 \end{aligned}
 \end{equation}
 All these integrals use $d\mu\,dt$. The bottom cutoff
 $h_\tau(t)$ in the proof of \eqref{eq:reverse-local-green} gives
 the boundary term by $u(\cdot,t)\to f$ and
 $v(\cdot,t)\to\varphi$ in $L^2$ on the patch. Equations
 \eqref{eq:lp-green-gradient}--\eqref{eq:lp-green-cutoffs}
 give integrable bounds for the other terms as $\tau\downarrow0$.
 Thus
 $$
 \left|\int_M\chi_jf\overline\varphi\,d\mu\right|
 \leq C_q\left(\|A\|_{L^p(3B_j)}
 +\|g\|_{L^p(B_j)}\right)\|\varphi\|_{L^q(X_j)}.
 $$
 Taking the supremum over test functions with $L^q$ norm at most
 one gives
 \begin{equation}\label{eq:lp-local-boundary}
 \|\chi_jf\|_p\leq C_q
 \left(\|A\|_{L^p(3B_j)}+\|g\|_{L^p(B_j)}\right).
 \end{equation}
 Since $\sum_j\chi_j=1$ and the enlarged balls have bounded
 overlap, \eqref{eq:lp-trace} yields
 $$
 \begin{aligned}
 \|f\|_p^p
 &\leq C_p\sum_j\|\chi_jf\|_p^p\\
 &\leq C_{p,c}\sum_j
 \left(\int_{3B_j}A^p\,d\mu+\int_{B_j}|g|^p\,d\mu\right)\\
 &\leq C_{p,c}\bigl(\|A\|_p^p+\|g\|_p^p\bigr)
 \leq C_{p,c}\|A\|_p^p.
 \end{aligned}
 $$
 This proves the lower estimate on $L^p\cap L^2$.

 For $f\in L^p(M)$, choose $f_k\in C_c^\infty(M)$ with
 $f_k\to f$ in $L^p$. The upper estimate and the triangle
 inequality on each cone give
 $$
 \|\Sc f_k-\Sc f\|_p
 \leq\|\Sc(f_k-f)\|_p\leq C_{p,c}\|f_k-f\|_p\longrightarrow0.
 $$
 Passing to the limit in the lower estimate proves
 \eqref{eq:reverse-Lp}.
\end{proof}

	\subsection{Real interpolation}
	
	For a compatible pair $(X_0,X_1)$, write
	$$
	K(t,f;X_0,X_1)=\inf_{f=f_0+f_1}
	\bigl(\|f_0\|_{X_0}+t\|f_1\|_{X_1}\bigr),
	$$
	and use the real interpolation norm
	$$
	\|f\|_{(X_0,X_1)_{\theta,p}}
	=\left(\int_0^\infty
	[t^{-\theta}K(t,f;X_0,X_1)]^p\frac{\dd t}{t}\right)^{1/p}.
	$$
	
	Theorem \ref{thm:main} identifies $L^2(M)$ with
	$H^2_{S_c}(M)$. The identity to be proved is
	\begin{equation}\label{eq:interpolation}
		( H^1_{S_c}(M) ,L^2(M))_{\theta,p}=L^p(M),
		\qquad \theta=2-2/p,\quad 1<p<2,
	\end{equation}
	with equivalent norms. This is Theorem
	\ref{thm:hardy-characterization}(iii). The same assertion holds
	for the atomic space by \eqref{eq:prior-H1}.
	
	\section{A weighted Green estimate}
	
	Let $G$ be the minimal positive Green operator of $\Delta$.
	For $f\in C_c^\infty(M)$, the function $b=Gf$ is smooth and
	$\Delta b=f$. It need not belong to $L^2(M)$.
	
	\begin{lemma}\label{lem:weighted-G}
		For $1<p<n/2$,
		\begin{equation}\label{eq:weighted-G}
			\|\rho^{-2}Gf\|_p\leq C_p\|f\|_p
			\qquad(f\in C_c^\infty(M)).
		\end{equation}
	\end{lemma}
	
	\begin{proof}
		We first show that, for $0<\alpha<n-2$,
		\begin{equation}\label{eq:G-power}
			G(\rho^{-\alpha-2})(x)\leq C_\alpha\rho(x)^{-\alpha}.
		\end{equation}
		Choose a smooth positive function $v_\alpha$ equal to
		$r^{-\alpha}$ far out on both ends. On the end of dimension
		$d=m$ or $d=n$,
		$$
		\Delta(r^{-\alpha})=\alpha(d-2-\alpha)r^{-\alpha-2}>0.
		$$
		There are $k_\alpha\in C_c^\infty(M)$, $k_\alpha\geq0$, and
		$c_\alpha>0$ such that
		$$
		\Delta v_\alpha+k_\alpha\geq c_\alpha\rho^{-\alpha-2}.
		$$
		The heat kernel estimates in \cite{GS}, integrated in time,
		give $Gk_\alpha(x)\leq C r^{2-d}$ on the end of dimension $d$.
		Consequently
		$$
		w_\alpha=v_\alpha+Gk_\alpha\simeq\rho^{-\alpha},\qquad
		\Delta w_\alpha\geq c_\alpha\rho^{-\alpha-2},\qquad w_\alpha>0.
		$$
		On each bounded smooth exhaustion domain, comparison with
		$c_\alpha^{-1}w_\alpha$ bounds the Dirichlet Green potential of
		$\rho^{-\alpha-2}$. Increasing the domains gives
		\eqref{eq:G-power}.
		
		Let $p'=p/(p-1)$ and choose
		\begin{equation}\label{eq:beta}
			\frac2{p'}<\beta<\min\left\{\frac n{p'},\frac{n-2}{p}\right\}.
		\end{equation}
		This interval is nonempty exactly when $p<n/2$.
		Put $h(x)=\rho(x)^{-\beta}$. Also put
		$J(x,y)=\rho(x)^{-2}G(x,y)$. By symmetry of $G$ and
		\eqref{eq:G-power},
		\begin{align}
			\int_M J(x,y)h(y)^{p'}\dd\mu(y)
			&\leq C h(x)^{p'},\label{eq:schur1}\\
			\int_M J(x,y)h(x)^p\dd\mu(x)
			&\leq C h(y)^p.\label{eq:schur2}
		\end{align}
		For \eqref{eq:schur1} use $\alpha=\beta p'-2$; for
		\eqref{eq:schur2} use $\alpha=\beta p$.
		Both lie in $(0,n-2)$ by \eqref{eq:beta}.
		H\"older's inequality now gives
		$$
		\left|\int_M J(x,y)f(y)\dd\mu(y)\right|^p
		\leq C h(x)^p\int_M J(x,y)|f(y)|^p h(y)^{-p}\dd\mu(y).
		$$
		Integrate in $x$ and apply \eqref{eq:schur2}. This proves
		\eqref{eq:weighted-G}.
	\end{proof}

	\section{Local estimates and a decomposition}
	
	Fix a sufficiently small $a>0$. Define the uncentered maximal
	operator
	\begin{equation}\label{eq:admissible-max}
		\mathcal M_a F(x)
		=\sup_{\substack{x\in B(z,r)\\0<r\leq a\rho(z)}}
		\frac1{\mu(B(z,r))}\int_{B(z,r)}|F|\dd\mu.
	\end{equation}
	The fivefold covering lemma and doubling for these balls give
	\begin{equation}\label{eq:maximal-p}
		\mu\{\mathcal M_aF>\lambda\}\leq C\lambda^{-1}\|F\|_1,
		\qquad \|\mathcal M_a F\|_s\leq C_s\|F\|_s\quad(s>1).
	\end{equation}
	The enlargement of each selected ball is still within the range
	of local doubling, since $a$ is fixed sufficiently small.
	
	\begin{lemma}\label{lem:q-atom}
		If $1<q<2$, $v\in C_c^\infty(B)$, and $r_B\leq\delta\rho(x_B)$,
		then
		\begin{equation}\label{eq:q-atom}
			\|\Delta v\|_{ H^1_{S_c}(M) }
			\leq C_q\mu(B)^{1-1/q}
			\bigl(\|\Delta v\|_q+r_B^{-2}\|v\|_q\bigr).
		\end{equation}
	\end{lemma}
	
	\begin{proof}
		On $2B$, H\"older's inequality and \eqref{eq:prior-N} give
		$$
		\int_{2B}\mathcal N_c^\#(\Delta v)\dd\mu
		\leq C_q\mu(B)^{1-1/q}\|\Delta v\|_q.
		$$
		On $M\setminus2B$, use $P_t\Delta v=\Delta P_tv$ and
		\eqref{eq:prior-kernel}:
		$$
		\int_{M\setminus2B}\mathcal N_c^\#(\Delta v)\dd\mu
		\leq Cr_B^{-2}\|v\|_1
		\leq Cr_B^{-2}\mu(B)^{1-1/q}\|v\|_q.
		$$
		The second inequality in \eqref{eq:prior-N} proves the result.
	\end{proof}
	
	We record the local elliptic estimates used below. If
	$B=B(x,r)$ and $32r\leq a_0\rho(x)$, then $8B$ lies in a
	complete model $X$ whose Laplacian agrees with $\Delta$ there.
	On the ends take $X=\mathbb R^m$ or
	$X=\mathbb R^n\times\mathbb S^{m-n}$.
	Near the compact part take one of finitely many smooth metrics
	on $\mathbb R^m$ equal to the Euclidean metric outside a compact
	set. Let $G_X$ be its positive Green kernel. For $1<q<2$ put
	$$
	v=G_X(\mathbf1_{8B}f).
	$$
	Then $\Delta v=f$ on $8B$, and
	\begin{equation}\label{eq:local-potential}
		\|v\|_{L^q(4B)}+r\|\nabla v\|_{L^q(4B)}
		\leq C_qr^2\|f\|_{L^q(8B)}.
	\end{equation}
	Here is a direct justification of the uniform constant.
	On the product models, integration of their heat kernels gives,
	with $d(x,y)=s$,
	$$
	|G_X(x,y)|\leq C
	\begin{cases}s^{2-m},&s\leq1,\\s^{2-d},&s>1,\end{cases}
	\quad
	|\nabla_xG_X(x,y)|\leq C
	\begin{cases}s^{1-m},&s\leq1,\\s^{1-d},&s>1,\end{cases}
	$$
	where $d=m$ or $d=n$. The ball volumes have orders $s^m$ for
	$s\leq1$ and $s^d$ for $s>1$. Therefore, in either variable,
	$$
	\int_{d(x,y)<Cr}|G_X(x,y)|\dd\mu\leq Cr^2,
	\qquad
	\int_{d(x,y)<Cr}|\nabla_xG_X(x,y)|\dd\mu\leq Cr.
	$$
	For the finitely many compact coordinate models, $r$ is bounded
	and the same estimates follow from the local Green singularity.
	For completeness, let $K(x,y)$ be either of these kernels,
	restricted to $4B\times8B$, and put
	$$
	A_K=\sup_{x\in4B}\int_{8B}|K(x,y)|\dd\mu(y),\qquad
	B_K=\sup_{y\in8B}\int_{4B}|K(x,y)|\dd\mu(x).
	$$
	H\"older's inequality followed by Tonelli gives Schur's estimate:
	$$
	\int_{4B}\left|\int_{8B}K(x,y)f(y)\dd\mu(y)\right|^q\dd\mu(x)
	\leq A_K^{q-1}B_K\|f\|_{L^q(8B)}^q.
	$$
	For $G_X$, both constants are at most $Cr^2$; for
	$\nabla_xG_X$, both are at most $Cr$. This proves
	\eqref{eq:local-potential}.
	
	For $\Delta h=0$ on $4B$, the local harmonic estimates give
	\begin{equation}\label{eq:local-harmonic}
		\sup_{2B}(|h|+r|\nabla h|)
		\leq C_q\left(\frac1{\mu(4B)}\int_{4B}|h|^q\dd\mu\right)^{1/q}.
	\end{equation}
	Subtracting the potential in \eqref{eq:local-potential} from
	a solution of $\Delta b=f$ gives, in particular,
	\begin{equation}\label{eq:local-grad-b}
		r\|\nabla b\|_{L^q(2B)}
		\leq C_q\bigl(\|b\|_{L^q(4B)}+r^2\|f\|_{L^q(8B)}\bigr).
	\end{equation}
	
	\begin{proposition}\label{prop:decomposition}
		Fix $1<q<p<\min\{2,n/2\}$ and $f\in C_c^\infty(M)$. Set
		\begin{equation}\label{eq:F-Q}
			b=Gf,\qquad F=|f|+\rho^{-2}|b|,\qquad
			Q=(\mathcal M_a(F^q))^{1/q}.
		\end{equation}
		For every $\lambda>0$ there is a decomposition
		$f=B_\lambda+g_\lambda$ such that
		\begin{align}
			\|B_\lambda\|_{ H^1_{S_c}(M) }
			&\leq C_q\int_{\{Q>\lambda\}}Q\dd\mu,\label{eq:bad-est}\\
			\|g_\lambda\|_2^2
			&\leq C_q\left(\int_{\{Q\leq\lambda\}}Q^2\dd\mu
			+\lambda^2\mu\{Q>\lambda\}\right).\label{eq:good-est}
		\end{align}
		Moreover,
		\begin{equation}\label{eq:Q-p}
			\|Q\|_p\leq C_{p,q}\|f\|_p.
		\end{equation}
	\end{proposition}
	
	\begin{proof}
		Equation \eqref{eq:Q-p} follows from
		\eqref{eq:maximal-p} and Lemma \ref{lem:weighted-G}.
		Let $\Omega=\{Q>\lambda\}$. This set is open and has finite
		measure. If it is empty, take $B_\lambda=0$.
		
		Fix $0<\eta\leq10^{-3}\min\{1,a,a_0,\delta\}$.
		On $\Omega$ set
		$$
		r(x)=\frac1{100}\min\{d(x,M\setminus\Omega),\eta\rho(x)\}.
		$$
		Choose a covering $B_i=B(x_i,r_i)$, $r_i=r(x_i)$, with
		the following properties:
		\begin{equation}\label{eq:Whitney}
			\Omega=\bigcup_i B_i,\qquad 16B_i\subset\Omega,
			\qquad \sum_i\mathbf1_{16B_i}\leq C.
		\end{equation}
		If $2B_i\cap2B_j\ne\varnothing$, then $r_i\simeq r_j$.
		These assertions follow by taking a maximal disjoint family of
		$B(x,r(x)/5)$: the function $r$ is $1/100$-Lipschitz, so
		intersecting enlargements have comparable radii; local doubling
		bounds their overlap. Choose a smooth partition on $\Omega$ with
		\begin{equation}\label{eq:partition}
			\sum_i\varphi_i=1,\qquad
			\supp\varphi_i\subset2B_i,\qquad
			|\varphi_i|+r_i|\nabla\varphi_i|+r_i^2|\Delta\varphi_i|\leq C.
		\end{equation}
		Such cutoffs exist on the two product models and on the finitely
		many compact coordinate neighborhoods. Dividing their sum gives
		\eqref{eq:partition}; \eqref{eq:Whitney} bounds its derivatives.
		
		Call an index $i$ a boundary index if
		$d(x_i,M\setminus\Omega)\leq\eta\rho(x_i)$.
		For such an index let
		$$
		v_i=G_{X_i}(\mathbf1_{8B_i}f),\qquad h_i=b-v_i
		\quad\text{on }8B_i.
		$$
		Then $\Delta h_i=0$ on $8B_i$. For every other index put
		$v_i=b$ and $h_i=0$. Define
		\begin{equation}\label{eq:pieces}
			b_i=\varphi_i v_i,\qquad a_i=\Delta b_i.
		\end{equation}
		Each $b_i$ is smooth and compactly supported in $2B_i$.
		For a boundary index, \eqref{eq:local-potential} and
		\eqref{eq:partition} give
		$$
		\|a_i\|_q+r_i^{-2}\|b_i\|_q
		\leq C_q\|f\|_{L^q(8B_i)}.
		$$
		For any other index, $r_i=\eta\rho(x_i)/100$ and
		$\rho\simeq r_i/\eta$ on $8B_i$.
		Equations \eqref{eq:local-grad-b} and \eqref{eq:partition} give
		the bound by $C_{q,\eta}\|F\|_{L^q(8B_i)}$ instead.
		Consequently, Lemma \ref{lem:q-atom} implies
		\begin{equation}\label{eq:piece-H1}
			\|a_i\|_{ H^1_{S_c}(M) }
			\leq C_q\mu(B_i)^{1-1/q}\|F\|_{L^q(8B_i)}
			\leq C_q\int_{B_i}Q\dd\mu.
		\end{equation}
		The last inequality follows because $8B_i$ is one of the balls
		allowed in \eqref{eq:admissible-max} for each $x\in B_i$.
		By \eqref{eq:Whitney},
		$$
		\sum_i\|a_i\|_{ H^1_{S_c}(M) }
		\leq C_q\int_\Omega Q\dd\mu
		\leq C_q\lambda^{1-p}\|Q\|_p^p<\infty.
		$$
		Thus $B_\lambda:=\sum_i a_i$ converges in $ H^1_{S_c}(M) $ and
		in $L^1(M)$. This proves \eqref{eq:bad-est}. Define
		$g_\lambda=f-B_\lambda$ as a function in $L^1+L^2$.
		Since every $a_i$ is supported in $\Omega$,
		\begin{equation}\label{eq:good-outside}
			g_\lambda=f\quad\text{almost everywhere on }M\setminus\Omega.
		\end{equation}
		
		We prove $|g_\lambda|\leq C_q\lambda$ on $\Omega$.
		The sums below are locally finite in this open set, so
		\begin{equation}\label{eq:good-inside}
			g_\lambda=\Delta\left(\sum_i\varphi_i h_i\right).
		\end{equation}
		At a point $x\in\Omega$, only indices with $x\in2B_i$ can
		contribute. If none is a boundary index, every $h_i=0$
		there and \eqref{eq:good-inside} is zero.
		Otherwise choose one boundary index $k$ with $x\in2B_k$.
		There is $z_k\in M\setminus\Omega$ with
		$d(z_k,x_k)\leq101r_k$. The radii of all contributing balls
		are comparable to $r_k$.
		More precisely, the Lipschitz bound for $r$ gives
		$$
		|r_i-r_k|\leq\frac{d(x_i,x_k)}{100}
		<\frac{r_i+r_k}{50},\qquad
		\frac12r_k\leq r_i\leq2r_k.
		$$
		The ball $D_k=B(x_k,128r_k)$ contains $z_k$ and all these $8B_i$:
		$$
		d(x_k,y)<2(r_i+r_k)+8r_i\leq22r_k
		\quad(y\in8B_i),\qquad 128r_k\leq a\rho(x_k).
		$$
		It is therefore allowed in \eqref{eq:admissible-max}.
		Local doubling gives $\mu(D_k)\leq C\mu(B_i)$.
		Since $Q(z_k)\leq\lambda$,
		\begin{equation}\label{eq:average-goodpoint}
			\|F\|_{L^q(8B_i)}\leq C\lambda\mu(B_i)^{1/q}.
		\end{equation}
		
		Put $D=B(x,r_k/4)$. Then $4D\subset4B_i\cap4B_k$ and
		$\mu(4D)\simeq\mu(B_k)$ for every contributing index $i$.
		If $i$ is a boundary index, $h_i-h_k=v_k-v_i$ is harmonic on
		$4D$. Equations \eqref{eq:local-potential} and
		\eqref{eq:average-goodpoint} give
		$$
		\|h_i-h_k\|_{L^q(4D)}
		\leq\|v_i\|_{L^q(4B_i)}+\|v_k\|_{L^q(4B_k)}
		\leq C_q\lambda r_k^2\mu(B_k)^{1/q}.
		$$
		Apply \eqref{eq:local-harmonic} on $D$ to obtain
		\begin{equation}\label{eq:h-difference}
			|h_i(x)-h_k(x)|+r_k|\nabla h_i(x)-\nabla h_k(x)|
			\leq C_q\lambda r_k^2.
		\end{equation}
		For any other index, $h_i=0$ and $r_i\simeq\eta\rho$.
		On $4D$, equation \eqref{eq:average-goodpoint} gives
		$$
		\|b\|_{L^q(4D)}\leq C_\eta r_k^2\|F\|_{L^q(8B_i)}
		\leq C_\eta\lambda r_k^2\mu(B_k)^{1/q}.
		$$
		Since $h_k=b-v_k$ is harmonic on $4D$, we have
		$$
		\|h_k\|_{L^q(4D)}
		\leq\|b\|_{L^q(4D)}+\|v_k\|_{L^q(4B_k)}
		\leq C_{q,\eta}\lambda r_k^2\mu(B_k)^{1/q}.
		$$
		Equation \eqref{eq:local-harmonic} proves
		\eqref{eq:h-difference} also in this case.
		
		Use $\sum_i\nabla\varphi_i=0$, $\sum_i\Delta\varphi_i=0$, and
		$\Delta h_i=0$ near $x$. Equations \eqref{eq:good-inside} and
		\eqref{eq:h-difference} yield
		\begin{align*}
			|g_\lambda(x)|
			&=\left|\sum_i\bigl((\Delta\varphi_i)(h_i-h_k)
			-2\nabla\varphi_i\cdot\nabla(h_i-h_k)\bigr)(x)\right|
			\leq C_q\lambda.
		\end{align*}
		Together with \eqref{eq:good-outside} and $|f|\leq Q$ almost
		everywhere, this proves \eqref{eq:good-est}.
		The atomic series converges in $L^1$ to $B_\lambda$, and the
		last estimate gives $g_\lambda\in L^2$.
	\end{proof}
	
	\section{The interpolation identity and operator bounds}
	\label{sec:L2-interpolation}
	
	\begin{proof}[Proof of Theorem \ref{thm:hardy-characterization}(iii)]
		The inclusion $ H^1_{S_c}(M) \hookrightarrow L^1$ and real
		interpolation of Lebesgue spaces give, for $1<p<2$,
		\begin{equation}\label{eq:easy-inclusion}
			( H^1_{S_c}(M) ,L^2)_{\theta,p}\hookrightarrow L^p,
			\qquad \theta=2-2/p.
		\end{equation}
		We prove the reverse inclusion first for
		$1<p<\min\{2,n/2\}$.
		Choose $1<q<p$ and let $Q$ be as in \eqref{eq:F-Q}.
		Proposition \ref{prop:decomposition} gives, for every
		$t,\lambda>0$,
		\begin{equation}\label{eq:K-level}
			\begin{split}
				K(t,f; H^1_{S_c}(M) ,L^2)
				\leq C_q\left\{
				\int_{\{Q>\lambda\}}Q\dd\mu
				+t\left(\int_{\{Q\leq\lambda\}}Q^2\dd\mu
				+\lambda^2\mu\{Q>\lambda\}\right)^{1/2}\right\}.
			\end{split}
		\end{equation}
		Write
		$$
		d_Q(\alpha)=\mu\{Q>\alpha\},\qquad
		Q^*(s)=\inf\{\alpha\geq0:d_Q(\alpha)\leq s\}.
		$$
		This is the right-continuous decreasing rearrangement.
		Fix $t>0$, put $a=t^2$, and first suppose that
		$\lambda:=Q^*(a)>0$. Then
		$$
		\ell:=d_Q(\lambda)\leq a,\qquad
		Q^*(s)=\lambda\quad\text{for almost every }\ell<s<a.
		$$
		Indeed, $s<a$ gives $Q^*(s)\geq Q^*(a)=\lambda$, whereas
		$s>\ell$ gives $Q^*(s)\leq\lambda$.
		Integration of the distribution functions gives
		$$
		\int_{\{Q>\lambda\}}Q\dd\mu
		=\int_0^\ell Q^*(s)\dd s\leq\int_0^aQ^*(s)\dd s,
		$$
		and
		$$
		\int_{\{Q\leq\lambda\}}Q^2\dd\mu+\lambda^2\ell
		=\int_\ell^\infty Q^*(s)^2\dd s+\lambda^2\ell
		=\int_a^\infty Q^*(s)^2\dd s+\lambda^2a.
		$$
		Since $t\sqrt a\,\lambda=a\lambda\leq\int_0^aQ^*(s)\dd s$,
		equation \eqref{eq:K-level} yields
		\begin{equation}\label{eq:K-rearranged}
			K(t,f; H^1_{S_c}(M) ,L^2)
			\leq C_q\left\{\int_0^{t^2}Q^*(s)\dd s
			+t\left(\int_{t^2}^\infty Q^*(s)^2\dd s\right)^{1/2}
			\right\}.
		\end{equation}
		If $Q^*(a)=0$, then $d_Q(\alpha)\leq a$ for every $\alpha>0$,
		so monotone convergence gives $\mu\{Q>0\}\leq a$. Therefore
		$$
		\|Q\|_1\leq a^{1-1/p}\|Q\|_p<\infty,
		\qquad \int_0^aQ^*(s)\dd s=\|Q\|_1,
		\qquad Q^*(s)=0\quad(s>a).
		$$
		Apply \eqref{eq:K-level} with $\lambda>0$ and let
		$\lambda\downarrow0$. Its terms satisfy
		$$
		\int_{\{Q>\lambda\}}Q\dd\mu\longrightarrow\|Q\|_1,
		$$
		and
		$$
		\int_{\{Q\leq\lambda\}}Q^2\dd\mu
		+\lambda^2\mu\{Q>\lambda\}
		\leq2\lambda^{2-p}\|Q\|_p^p\longrightarrow0.
		$$
		This proves \eqref{eq:K-rearranged} also when $Q^*(a)=0$.
		
		We now integrate \eqref{eq:K-rearranged} directly.
		Put $s=t^2$ and use $\theta=2-2/p$. The first term gives
		$$
		\int_0^\infty
		\left[t^{-\theta}\int_0^{t^2}Q^*(v)\dd v\right]^p\frac{\dd t}{t}
		=\frac12\int_0^\infty
		\left[\frac1s\int_0^sQ^*(v)\dd v\right]^p\dd s.
		$$
		Minkowski's inequality gives its bound, since
		$$
		\left\|\frac1s\int_0^sQ^*(v)\dd v\right\|_{L^p(ds)}
		=\left\|\int_0^1Q^*(su)\dd u\right\|_{L^p(ds)}
		\leq\int_0^1u^{-1/p}\dd u\,\|Q\|_p
		=p'\|Q\|_p.
		$$
		Here $p'=p/(p-1)$ and $p>1$.
		
		The second term gives
		$$
		\int_0^\infty
		\left[t^{1-\theta}
		\left(\int_{t^2}^\infty Q^*(v)^2\dd v\right)^{1/2}\right]^p
		\frac{\dd t}{t}
		=\frac12\int_0^\infty
		\left[\frac1s\int_s^\infty Q^*(v)^2\dd v\right]^{p/2}\dd s.
		$$
		Since $Q^*$ is decreasing,
		$$
		\int_s^\infty Q^*(v)^2\dd v
		=\sum_{k\geq0}\int_{2^ks}^{2^{k+1}s}Q^*(v)^2\dd v
		\leq\sum_{k\geq0}2^ks\,Q^*(2^ks)^2.
		$$
		Use $(\sum_k a_k)^{p/2}\leq\sum_k a_k^{p/2}$ for $a_k\geq0$,
		followed by Tonelli and the change of variable $v=2^ks$:
		\begin{align*}
			\int_0^\infty
			\left[\frac1s\int_s^\infty Q^*(v)^2\dd v\right]^{p/2}\dd s
			&\leq\sum_{k\geq0}2^{kp/2}
			\int_0^\infty Q^*(2^ks)^p\dd s
			=\sum_{k\geq0}2^{k(p/2-1)}\|Q\|_p^p
			=\frac{\|Q\|_p^p}{1-2^{p/2-1}}.
		\end{align*}
		The denominator is positive because $p<2$.
		Combining the two bounds with \eqref{eq:K-rearranged} and
		\eqref{eq:Q-p} gives
		\begin{equation}\label{eq:first-range}
			\|f\|_{( H^1_{S_c}(M) ,L^2)_{\theta,p}}
			\leq C_{p,q}\|Q\|_p\leq C_{p,q}\|f\|_p
			\qquad(f\in C_c^\infty(M)).
		\end{equation}
		Density of $C_c^\infty$ in $L^p$, completeness of the
		interpolation space, and \eqref{eq:easy-inclusion} extend
		\eqref{eq:first-range} to all $f\in L^p$.
		
		Now fix any $1<p<2$. Choose
		$$
		1<p_0<\min\{p,n/2\},\qquad
		\theta_0=2-2/p_0,\qquad
		\eta=\frac{\theta-\theta_0}{1-\theta_0}\in(0,1).
		$$
		The case just proved gives
		$( H^1_{S_c}(M) ,L^2)_{\theta_0,p_0}=L^{p_0}$.
		The real reiteration theorem
		\cite[Theorem 3.5.3]{BL} gives
		\begin{align*}
			( H^1_{S_c}(M) ,L^2)_{\theta,p}
			&=\bigl(( H^1_{S_c}(M) ,L^2)_{\theta_0,p_0},L^2\bigr)_{\eta,p}
			=(L^{p_0},L^2)_{\eta,p}=L^p,
			\qquad \frac1p=\frac{1-\eta}{p_0}+\frac\eta2.
		\end{align*}
		This proves \eqref{eq:interpolation}. Equation
		\eqref{eq:prior-H1} proves the assertion for the atomic space.
		For $n=3$, the first part uses $1<p<3/2$; the last argument
		allows every $1<p<2$ by taking $1<p_0<\min\{p,3/2\}$.
	\end{proof}
	
	\begin{corollary}\label{cor:operator}
		Let $T$ have compatible linear extensions from $ H^1_{S_c}(M) $ to
		$L^1(M)$ and from $L^2(M)$ to $L^2(M)$. Suppose
		$$
		\|Tf\|_1\leq A_1\|f\|_{ H^1_{S_c}(M) },\qquad
		\|Tf\|_2\leq A_2\|f\|_2.
		$$
		Then, for $1<p<2$,
		\begin{equation}\label{eq:T-p}
			\|Tf\|_p\leq C_p A_1^{\,2/p-1}A_2^{\,2-2/p}\|f\|_p.
		\end{equation}
	\end{corollary}
	
	\begin{proof}
		Compatibility means that the two extensions agree on
		$ H^1_{S_c}(M) \cap L^2$. They therefore define one linear map on
		$ H^1_{S_c}(M) +L^2$. Assume first that $A_1,A_2>0$.
		For $f=f_0+f_1$, $f_0\in H^1_{S_c}(M) $, $f_1\in L^2$, we have
		$$
		K(t,Tf;L^1,L^2)
		\leq A_1\|f_0\|_{ H^1_{S_c}(M) }+tA_2\|f_1\|_2.
		$$
		Take the infimum over all such decompositions:
		$$
		K(t,Tf;L^1,L^2)
		\leq A_1 K(tA_2/A_1,f; H^1_{S_c}(M) ,L^2).
		$$
		Integration with $s=tA_2/A_1$ gives
		$$
		\|Tf\|_{(L^1,L^2)_{\theta,p}}
		\leq A_1^{1-\theta}A_2^\theta
		\|f\|_{( H^1_{S_c}(M) ,L^2)_{\theta,p}}.
		$$
		Apply Theorem \ref{thm:hardy-characterization}(iii),
		$(L^1,L^2)_{\theta,p}=L^p$, and $\theta=2-2/p$.
		If either bound is zero, replace $A_i$ by $A_i+\epsilon$ and
		let $\epsilon\downarrow0$.
	\end{proof}

	\section{Spectral multipliers}
	\label{sec:spectral-multipliers}
	
	We prove Theorem \ref{thm:spectral-multiplier}. We keep the
	exact ends, $K$, $\rho$, and $\Delta$ fixed above. The complete
	models are $M_m=\R^m$ and $M_n=\R^n\times S^{m-n}$. Write $d_i$
	for their distances and $\Delta_i$ for their nonnegative Laplacians.
	For a compactly supported $C^q$ function, put
	$$
	\norm G_{C^q}=\max_{0\leq j\leq q}\sup|G^{(j)}|.
	$$
	The term $j=0$ is included.
	
	For the fixed $\delta$ in Theorem \ref{thm:hardy-characterization},
	\begin{equation}\label{eq:hw-local-doubling}
		\mu(2B)\leq C\mu(B),\qquad
		B=B(x_B,r),\quad0<r\leq\delta\rho(x_B).
	\end{equation}
	When $\rho(x_B)$ is large, $2B$ lies in one exact end and the
	product ball volumes give this bound. When $\rho(x_B)$ is
	bounded, \eqref{eq:volume-upper} and
	\eqref{eq:bounded-volume-lower} give the same bound.
	
	\begin{prop}\label{prop:mihlin-atom}
		Assume the multiplier hypotheses of Theorem
		\ref{thm:spectral-multiplier}, and fix $\delta$ as in
		Theorem \ref{thm:hardy-characterization}.
		For each integer $J\geq1$, let $a=\Delta^Jb$ be associated with
		$B=B(x_B,r)$, where
		\begin{equation}\label{eq:mihlin-2}
			\begin{gathered}
				b\in\operatorname{Dom}(\Delta^J),\qquad r\leq\delta\rho(x_B),\\
				\supp \Delta^j b\subset B,\qquad
				\|\Delta^j b\|_2\leq r^{2(J-j)}\mu(B)^{-1/2},\qquad0\leq j\leq J.
			\end{gathered}
		\end{equation}
		Then
		\begin{equation}\label{eq:mihlin-5}
			\|\mathfrak m(\Delta)a\|_1\leq C_{J,N,\delta}A_N(\mathfrak m).
		\end{equation}
		The constant in \eqref{eq:mihlin-5} is independent of the atom,
		its radius, and its end. For $J=1$, these are the atoms in
		\eqref{eq:atoms}.
	\end{prop}
	
	Subsection \ref{subsec:mihlin-atom} proves this proposition.
	For $J=1$, these are exactly the atoms in
	Theorem \ref{thm:hardy-characterization}(i).
	
	The kernel estimate used throughout is the following. Here
	$K_T(x,z)$ denotes the kernel of $T$ with integration in $z$.
	\begin{proposition}\label{prop:mihlin-bands}
		Let $N$ be an integer with $N>(m+3)/2$, and put $\varepsilon=1/8$.
		For $F\in C_c^N((1/4,4))$ and every $s>0$,
		\begin{align}
			\sup_z\int_M|K_{F(s^2\Delta)}(x,z)|\dd\mu(x)
			&\leq C_N\norm F_{C^N},\label{eq:mihlin-6}\\
			\sup_{z\in B}\int_{M\setminus2B}|K_{F(s^2\Delta)}(x,z)|\dd\mu(x)
			&\leq C_{N,\delta}\norm F_{C^N}(s/r)^\varepsilon,
			\quad s<r\leq\delta\rho(x_B).
			\label{eq:mihlin-7}
		\end{align}
	\end{proposition}
	
	\subsection{The heat kernel estimates and finite propagation}
	Let $h_t(x,z)$ be the real heat kernel. For $s\geq1$, the nine
	ordered cases correspond to $x,z\in K,E_m,E_n$. By symmetry it
	is enough to list six cases. The bounds of Grigor'yan and
	Saloff-Coste \cite[Theorem 6.6 and Corollary 6.13]{GS} give
	\begin{align*}
		h_{s^2}(x,z)&\leq Cs^{-n}e^{-cd(x,z)^2/s^2},
		&&x,z\in K,\\
		h_{s^2}(x,z)&\leq C\bigl(s^{-m}+s^{-n}\rho(x)^{2-m}\bigr)
		e^{-cd(x,z)^2/s^2},
		&&x\in E_m, z\in K,\\
		h_{s^2}(x,z)&\leq Cs^{-n}\bigl(1+\rho(x)^{2-n}\bigr)
		e^{-cd(x,z)^2/s^2},
		&&x\in E_n, z\in K,\\
		h_{s^2}(x,z)&\leq C\bigl(s^{-n}\rho(x)^{2-m}
		+s^{-m}\rho(z)^{2-n}\bigr)
		e^{-cd(x,z)^2/s^2},
		&&x\in E_m, z\in E_n,\\
		h_{s^2}(x,z)&\leq Cs^{-m}e^{-cd(x,z)^2/s^2}
		+Cs^{-n}\rho(x)^{2-m}\rho(z)^{2-m}
		e^{-c(\rho(x)^2+\rho(z)^2)/s^2},
		&&x,z\in E_m,\\
		h_{s^2}(x,z)&\leq Cs^{-n}e^{-cd(x,z)^2/s^2}
		+Cs^{-n}\rho(x)^{2-n}\rho(z)^{2-n}
		e^{-c(\rho(x)^2+\rho(z)^2)/s^2},
		&&x,z\in E_n.
	\end{align*}
	For the additional terms in the last two lines,
	$$
	\rho(x)^{2-m}\rho(z)^{2-m}\leq1,\qquad
	\rho(x)^{2-n}\rho(z)^{2-n}\leq1,\qquad
	d(x,z)^2\leq C\bigl(\rho(x)^2+\rho(z)^2\bigr).
	$$
	Since $s^{-m}\leq s^{-n}$ for $s\geq1$, all six bounds imply
	\begin{equation}\label{eq:hw-coarse-heat}
		0\leq h_{s^2}(x,z)\leq Cs^{-n}e^{-cd(x,z)^2/s^2},
		\qquad s\geq1.
	\end{equation}
	Lemma 5.9 of \cite{GS} gives the bounded-time Gaussian bound.
	Together with \eqref{eq:bounded-volume-lower}, it gives
	\begin{equation}\label{eq:hw-bounded-time-heat}
		0\leq h_{s^2}(x,z)\leq C_Ts^{-m}e^{-c_Td(x,z)^2/s^2},
		\qquad0<s\leq T<\infty.
	\end{equation}
	For each fixed $s_*>0$, \eqref{eq:hw-coarse-heat} therefore holds
	for $s\geq s_*$, with a constant depending on $s_*$.
	The heat semigroup is sub-Markovian, so
	\begin{equation}\label{eq:hw-mass}
		\int_M h_t(x,z)\dd\mu(x)\leq1.
	\end{equation}
	\begin{align}
		\norm{h_{s^2}(\cdot,z)}_2^2
		&\leq\sup_x h_{s^2}(x,z)\int_Mh_{s^2}(x,z)\dd\mu(x)
		\leq Cs^{-n},\label{eq:hw-heat-column}\\
		\norm{\mathbf1_{\{d(\cdot,z)>R\}}h_{s^2}(\cdot,z)}_2^2
		&\leq\sup_{d(x,z)>R}h_{s^2}(x,z)
		\int_Mh_{s^2}(x,z)\dd\mu(x)\nonumber\\
		&\leq Cs^{-n}e^{-cR^2/s^2},\qquad R\geq s.
		\label{eq:hw-heat-tail}
	\end{align}
	
	The wave equation has finite propagation:
	\begin{equation}\label{eq:hw-finite-propagation}
		\supp K_{\cos(t\sqrt \Delta)}\subset\{(x,z):d(x,z)\leq|t|\}.
	\end{equation}
	Also $\norm{\cos(t\sqrt \Delta)}_{2\to2}\leq1$.
	
	For completeness, let $u(t)=\cos(t\sqrt \Delta)f$, where
	$f\in C_c^\infty(M)$. For $t\geq0$, set
	$$
	r_f(x)=d(x,\supp f),\qquad
	\Omega_t=\{x:r_f(x)>t\},\qquad
	E(t)=\frac12\int_{\Omega_t}
	(|\partial_tu|^2+|\nabla u|^2)\dd\mu.
	$$
	The wave equation gives
	$$
	\partial_t\frac{|\partial_tu|^2+|\nabla u|^2}{2}
	=\operatorname{div}\Re(\overline{\partial_tu}\nabla u).
	$$
	Since $|\nabla r_f|=1$ almost everywhere outside $\supp f$,
	the coarea formula and integration by parts give, for almost
	every $t>0$,
	\begin{align*}
		E'(t)&=-\int_{\{r_f=t\}}
		\left\{\frac{|\partial_tu|^2+|\nabla u|^2}{2}
		+\Re(\overline{\partial_tu}\nabla u\cdot\nabla r_f)\right\}
		\dd\sigma,\\
		E'(t)&\leq0,\qquad E(0)=0,\qquad E(t)=0.
	\end{align*}
	One may first insert a compact spatial cutoff and approximate
	the Lipschitz function $r_f$; the finite energy of $u$ permits
	passage to the limit. If $r_f(x)>t$, then
	$$
	u(x,t)=f(x)+\int_0^t\partial_\tau u(x,\tau)\dd\tau=0.
	$$
	Evenness in $t$ and density give \eqref{eq:hw-finite-propagation}.
	
	The Fourier support form of \eqref{eq:hw-finite-propagation} is
	Lemma 2.1 of Cowling and Sikora \cite{CS}.
	
	The volume bound \eqref{eq:volume-upper} will be used below.
	
	\subsection{Spectral kernels on a fixed compact set}
	\begin{lemma}\label{lem:hw-point}
		Fix $b_0<\infty$ and $s_*>0$. If $G$ is a bounded Borel
		function supported in $[0,b_0]$,
		then, for $s\geq s_*$,
		\begin{equation}\label{eq:hw-point-uniform}
			|K_{G(s^2\Delta)}(x,z)|\leq Cs^{-n}\norm G_\infty.
		\end{equation}
		If $G\in C^q([0,\infty))$, with compact support and $q\geq1$, then
		\begin{equation}\label{eq:hw-point-tail}
			|K_{G(s^2\Delta)}(x,z)|\leq C_qs^{-n}\norm G_{C^q}
			(1+d(x,z)/s)^{-(q-1/2)}.
		\end{equation}
	\end{lemma}
	\begin{proof}
		Set $H(u)=e^{2u^2}G(u^2)$ on $\R$. Then
		$$
		G(s^2\Delta)=e^{-s^2\Delta}H(s\sqrt \Delta)e^{-s^2\Delta},\qquad
		\norm H_\infty\leq C\norm G_\infty.
		$$
		Write $h_z=h_{s^2}(\cdot,z)$, and take the $L^2$ inner product
		to be linear in its first variable. Then
		$$
		K_{G(s^2\Delta)}(x,z)=\langle H(s\sqrt \Delta)h_z,h_x\rangle,\qquad
		|K_{G(s^2\Delta)}(x,z)|
		\leq\norm H_\infty\norm{h_z}_2\norm{h_x}_2
		\leq Cs^{-n}\norm G_\infty.
		$$
		This is a continuous kernel. Indeed, the semigroup identity gives
		\begin{align*}
			\norm{h_x-h_{x'}}_2^2
			&=h_{2s^2}(x,x)+h_{2s^2}(x',x')
			-2h_{2s^2}(x,x')\longrightarrow0
			\quad(x'\longrightarrow x).
		\end{align*}
		For $q\geq1$, Plancherel and Cauchy--Schwarz give
		\begin{align}
			\int_\R(1+|u|)^{2q}|\widehat H(u)|^2\dd u
			&\leq C_q\norm G_{C^q}^2,\nonumber\\
			\int_{|u|>T}|\widehat H(u)|\dd u
			&\leq C_q\norm G_{C^q}(1+T)^{-(q-1/2)}.
			\label{eq:hw-Fourier-tail}
		\end{align}
		Use $\widehat H(u)=\int_{\R}e^{-iuv}H(v)\dd v$. Since $H$ is
		even,
		$$
		H(s\sqrt \Delta)=\frac1{2\pi}\int_{\R}\widehat H(u)
		\cos(us\sqrt \Delta)\dd u.
		$$
		Let $R=d(x,z)\geq4s$, and put
		$h_x^0=\mathbf1_{B(x,R/4)}h_x$, $h_x^1=h_x-h_x^0$, with the same
		notation at $z$. By \eqref{eq:hw-finite-propagation},
		$$
		\langle\cos(us\sqrt \Delta)h_z^0,h_x^0\rangle=0
		\quad\text{if }|u|<R/(2s).
		$$
		By \eqref{eq:hw-heat-column}--\eqref{eq:hw-heat-tail},
		$$
		\norm{h_x^1}_2\norm{h_z}_2+
		\norm{h_x^0}_2\norm{h_z^1}_2
		\leq Cs^{-n}e^{-cR^2/s^2}.
		$$
		Consequently,
		\begin{align*}
			|K_{G(s^2\Delta)}(x,z)|
			&\leq Cs^{-n}\left\{
			\int_{|u|>R/(2s)}|\widehat H(u)|\dd u
			+e^{-cR^2/s^2}\norm{\widehat H}_1\right\}\\
			&\leq C_qs^{-n}\norm G_{C^q}(1+R/s)^{-(q-1/2)}.
		\end{align*}
		For $R<4s$, \eqref{eq:hw-point-uniform} gives the same bound,
		since $(1+R/s)^{-(q-1/2)}\geq5^{-(q-1/2)}$.
	\end{proof}
	
	Choose relatively compact open sets $U\Subset U'\Subset M$.
	They contain the supports of the coefficients of the
	commutators, the boundaries of the exact ends, and the compact
	cutoff used below. Fix these sets before choosing $s$ and put
	$D(z)=d(z,U')$.
	\begin{lemma}\label{lem:hw-compact}
		Let $C_0$ be a differential operator of order at most one with
		smooth coefficients supported in $U$. Under the hypotheses of
		Lemma \ref{lem:hw-point},
		\begin{align}
			\int_U|C_0K_{G(s^2\Delta)}(w,z)|\dd\mu(w)
			&\leq Cs^{-n}\norm G_\infty,\label{eq:hw-compact-uniform}\\
			\int_U|C_0K_{G(s^2\Delta)}(w,z)|\dd\mu(w)
			&\leq C_qs^{-n}\norm G_{C^q}
			(1+D(z)/s)^{-(q-1/2)}.
			\label{eq:hw-compact-tail}
		\end{align}
		The constants are uniform for $s\geq s_*>0$.
	\end{lemma}
	\begin{proof}
		Choose $\eta\in C_c^\infty(U')$ with $\eta=1$ on $U$.
		For $v\in\operatorname{Dom} \Delta$, integration by parts and Young's inequality give
		\begin{align*}
			\int_M\eta^2|\nabla v|^2\dd\mu
			&=\Re\langle \Delta v,\eta^2v\rangle
			-2\Re\int_M\eta\overline v
			\nabla v\cdot\nabla\eta\dd\mu\\
			&\leq\norm{\Delta v}_{L^2(U')}\norm v_{L^2(U')}
			+\frac12\int_M\eta^2|\nabla v|^2\dd\mu
			+2\norm{\nabla\eta}_\infty^2\norm v_{L^2(U')}^2.
		\end{align*}
		Hence
		\begin{equation}\label{eq:hw-elliptic}
			\begin{aligned}
				\norm{\nabla v}_{L^2(U)}^2
				&\leq2\norm{\Delta v}_{L^2(U')}\norm v_{L^2(U')}
				+4\norm{\nabla\eta}_\infty^2\norm v_{L^2(U')}^2,\\
				\norm{C_0v}_{L^1(U)}
				&\leq C\bigl(\norm v_{L^2(U')}+\norm{\Delta v}_{L^2(U')}\bigr).
			\end{aligned}
		\end{equation}
		For the kernel column, the single heat factorization gives
		$$
		v=K_{G(s^2\Delta)}(\cdot,z)
		=[e^uG(u)](s^2\Delta)h_z\in\operatorname{Dom} \Delta,\qquad
		\Delta v=s^{-2}K_{(uG(u))(s^2\Delta)}(\cdot,z).
		$$
		The support interval and $s\geq s_*$ give
		$$
		\norm{uG(u)}_\infty\leq b_0\norm G_\infty,\qquad
		\norm{uG(u)}_{C^q}\leq C_{q,b_0}\norm G_{C^q},\qquad
		s^{-2}\leq s_*^{-2}.
		$$
		Lemma \ref{lem:hw-point}, applied to $G$ and $uG(u)$ on $U'$, yields
		\begin{align*}
			\norm v_{L^2(U')}+\norm{\Delta v}_{L^2(U')}
			&\leq Cs^{-n}\norm G_\infty,\\
			\norm v_{L^2(U')}+\norm{\Delta v}_{L^2(U')}
			&\leq C_qs^{-n}\norm G_{C^q}
			(1+D(z)/s)^{-(q-1/2)}.
		\end{align*}
		Substitution in \eqref{eq:hw-elliptic} proves both estimates.
		The derivative $C_0v$ is its weak derivative on $U$.
		For $f\in C_c^\infty(M)$, the kernel identity also holds after
		applying $C_0$:
		$$
		C_0G(s^2\Delta)f=\int_M C_0K_{G(s^2\Delta)}(\cdot,z)f(z)\dd\mu(z)
		\quad\text{in }L^2(U).
		$$
		Indeed, $\norm{C_0w}_{L^2(U)}
		\leq C(\norm w_2+\norm{\Delta^{1/2}w}_2)$, and
		$[e^uG(u)](s^2\Delta)$ maps $L^2$ boundedly into $\operatorname{Dom}(\Delta^{1/2})$.
		The single heat factorization therefore allows the integral
		to be taken in that form domain.
	\end{proof}
	
	We will use the following consequence with a fixed smooth
	$\theta$ supported in $[1/16,16]$. For
	$G_v(u)=e^{ivu}\theta(u)$, \eqref{eq:hw-compact-uniform} is independent
	of $v$, while \eqref{eq:hw-compact-tail} with $q=2$ gives
	$$
	\int_U|C_0K_{G_v(s^2\Delta)}(w,z)|\dd\mu(w)
	\leq Cs^{-n}\min\{1,(1+|v|)^2(1+D(z)/s)^{-3/2}\}.
	$$
	Since $\min\{1,t\}\leq t^{1/12}$ for $t\geq0$,
	\begin{equation}\label{eq:hw-phase-compact}
		\int_U|C_0K_{G_v(s^2\Delta)}(w,z)|\dd\mu(w)
		\leq Cs^{-n}(1+|v|)^{1/6}(1+D(z)/s)^{-1/8}.
	\end{equation}
	
	\subsection{The two product models}
	Denote the nonnegative model Laplacians by $\Delta_m,\Delta_n$.
	For either model, let
	$V_i(t)$ be the volume of a ball of radius $t$. Then
	$$
	V_m(t)\approx t^m,\qquad
	V_n(t)\approx\begin{cases}t^m,&t\leq1,\\t^n,&t\geq1,\end{cases}
	\qquad V_i(R)\leq C(R/t)^mV_i(t)\quad(R\geq t).
	$$
	\begin{lemma}\label{lem:hw-model-bands}
		For an integer $q>(m+1)/2$ and a fixed support interval
		$[0,b_0]$, the model kernels satisfy
		\begin{align}
			\sup_z\int_{M_i}|K_{G(t^2\Delta_i)}(x,z)|\dd\mu_i(x)
			&\leq C_q\norm G_{C^q},\label{eq:hw-model-L1}\\
			\sup_z\int_{d_i(x,z)>R}|K_{G(t^2\Delta_i)}(x,z)|\dd\mu_i(x)
			&\leq C_q\norm G_{C^q}(t/R)^{q-1/2-m/2},
			\quad R\geq t.
			\label{eq:hw-model-tail}
		\end{align}
		Here $t>0$, and $G$ is a $C^q$ function supported in $[0,b_0]$.
	\end{lemma}
	\begin{proof}
		The model heat kernels satisfy
		$$
		h^{(i)}_{t^2}(x,z)\leq CV_i(t)^{-1}e^{-cd_i(x,z)^2/t^2}.
		$$
		Their mass is at most one. Consequently their heat columns
		obey
		$$
		\norm{h^{(i)}_{t^2}(\cdot,z)}_2\leq CV_i(t)^{-1/2},\qquad
		\norm{\mathbf1_{\{d_i(\cdot,z)>R\}}h^{(i)}_{t^2}(\cdot,z)}_2
		\leq CV_i(t)^{-1/2}e^{-cR^2/t^2}.
		$$
		Use
		$$
		G(t^2\Delta_i)=H(t\sqrt{\Delta_i})e^{-t^2\Delta_i},\qquad
		H(u)=e^{u^2}G(u^2).
		$$
		For $R\geq4t$, put $h^0=\mathbf1_{B_i(z,R/4)}h^{(i)}_{t^2}(\cdot,z)$
		and $h^1=h^{(i)}_{t^2}(\cdot,z)-h^0$. Finite propagation gives
		$$
		\mathbf1_{\{d_i(\cdot,z)>R\}}\cos(ut\sqrt{\Delta_i})h^0=0
		\quad(|u|<R/(2t)).
		$$
		Fourier inversion and \eqref{eq:hw-Fourier-tail} therefore give
		\begin{align*}
			\norm{\mathbf1_{\{d_i(\cdot,z)>R\}}K_{G(t^2\Delta_i)}(\cdot,z)}_2&\leq CV_i(t)^{-1/2}
			\left(\int_{|u|>R/(2t)}|\widehat H(u)|\dd u
			+e^{-cR^2/t^2}\norm{\widehat H}_1\right)\\
			&\leq C_qV_i(t)^{-1/2}\norm G_{C^q}
			(1+R/t)^{-(q-1/2)}.
		\end{align*}
		For $t\leq R<4t$ the same estimate follows from
		$\norm{K_{G(t^2\Delta_i)}(\cdot,z)}_2
		\leq\norm H_\infty\norm{h^{(i)}_{t^2}(\cdot,z)}_2$.
		Thus, for all $R\geq t$,
		\begin{equation}\label{eq:hw-model-L2-tail}
			\norm{\mathbf1_{\{d_i(\cdot,z)>R\}}K_{G(t^2\Delta_i)}(\cdot,z)}_2
			\leq C_qV_i(t)^{-1/2}\norm G_{C^q}
			(1+R/t)^{-(q-1/2)}.
		\end{equation}
		Put $p_q=q-1/2-m/2>0$. Cauchy--Schwarz on the annuli gives
		\begin{align*}
			\int_{d_i(x,z)>R}|K_{G(t^2\Delta_i)}(x,z)|\dd\mu_i(x)
			&\leq\sum_{j=0}^\infty V_i(2^{j+1}R)^{1/2}
			\norm{\mathbf1_{\{d_i(\cdot,z)>2^jR\}}K_{G(t^2\Delta_i)}(\cdot,z)}_2\\
			&\leq C_q\norm G_{C^q}(t/R)^{p_q}
			\sum_{j=0}^\infty2^{-jp_q}\\
			&=\frac{C_q}{1-2^{-p_q}}\norm G_{C^q}(t/R)^{p_q}.
		\end{align*}
		On the remaining ball,
		$$
		\int_{B_i(z,t)}|K_{G(t^2\Delta_i)}(x,z)|\dd\mu_i(x)
		\leq V_i(t)^{1/2}\norm{K_{G(t^2\Delta_i)}(\cdot,z)}_2
		\leq C\norm G_\infty.
		$$
		These are \eqref{eq:hw-model-tail} and \eqref{eq:hw-model-L1}.
	\end{proof}
	
	Let $\nu_1>0$ be the first nonzero eigenvalue of the Laplacian
	on $S^{m-n}$.
	\begin{lemma}\label{lem:hw-model-phase}
		For a fixed $\theta\in C_c^\infty((1/16,16))$ and all
		$s>4/\sqrt{\nu_1}$,
		\begin{equation}\label{eq:hw-model-phase}
			\norm{e^{iu s^2\Delta_i}\theta(s^2\Delta_i)}_{L^1(M_i)\to L^1(M_i)}
			\leq C(1+|u|)^{m/2},\qquad u\in\R.
		\end{equation}
	\end{lemma}
	\begin{proof}
		For $M_n$, our choice of $s$ gives $16/s^2<\nu_1$.
		Writing $\Delta_n=-\sum_{k=1}^n\partial_{x_k}^2+\Delta_{S^{m-n}}$, we have
		$$
		\theta\bigl(s^2(|\xi|^2+\nu_j)\bigr)=0\quad(j\geq1),\qquad
		K_{\Pi_0}(\omega,\omega')=\mu(S^{m-n})^{-1},\qquad
		\norm{\Pi_0}_{1\to1}=1.
		$$
		It remains to work on $\R^d$, $d=m$ or $n$. Put
		$$
		k_u=\mathcal F^{-1}\!\left[e^{iu|\xi|^2}\theta(|\xi|^2)\right].
		$$
		The kernel at length $s$ is $s^{-d}k_u((x-z)/s)$, so its
		$L^1$ norm is $\norm{k_u}_1$. For each integer $q$,
		Plancherel gives
		\begin{align*}
			\norm{k_u}_2&\leq C,\\
			\norm{|x|^qk_u}_2
			&\leq C_q\sum_{|\gamma|=q}
			\norm{\partial_\xi^\gamma
				(e^{iu|\xi|^2}\theta(|\xi|^2))}_2
			\leq C_q(1+|u|)^q.
		\end{align*}
		Choose $q>d/2$ and set $R=1+|u|$. Cauchy--Schwarz gives
		\begin{align*}
			\norm{k_u}_1
			&\leq\left(\int_{\R^d}(1+|x|/R)^{-2q}\dd x\right)^{1/2}
			\norm{(1+|x|/R)^q k_u}_2\\
			&\leq C_qR^{d/2}
			\bigl(\norm{k_u}_2+R^{-q}\norm{|x|^qk_u}_2\bigr)
			\leq C_qR^{d/2}\leq C_qR^{m/2}.
		\end{align*}
		The proof is complete.
	\end{proof}
	
	\subsection{The compact commutator comparison}
	Choose $0\leq\chi_i\leq1$, $i=m,n$, supported in the exact
	part of $E_i$ and equal to one outside a larger compact set.
	Extend $\chi_if$ by zero into the missing part of $M_i$, and put
	$$
	X_if=\chi_if,\qquad
	C_if=(\Delta\chi_i)f-2\nabla\chi_i\cdot\nabla f.
	$$
	The coefficients of $C_i$ are supported in $U\cap E_i$, where
	the metrics agree. Thus Lemma \ref{lem:hw-compact} applies, and
	\begin{equation}\label{eq:hw-commutator}
		\norm{C_if}_2\leq C\bigl(\norm f_2+\norm{\Delta^{1/2}f}_2\bigr),
		\qquad \norm{X_i}_{1\to1},\norm{X_i}_{2\to2}\leq1.
	\end{equation}
	The adjoint $X_i^*$ is restriction to the exact end followed by
	multiplication by $\chi_i$.
	
	Here is the domain identity. The energy forms are
	$$
	\mathfrak q(f,g)=\int_M\nabla f\cdot\nabla\overline g\dd\mu,
	\qquad
	\mathfrak q_i(v,w)=\int_{M_i}\nabla v\cdot\nabla\overline w\dd\mu_i.
	$$
	The dot denotes the metric extended bilinearly. Boundedness of
	$\chi_i$ and $\nabla\chi_i$ gives
	$$
	X_i\operatorname{Dom}(\Delta^{1/2})\subset\operatorname{Dom}(\Delta_i^{1/2}),\qquad
	X_i^*\operatorname{Dom}(\Delta_i^{1/2})\subset\operatorname{Dom}(\Delta^{1/2}).
	$$
	For $f\in\operatorname{Dom} \Delta$ and $\varphi\in\operatorname{Dom}(\Delta_i^{1/2})$, integration by
	parts gives
	\begin{align*}
		\mathfrak q_i(X_if,\varphi)-\mathfrak q(f,X_i^*\varphi)
		&=\int\bigl(f\nabla\chi_i\cdot\nabla\overline\varphi
		-\overline\varphi\nabla f\cdot\nabla\chi_i\bigr)\dd\mu\\
		&=\int\bigl((\Delta\chi_i)f
		-2\nabla\chi_i\cdot\nabla f\bigr)
		\overline\varphi\dd\mu
		=\langle C_if,\varphi\rangle.
	\end{align*}
	The integrals have compact support in the common part of the two
	manifolds. The definition of the operator associated with a
	closed form therefore gives
	\begin{equation}\label{eq:hw-commutator-domain}
		\mathfrak q_i(X_if,\varphi)
		=\langle X_i\Delta f+C_if,\varphi\rangle,\qquad
		X_if\in\operatorname{Dom} \Delta_i,\qquad \Delta_iX_if-X_i\Delta f=C_if.
	\end{equation}
	The form identities follow first for smooth functions and then
	by approximation in the form norms.
	
	Fix $\psi_0\in C_c^\infty([0,\infty))$ with
	$$
	\psi_0=1\text{ on }[0,1],\qquad \supp\psi_0\subset[0,2].
	$$
	For this fixed operator and every $A_0>0$,
	\begin{align}
		\sup_z\int_M|K_{\psi_0(\Delta)}(w,z)|\dd\mu(w)&\leq C,
		\label{eq:hw-fixed-L1}\\
		\sup_z\int_{d(w,z)>R}|K_{\psi_0(\Delta)}(w,z)|\dd\mu(w)
		&\leq C_{A_0}R^{-A_0},\qquad R\geq1.
		\label{eq:hw-fixed-tail}
	\end{align}
	To verify these bounds, choose an integer $q>A_0+(m+1)/2$ and
	set $H(u)=e^{u^2}\psi_0(u^2)$. The heat column split and finite
	propagation in the proof of \eqref{eq:hw-model-L2-tail}, with
	\eqref{eq:hw-bounded-time-heat} at $s=1$, give
	$$
	\norm{K_{\psi_0(\Delta)}(\cdot,z)}_2\leq C,\qquad
	\norm{\mathbf1_{\{d(\cdot,z)>R\}}K_{\psi_0(\Delta)}(\cdot,z)}_2
	\leq C_q(1+R)^{-(q-1/2)}.
	$$
	Put $p_q=q-1/2-m/2>A_0$. By \eqref{eq:volume-upper},
	\begin{align*}
		\int_{d(w,z)>R}|K_{\psi_0(\Delta)}(w,z)|\dd\mu(w)
		&\leq C_q\sum_{j=0}^\infty(2^{j+1}R)^{m/2}
		(2^jR)^{-(q-1/2)}
		\leq\frac{C_q}{1-2^{-p_q}}R^{-p_q}\leq C_{A_0}R^{-A_0}.
	\end{align*}
	The integral on $B(z,1)$ is at most
	$\mu(B(z,1))^{1/2}\norm{K_{\psi_0(\Delta)}(\cdot,z)}_2\leq C$.
	
	For all sufficiently large $s$, put $F_s(u)=F(s^2u)$ and
	$$
	D_{i,s}=F_s(\Delta_i)X_i\psi_0(\Delta)-X_iF_s(\Delta).
	$$
	This maps functions on $M$ to functions on $M_i$.
	\begin{proposition}\label{prop:hw-comparison}
		There is $s_0<\infty$ such that, for $s\geq s_0$,
		\begin{equation}\label{eq:hw-comparison-column}
			\int_{M_i}|K_{D_{i,s}}(x,z)|\dd\mu_i(x)
			\leq C_N\norm F_{C^N}(1+D(z)/s)^{-1/8}.
		\end{equation}
	\end{proposition}
	
	\subsubsection{A formula for disjoint spectral intervals}
	Let $A,B$ be nonnegative self-adjoint operators and $X$ a bounded
	map with $X\operatorname{Dom} B\subset\operatorname{Dom} A$. Suppose that
	$C=AX-XB$ is bounded after restriction to each bounded
	$B$-spectral interval. These hypotheses hold by
	\eqref{eq:hw-commutator}--\eqref{eq:hw-commutator-domain}.
	Restrict both variables to bounded spectral intervals. If
	$P(A)$ is supported in $0\leq\lambda\leq a$ and $Q(B)$ in
	$\mu\geq b>0$, where $a/b\leq1/2$, then
	\begin{equation}\label{eq:hw-separated-low-high}
		P(A)XQ(B)=-\sum_{k=0}^\infty A^kP(A)C B^{-k-1}Q(B).
	\end{equation}
	If $P(A)$ is supported in $\lambda\geq b$ and $Q(B)$ in
	$0\leq\mu\leq a$, then
	\begin{equation}\label{eq:hw-separated-high-low}
		P(A)XQ(B)=\sum_{k=0}^\infty A^{-k-1}P(A)C B^kQ(B).
	\end{equation}
	Each inverse power is multiplied by its cutoff away from zero
	and is defined to be zero elsewhere.
	
	Put $Y=P(A)XQ(B)$ and $C'=P(A)CQ(B)$. The identity
	$AY-YB=C'$ gives, in the first case,
	$$
	Y=AYB^{-1}-C'B^{-1}
	=-\sum_{k=0}^{K-1}A^kC'B^{-k-1}+A^KYB^{-K}.
	$$
	In the second case,
	$$
	Y=A^{-1}C'+A^{-1}YB
	=\sum_{k=0}^{K-1}A^{-k-1}C'B^k+A^{-K}YB^K.
	$$
	The two remainders satisfy, respectively,
	$$
	\norm{A^KYB^{-K}}_{2\to2}\leq(a/b)^K\norm Y_{2\to2},\qquad
	\norm{A^{-K}YB^K}_{2\to2}\leq(a/b)^K\norm Y_{2\to2}.
	$$
	Let $K\to\infty$. Upper intervals extending to infinity will be
	split into bounded intervals below.
	
	\subsubsection{The intervals that meet}
	Choose $\theta\in C_c^\infty((1/16,16))$ with
	$0\leq\theta\leq1$ and $\theta=1$ on $[1/8,8]$. Set
	$$
	R_U=\sup_{w\in U'}\rho(w),\qquad
	s_0\geq\max\{4,R_U+1,8/\sqrt{\nu_1}\}.
	$$
	For $s\geq s_0$, we have $16/s^2\leq1$ and $16/s^2<\nu_1$.
	Hence $\psi_0(\Delta)\theta(s^2\Delta)=\theta(s^2\Delta)$, and
	Lemma \ref{lem:hw-model-phase} applies.
	Write $A=\Delta_i$, $B=\Delta$, $X=X_i$, and $C=C_i$ until the end of
	this section. Put $\theta_A=\theta(s^2A)$ and
	$\theta_B=\theta(s^2B)$.
	
	The domain identity \eqref{eq:hw-commutator-domain} gives
	\begin{align*}
		\frac{\dd}{\dd v}
		\left(e^{i(u-v)s^2A}\theta_AX\theta_Be^{iv s^2B}\right)
		&=-is^2e^{i(u-v)s^2A}\theta_A C\theta_Be^{iv s^2B},\\
		\theta_A\bigl(e^{iu s^2A}X-Xe^{iu s^2B}\bigr)\theta_B
		&=is^2\int_0^u e^{i(u-v)s^2A}\theta_A
		C\theta_Be^{iv s^2B}\dd v.
	\end{align*}
	All differentiated factors have bounded spectral support.
	Since $F(\lambda)=(2\pi)^{-1}\int_{\R}\widehat F(u)e^{iu\lambda}\dd u$,
	\begin{align}
		&\theta_A\bigl(F_s(A)X-XF_s(B)\bigr)\theta_B\nonumber\\
		&\quad=\frac{is^2}{2\pi}\int_\R\widehat F(u)
		\int_0^u e^{i(u-v)s^2A}\theta_A
		C e^{iv s^2B}\theta_B\dd v\dd u.
		\label{eq:hw-Duhamel}
	\end{align}
	For $u<0$, the inner integral is oriented from $0$ to $u$.
	For fixed $s$, its $L^2$ operator norm is at most
	$$
	Cs^2\norm{C\theta_B}_{2\to2}
	\int_{\R}|u|\,|\widehat F(u)|\dd u<\infty.
	$$
	
	By Lemma \ref{lem:hw-model-phase} and \eqref{eq:hw-phase-compact},
	the $L^1$ norm of the kernel column in \eqref{eq:hw-Duhamel} is at most
	\begin{equation}\label{eq:hw-Duhamel-L1}
		Cs^{2-n}(1+D(z)/s)^{-1/8}
		\int_\R|\widehat F(u)|\,|u|
		(1+|u|)^{m/2+1/6}\dd u.
	\end{equation}
	Indeed, for $v$ between $0$ and $u$,
	$$
	(1+|u-v|)^{m/2}(1+|v|)^{1/6}
	\leq(1+|u|)^{m/2+1/6}.
	$$
	The left factor acts on the model; $C$ has compact support.
	For integers $m,N$ with $N>(m+3)/2$, we have $N-m/2\geq2$.
	Plancherel and Cauchy--Schwarz give
	\begin{align*}
		\int_\R |\widehat F(u)|(1+|u|)^{m/2+7/6}\dd u
		&\leq\norm{(1+|u|)^N\widehat F(u)}_2
		\left(\int_\R(1+|u|)^{2(m/2+7/6-N)}\dd u\right)^{1/2}\\
		&\leq C_N\norm F_{C^N}
		\left(\int_\R(1+|u|)^{-5/3}\dd u\right)^{1/2}\\
		&\leq C_N\norm F_{C^N}.
	\end{align*}
	Since $s^{2-n}\leq1$, \eqref{eq:hw-Duhamel-L1} has the bound in
	\eqref{eq:hw-comparison-column}.
	
	\subsubsection{The four disjoint cases}
	Split $1-\theta(u)=\alpha_-(u)+\alpha_+(u)$ on $u\geq0$, where
	$$
	\supp\alpha_-\subset[0,1/8],\qquad
	\supp\alpha_+\subset[8,\infty),\qquad
	\alpha_-=1\text{ near }0,\qquad \alpha_+=1\text{ on }[16,\infty).
	$$
	Since $\theta F=F$ and $\psi_0(B)\theta_B=\theta_B$,
	\begin{align}
		D_{i,s}
		&=\theta_AD_{i,s}\theta_B+(I-\theta_A)D_{i,s}
		+\theta_AD_{i,s}(I-\theta_B)\nonumber\\
		&=\theta_A(F_s(A)X-XF_s(B))\theta_B\nonumber\\
		&\quad-\alpha_-(s^2A)XF_s(B)
		-\alpha_+(s^2A)XF_s(B)\nonumber\\
		&\quad+F_s(A)X\alpha_-(s^2B)
		+F_s(A)X\alpha_+(s^2B)\psi_0(B).
		\label{eq:hw-four-parts}
	\end{align}
	Here $\psi_0(B)\alpha_-(s^2B)=\alpha_-(s^2B)$.
	
	The model estimate uses $N$ derivatives and the compact estimate
	uses two. We keep the resulting factor $(1+k)^{N+2}$ explicit.
	On a fixed interval $[a,b]\Subset(0,\infty)$, the Leibniz formula gives
	\begin{align*}
		\norm{u^kG(u)}_{C^q}
		&\leq C_q(1+k)^q b^k\norm G_{C^q},\\
		\norm{u^{-k-1}G(u)}_{C^q}
		&\leq C_q(1+k)^q a^{-k-1}\norm G_{C^q}.
	\end{align*}
	Also,
	$$
	\norm{u^k\alpha_-(u)}_{C^q}\leq C_q(1+k)^q8^{-k},\qquad
	\sum_{k=0}^\infty(1+k)^{N+2}2^{-k}<\infty.
	$$
	Constants such as $a^{-q}$ and $8^q$ do not depend on $k$.
	
	For the two upper intervals, choose a smooth nonincreasing
	$\zeta$ equal to one on $[0,16]$ and zero on $[32,\infty)$.
	Put
	$$
	\eta_0(u)=\zeta(u),\qquad
	\eta_h(u)=\zeta(2^{-h}u)-\zeta(2^{1-h}u)\quad(h\geq1),\qquad
	\beta_h(\lambda)=\alpha_+(s^2\lambda)\eta_h(s^2\lambda).
	$$
	Telescoping gives
	$$
	\sum_{h=0}^H\beta_h(\lambda)
	=\alpha_+(s^2\lambda)\zeta(2^{-H}s^2\lambda)
	\longrightarrow\alpha_+(s^2\lambda).
	$$
	Set
	\begin{equation}\label{eq:hw-packet-endpoints}
		\ell_h=8\cdot2^hs^{-2},\qquad
		\nu_h=16\cdot2^hs^{-2}=2\ell_h.
	\end{equation}
	Then $\supp\beta_h\subset[\ell_h,4\ell_h]$ and, for each fixed $q$,
	$$
	\norm{\beta_h(\nu_hu)}_{C^q}\leq C_q,\qquad
	\norm{(\nu_hu)^{-k-1}\beta_h(\nu_hu)}_{C^q}
	\leq C_q(1+k)^q\ell_h^{-k-1}.
	$$
	For $h\geq1$, $\alpha_+$ is already one on the support of
	$\eta_h$; for $h=0$, all its arguments lie in a fixed interval.
	
	\medskip\noindent\emph{First case: $A$ is below the band.}
	By \eqref{eq:hw-separated-low-high}, the terms, apart from a sign, are
	$$
	A^k\alpha_-(s^2A)\,C\,B^{-k-1}F_s(B).
	$$
	Lemmas \ref{lem:hw-model-bands} and \ref{lem:hw-compact} give
	\begin{align*}
		\norm{A^k\alpha_-(s^2A)}_{1\to1}
		&\leq C_N(1+k)^Ns^{-2k}8^{-k},\\
		\int_U|K_{CB^{-k-1}F_s(B)}(w,z)|\dd\mu(w)
		&\leq C(1+k)^2s^{2k+2-n}4^{k+1}\norm F_{C^N}
		(1+D(z)/s)^{-1/8}.
	\end{align*}
	Their product is at most
	\begin{equation}\label{eq:hw-first-case}
		C_N(1+k)^{N+2}2^{-k}s^{2-n}\norm F_{C^N}
		(1+D(z)/s)^{-1/8}.
	\end{equation}
	The compact estimate with $q=2$ has decay $3/2$, which implies
	the displayed decay $1/8$. Sum in $k$.
	
	\medskip\noindent\emph{Second case: $A$ is above the band.}
	Use \eqref{eq:hw-separated-high-low} on each $\beta_h$.
	The terms are
	$$
	A^{-k-1}\beta_h(A)\,C\,B^kF_s(B).
	$$
	The two factors satisfy
	\begin{align*}
		\norm{A^{-k-1}\beta_h(A)}_{1\to1}
		&\leq C_N(1+k)^N\ell_h^{-k-1},\\
		\int_U|K_{CB^kF_s(B)}(w,z)|\dd\mu(w)
		&\leq C(1+k)^2s^{-n}(4s^{-2})^k\norm F_{C^N}
		(1+D(z)/s)^{-1/8}.
	\end{align*}
	Since $4s^{-2}/\ell_h=2^{-h-1}\leq1/2$, their product is at most
	$$
	C_N(1+k)^{N+2}2^{-k}\ell_h^{-1}s^{-n}\norm F_{C^N}
	(1+D(z)/s)^{-1/8}.
	$$
	Hence the sum in $k,h$ is at most
	$C_Ns^{2-n}\norm F_{C^N}(1+D(z)/s)^{-1/8}$, because
	$$
	\sum_{h=0}^\infty\ell_h^{-1}
	=\frac{s^2}{8}\sum_{h=0}^\infty2^{-h}=\frac{s^2}{4}.
	$$
	
	\medskip\noindent\emph{Third case: $B$ is below the band.}
	By \eqref{eq:hw-separated-high-low}, the terms are
	$$
	A^{-k-1}F_s(A)\,C\,B^k\alpha_-(s^2B).
	$$
	The two factors satisfy
	\begin{align*}
		\norm{A^{-k-1}F_s(A)}_{1\to1}
		&\leq C_N(1+k)^Ns^{2k+2}4^{k+1}\norm F_{C^N},\\
		\int_U|K_{CB^k\alpha_-(s^2B)}(w,z)|\dd\mu(w)
		&\leq C(1+k)^2s^{-2k-n}8^{-k}(1+D(z)/s)^{-1/8}.
	\end{align*}
	Their product has the bound \eqref{eq:hw-first-case}, including
	the whole interval at zero. Sum in $k$.
	
	\medskip\noindent\emph{Fourth case: $B$ is above the band.}
	Set $\gamma_h(\mu)=\beta_h(\mu)\psi_0(\mu)$.
	If $\gamma_h\not\equiv0$, then $\ell_h\leq2$, hence $\nu_h\leq4$.
	For each fixed $q$,
	\begin{align*}
		\partial_u^q\gamma_h(\nu_hu)
		&=\sum_{j=0}^q\binom qj
		\partial_u^{q-j}\beta_h(\nu_hu)\,
		\nu_h^j\psi_0^{(j)}(\nu_hu),\\
		\norm{\gamma_h(\nu_hu)}_{C^q}&\leq C_q,\\
		\norm{(\nu_hu)^{-k-1}\gamma_h(\nu_hu)}_{C^2}
		&\leq C(1+k)^2\ell_h^{-k-1}.
	\end{align*}
	By \eqref{eq:hw-separated-low-high}, the terms are
	$$
	A^kF_s(A)\,C\,B^{-k-1}\gamma_h(B).
	$$
	Lemma \ref{lem:hw-compact} applies with $s=\nu_h^{-1/2}\geq1/2$.
	The two factors satisfy
	\begin{align*}
		\norm{A^kF_s(A)}_{1\to1}
		&\leq C_N(1+k)^N(4s^{-2})^k\norm F_{C^N},\\
		\int_U|K_{CB^{-k-1}\gamma_h(B)}(w,z)|\dd\mu(w)
		&\leq C(1+k)^2\nu_h^{n/2}\ell_h^{-k-1}
		(1+D(z)\sqrt{\nu_h})^{-1/8}.
	\end{align*}
	Since $\nu_h=2\ell_h$, $4s^{-2}/\ell_h\leq1/2$ and
	$\sqrt{\nu_h}\geq4/s$, their product is at most
	$$
	C_N(1+k)^{N+2}2^{-k}\norm F_{C^N}\,
	\nu_h^{n/2-1}(1+D(z)/s)^{-1/8}.
	$$
	The decisive sum is
	\begin{equation}\label{eq:hw-transient-sum}
		\sum_{\substack{h\geq0\\\nu_h\leq4}}\nu_h^{n/2-1}
		\leq\frac{4^{n/2-1}}{1-2^{-(n/2-1)}},\qquad n>2.
	\end{equation}
	This follows by summing backwards from the largest $\nu_h$;
	the ratio is $2^{-(n/2-1)}<1$. Sum in $h$ and $k$.
	
	\subsubsection{Convergence and the kernel identity}
	Each series in $k$ converges in $L^2$ operator norm by the
	remainder estimates, and its kernel columns converge absolutely
	in $L^1$ by the displayed bounds. In the second case,
	$$
	\sum_{h=0}^H\beta_h(A)XF_s(B)
	=\alpha_+(s^2A)\zeta(2^{-H}s^2A)XF_s(B)
	\longrightarrow\alpha_+(s^2A)XF_s(B)
	\quad\text{strongly in }L^2.
	$$
	These sums also converge in $L^1$ operator norm, since
	$\sum_h\ell_h^{-1}<\infty$. The fourth case has only finitely
	many nonzero terms in $h$ for each $s$.
	
	The integral \eqref{eq:hw-Duhamel} converges in $L^2$ operator norm
	and in the supremum of the $L^1$ column norms by
	\eqref{eq:hw-Duhamel-L1}. For $f\in L^1\cap L^2$, let $g_1,g_2$
	be the two limits of any of these sums or integrals. For a set
	$\Omega$ of finite measure and a common truncation $T_Hf$,
	$$
	\norm{g_1-g_2}_{L^1(\Omega)}
	\leq\norm{g_1-T_Hf}_1+
	\mu_i(\Omega)^{1/2}\norm{T_Hf-g_2}_2\longrightarrow0.
	$$
	Thus the kernel formulas agree with their spectral identities.
	
	The kernel of $D_{i,s}$ is continuous. Its first term is the
	$L^2$ pairing of a model kernel row with
	$X_iK_{\psi_0(\Delta)}(\cdot,z)$; its second term is
	$\chi_i(x)K_{F_s(\Delta)}(x,z)$, extended by zero. The heat
	factorizations give continuous $L^2$ rows and columns, as in
	Lemma \ref{lem:hw-point}. Consequently, for $z_\ell\to z_0$,
	$$
	\int_{M_i}|K_{D_{i,s}}(x,z_0)|\dd\mu_i(x)
	\leq\liminf_{\ell\to\infty}
	\int_{M_i}|K_{D_{i,s}}(x,z_\ell)|\dd\mu_i(x).
	$$
	Choose $z_\ell$ where the almost everywhere bound holds and use
	the continuity of $D(z)$. This proves the bound at every $z_0$.
	Equation \eqref{eq:hw-four-parts} proves Proposition
	\ref{prop:hw-comparison}.
	
	\subsection{The global band estimates}
	\subsubsection{Bounded values of \texorpdfstring{$s$}{s}}
	For $0<s\leq s_0$, \eqref{eq:hw-bounded-time-heat} and
	\eqref{eq:hw-mass} give
	$$
	\norm{h_{s^2}(\cdot,z)}_2\leq C_{s_0}s^{-m/2},\qquad
	\norm{\mathbf1_{\{d(\cdot,z)>R\}}h_{s^2}(\cdot,z)}_2
	\leq C_{s_0}s^{-m/2}e^{-cR^2/s^2}.
	$$
	Use $F(s^2\Delta)=H(s\sqrt \Delta)e^{-s^2\Delta}$, with
	$H(u)=e^{u^2}F(u^2)$. For $R\geq4s$, split the heat column
	at $R/4$. For $s\leq R<4s$, use its full $L^2$ norm.
	Exactly the two estimates in the proof of
	\eqref{eq:hw-model-L2-tail} give
	\begin{align*}
		\norm{K_{F(s^2\Delta)}(\cdot,z)}_2
		&\leq C_{s_0}s^{-m/2}\norm F_{C^N},\\
		\norm{\mathbf1_{\{d(\cdot,z)>R\}}K_{F(s^2\Delta)}(\cdot,z)}_2
		&\leq C_{N,s_0}s^{-m/2}\norm F_{C^N}(1+R/s)^{-(N-1/2)}.
	\end{align*}
	Put $p_N=N-1/2-m/2>1/8$. By \eqref{eq:volume-upper},
	\begin{align*}
		\int_{d(x,z)>R}|K_{F(s^2\Delta)}(x,z)|\dd\mu(x)
		&\leq\sum_{j=0}^\infty C(2^{j+1}R)^{m/2}
		\norm{\mathbf1_{\{d(\cdot,z)>2^jR\}}K_{F(s^2\Delta)}(\cdot,z)}_2\\
		&\leq\frac{C_{N,s_0}}{1-2^{-p_N}}\norm F_{C^N}(s/R)^{p_N},
		\qquad R\geq s,\\
		\int_{B(z,s)}|K_{F(s^2\Delta)}(x,z)|\dd\mu(x)
		&\leq Cs^{m/2}\norm{K_{F(s^2\Delta)}(\cdot,z)}_2
		\leq C_{s_0}\norm F_{C^N}.
	\end{align*}
	Thus \eqref{eq:mihlin-6} holds. For \eqref{eq:mihlin-7}, use
	$$
	z\in B,\ x\notin2B\quad\Longrightarrow\quad d(x,z)>r>s,
	\qquad (s/r)^{p_N}\leq(s/r)^{1/8}.
	$$
	
	\subsubsection{Large values of \texorpdfstring{$s$}{s}}
	Let $\chi_0=1-\chi_m-\chi_n$. The end cutoffs have disjoint
	supports, so $0\leq\chi_0\leq1$, and $\supp\chi_0\subset U'$.
	Lemma \ref{lem:hw-point}, with $q=2$, gives
	\begin{equation}\label{eq:hw-compact-part}
		\int_M\chi_0(x)|K_{F(s^2\Delta)}(x,z)|\dd\mu(x)
		\leq C_N\norm F_{C^N} s^{-n}(1+D(z)/s)^{-1/8}.
	\end{equation}
	On each model the identity is
	\begin{equation}\label{eq:hw-comparison-identity}
		X_iF_s(\Delta)=F_s(\Delta_i)X_i\psi_0(\Delta)-D_{i,s}.
	\end{equation}
	By Fubini, \eqref{eq:hw-model-L1}, and \eqref{eq:hw-fixed-L1},
	\begin{align*}
		\int_{M_i}|K_{F_s(\Delta_i)X_i\psi_0(\Delta)}(x,z)|\dd\mu_i(x)
		&\leq\int_{E_i}\chi_i(w)|K_{\psi_0(\Delta)}(w,z)|\\
		&\qquad\times
		\left(\int_{M_i}|K_{F_s(\Delta_i)}(x,w)|\dd\mu_i(x)\right)
		\dd\mu(w)\\
		&\leq C_N\norm F_{C^N}.
	\end{align*}
	Together with Proposition \ref{prop:hw-comparison} and
	$\chi_0+\chi_m+\chi_n=1$, this proves \eqref{eq:mihlin-6}.
	
	Now let $s\geq s_0$, $s<r\leq\delta\rho(x_B)$, and $z\in B$.
	If $r\leq2s$, then
	$$
	\int_{M\setminus2B}|K_{F(s^2\Delta)}(x,z)|\dd\mu(x)
	\leq C_N\norm F_{C^N}
	\leq2^{1/8}C_N\norm F_{C^N}(s/r)^{1/8}.
	$$
	Suppose $r>2s$. Our fixed choice $s_0\geq R_U+1$ gives
	\begin{equation}\label{eq:hw-distance-to-core}
		\begin{aligned}
			D(z)&\geq\rho(z)-R_U\geq\rho(x_B)-r-R_U\\
			&\geq\rho(x_B)-\tfrac32r
			\geq(1-\tfrac32\delta)\rho(x_B)
			\geq\tfrac12\rho(x_B)\geq r/(2\delta)\geq8r.
		\end{aligned}
	\end{equation}
	The compact part and the comparison errors satisfy
	$$
	Cs^{-n}(1+D(z)/s)^{-1/8}
	+C(1+D(z)/s)^{-1/8}\leq C(s/r)^{1/8}.
	$$
	
	It remains to estimate the first term of
	\eqref{eq:hw-comparison-identity}. Split its variable $w$ at
	$d(w,z)=r/2$. The part with $d(w,z)>r/2$ has column norm at most
	\begin{align*}
		&C_N\norm F_{C^N}\int_{d(w,z)>r/2}|K_{\psi_0(\Delta)}(w,z)|\dd\mu(w)\leq C_N\norm F_{C^N} r^{-1/8}
		\leq C_N\norm F_{C^N}(s/r)^{1/8},
	\end{align*}
	by \eqref{eq:hw-fixed-tail} and $s\geq1$.
	
	For the other part, take $w\in E_i$ with $d(w,z)\leq r/2$.
	Then
	$$
	D(w)\geq D(z)-r/2\geq\tfrac{15}{2}r,\qquad
	d_i(w,\partial E_i)=d(w,\partial E_i)\geq D(w)>r/2.
	$$
	The equality holds because a path to the boundary may be stopped
	at its first meeting with the boundary; before that meeting the
	metrics agree. Any model path of length less than $r/2$ from
	$w$ therefore stays in $E_i$ and is a path in $M$. Consequently,
	\begin{align*}
		x\in E_i\setminus2B,\quad d_i(x,w)<r/2
		&\quad\Longrightarrow\quad
		d(x,x_B)\leq d_i(x,w)+d(w,z)+d(z,x_B)<2r,\\
		&\quad\Longrightarrow\quad\text{a contradiction}.
	\end{align*}
	Thus $d_i(x,w)\geq r/2$ on the set of integration. Put
	$p_N=N-1/2-m/2>1/8$. Since $r/2>s$, \eqref{eq:hw-model-tail} gives
	\begin{align*}
		&\int_{E_i\setminus2B}\int_{\substack{w\in E_i\\d(w,z)\leq r/2}}
		|K_{F_s(\Delta_i)}(x,w)|\chi_i(w)|K_{\psi_0(\Delta)}(w,z)|
		\dd\mu(w)\dd\mu_i(x)\\
		&\quad\leq C_N\norm F_{C^N}(2s/r)^{p_N}
		\int_M|K_{\psi_0(\Delta)}(w,z)|\dd\mu(w)\\
		&\quad\leq C_N\norm F_{C^N}(s/r)^{1/8}.
	\end{align*}
	Summing the compact part and the two end parts proves
	\eqref{eq:mihlin-7}. This proves Proposition \ref{prop:mihlin-bands}.
	
	\subsection{The estimate on every atom}
	\label{subsec:mihlin-atom}
	
	\begin{proof}[Proof of Proposition \ref{prop:mihlin-atom}]
		Choose a smooth nonincreasing function $\kappa$ with
		$\kappa=1$ on $[0,2]$ and $\kappa=0$ on $[3,\infty)$. Put
		$$
		s_j=2^j,\qquad \psi(u)=\kappa(u)-\kappa(4u)\geq0,
		\qquad\supp\psi\subset[1/2,3]\Subset(1/4,4).
		$$
		The partition follows by telescoping:
		\begin{equation}\label{eq:hw-dyadic-partition}
			0\leq\sum_{j=-H}^H\psi(s_j^2\lambda)
			=\kappa(4^{-H}\lambda)-\kappa(4^{H+1}\lambda)
			\longrightarrow1\qquad(\lambda>0).
		\end{equation}
		Put $F_j(u)=\mathfrak m(u/s_j^2)\psi(u)$. For $0\leq k\leq N$,
		$$
		F_j^{(k)}(u)=\sum_{\ell=0}^k\binom k\ell
		s_j^{-2\ell}\mathfrak m^{(\ell)}(u/s_j^2)
		\psi^{(k-\ell)}(u),\qquad
		s_j^{-2\ell}|\mathfrak m^{(\ell)}(u/s_j^2)|
		\leq A_N(\mathfrak m)u^{-\ell}.
		$$
		Since $1/2\leq u\leq3$ on the supports, the Leibniz formula gives
		\begin{equation}\label{eq:hw-uniform-functions}
			\norm{F_j}_{C^N}+\norm{u^JF_j(u)}_{C^N}\leq C_{J,N}A_N(\mathfrak m).
		\end{equation}
		For the atom in \eqref{eq:mihlin-2}, Cauchy--Schwarz gives
		\begin{equation}\label{eq:hw-atom-L1}
			\norm a_1\leq\mu(B)^{1/2}\norm a_2\leq1,\qquad
			\norm b_1\leq\mu(B)^{1/2}\norm b_2\leq r^{2J}.
		\end{equation}
		On $2B$, the spectral theorem and \eqref{eq:hw-local-doubling} give
		\begin{equation}\label{eq:hw-near-atom}
			\norm{\mathfrak m(\Delta)a}_{L^1(2B)}
			\leq\mu(2B)^{1/2}A_0(\mathfrak m)\norm a_2
			\leq CA_0(\mathfrak m).
		\end{equation}
		For $s_j<r$, Proposition \ref{prop:mihlin-bands} gives
		$$
		\norm{F_j(s_j^2\Delta)a}_{L^1(M\setminus2B)}
		\leq\int_B|a(z)|\int_{M\setminus2B}
		|K_{F_j(s_j^2\Delta)}(x,z)|\dd\mu(x)\dd\mu(z)
		\leq C_NA_N(\mathfrak m)(s_j/r)^{1/8}.
		$$
		For $s_j\geq r$, use $a=\Delta^Jb$ first:
		\begin{align*}
			F_j(s_j^2\Delta)a
			&=s_j^{-2J}(u^JF_j(u))(s_j^2\Delta)b,\\
			\norm{F_j(s_j^2\Delta)a}_1
			&\leq C_{J,N}A_N(\mathfrak m)s_j^{-2J}\norm b_1
			\leq C_{J,N}A_N(\mathfrak m)(r/s_j)^{2J}.
		\end{align*}
		Taking the last $j$ with $2^j<r$ and the first with $2^j\geq r$
		gives
		\begin{align}
			\sum_{j\in\mathbb Z}\norm{F_j(s_j^2\Delta)a}_{L^1(M\setminus2B)}
			&\leq C_{J,N}A_N(\mathfrak m)
			\left(\sum_{s_j<r}(s_j/r)^{1/8}
			+\sum_{s_j\geq r}(r/s_j)^{2J}\right)\nonumber\\
			&\leq C_{J,N}A_N(\mathfrak m)
			\left(\frac1{1-2^{-1/8}}+\frac1{1-2^{-2J}}\right).
			\label{eq:hw-dyadic-sum}
		\end{align}
		There is no $L^2$ zero eigenfunction. Indeed,
		\begin{equation}\label{eq:hw-no-zero}
			\Delta f=0,\ f\in\operatorname{Dom} \Delta
			\quad\Longrightarrow\quad
			\norm{\Delta^{1/2}f}_2^2=\langle \Delta f,f\rangle=0
			\quad\Longrightarrow\quad\nabla f=0
			\quad\Longrightarrow\quad f=0,
		\end{equation}
		since $M$ is connected and $\mu(M)=\infty$.
		The spectral theorem and \eqref{eq:hw-dyadic-partition} imply
		$$
		f_H:=\sum_{j=-H}^HF_j(s_j^2\Delta)a\longrightarrow\mathfrak m(\Delta)a
		\quad\text{in }L^2.
		$$
		By \eqref{eq:hw-dyadic-sum}, $f_H$ also converges in
		$L^1(M\setminus2B)$ to a function $g$. For each set
		$\Omega\subset M\setminus2B$ of finite measure,
		$$
		\norm{g-\mathfrak m(\Delta)a}_{L^1(\Omega)}
		\leq\norm{g-f_H}_{L^1(M\setminus2B)}
		+\mu(\Omega)^{1/2}\norm{f_H-\mathfrak m(\Delta)a}_2
		\longrightarrow0.
		$$
		Combine this with \eqref{eq:hw-near-atom} to obtain
		\eqref{eq:mihlin-5}. At frequencies tending to zero,
		$s_j\to\infty$ and the summable factor is $(r/s_j)^{2J}$.
		Already $J=1$ suffices.
	\end{proof}
	
	\subsection{Extension, consistency, and interpolation}
	\label{subsec:mihlin-atomic-Hardy}
	Let $H^1_{\mathrm{at},\Delta^J,\delta}$ consist of $L^1$ functions
	represented by the atoms \eqref{eq:mihlin-2} with summable
	coefficients, with norm
	$$
	\norm f_{H^1_{\mathrm{at},\Delta^J,\delta}}
	=\inf\left\{\sum_\nu|c_\nu|:
	f=\sum_\nu c_\nu a_\nu\text{ in }L^1\right\}.
	$$
	This is a Banach space. If $\mathcal A_J$ is the family of atoms,
	then
	$$
	S:\ell^1(\mathcal A_J)\longrightarrow L^1,\qquad
	S(c)=\sum_{a\in\mathcal A_J}c_a a,\qquad
	\norm{S(c)}_1\leq\norm c_{\ell^1},\qquad
	H^1_{\mathrm{at},\Delta^J,\delta}\simeq\ell^1(\mathcal A_J)/\ker S.
	$$
	The kernel is closed; the quotient norm is the displayed atomic
	norm. Finite atomic sums are dense.
	
	The kernel and spectral definitions of every band operator are
	consistent on $L^1\cap L^2$. To see this, let $R_1,R_2$ denote
	these definitions and choose $f_k\in C_c^\infty(M)$ with
	$\norm{f-f_k}_1+\norm{f-f_k}_2\to0$. For $\mu(\Omega)<\infty$,
	\begin{equation}\label{eq:hw-band-consistency}
		\norm{R_1f-R_2f}_{L^1(\Omega)}
		\leq\norm{R_1}_{1\to1}\norm{f-f_k}_1
		+\mu(\Omega)^{1/2}\norm{R_2}_{2\to2}\norm{f-f_k}_2
		\longrightarrow0.
	\end{equation}
	Here $R_1f_k=R_2f_k$ by the kernel factorization. Thus these
	operators are well defined on $L^1+L^2$.
	
	\begin{lemma}\label{lem:band-separation}
		Set $Q_j=\psi(s_j^2\Delta)$. If $h\in L^1+L^2$ and $Q_jh=0$ for
		all $j$, then $h=0$.
	\end{lemma}
	\begin{proof}
		For $0<t\leq1$, \eqref{eq:hw-bounded-time-heat} and Minkowski's
		inequality give
		$$
		\norm{e^{-t\Delta}h_1}_2
		\leq\sup_z\norm{h_t(\cdot,z)}_2\norm{h_1}_1
		\leq Ct^{-m/4}\norm{h_1}_1.
		$$
		If $h=h_1+h_2$, $h_1\in L^1$, $h_2\in L^2$, then
		$e^{-t\Delta}h\in L^2$. To justify commutation on $L^1$, take
		$f_k\in C_c^\infty$ with $f_k\to h_1$ in $L^1$. We have
		\begin{align*}
			\norm{Q_je^{-t\Delta}(f_k-h_1)}_2
			&\leq Ct^{-m/4}\norm{Q_j}_{2\to2}\norm{f_k-h_1}_1,\\
			\norm{e^{-t\Delta}Q_j(f_k-h_1)}_2
			&\leq Ct^{-m/4}\norm{Q_j}_{1\to1}\norm{f_k-h_1}_1.
		\end{align*}
		Spectral commutation on $f_k$, followed by these limits, yields
		$$
		Q_je^{-t\Delta}h=e^{-t\Delta}Q_jh=0\quad(j\in\mathbb Z),\qquad
		e^{-t\Delta}h=\lim_{H\to\infty}\sum_{j=-H}^HQ_je^{-t\Delta}h=0
		\quad\text{in }L^2.
		$$
		The last limit uses \eqref{eq:hw-dyadic-partition} and
		\eqref{eq:hw-no-zero}.
		
		For $v\in C_c^\infty(M)$, spectral calculus gives
		$$
		v-e^{-t\Delta}v=\int_0^t e^{-u\Delta}\Delta v\dd u,\qquad
		\norm{v-e^{-t\Delta}v}_1\leq t\norm{\Delta v}_1.
		$$
		Contraction and density give, for $h_1\in L^1$,
		$$
		\norm{h_1-e^{-t\Delta}h_1}_1
		\leq2\norm{h_1-v}_1+t\norm{\Delta v}_1.
		$$
		$$
		\limsup_{t\downarrow0}\norm{h_1-e^{-t\Delta}h_1}_1
		\leq2\norm{h_1-v}_1.
		$$
		Taking the infimum over $v\in C_c^\infty(M)$ gives zero.
		Strong convergence on $L^2$ follows from the spectral theorem.
		Hence, with the usual sum norm,
		$$
		\norm{h-e^{-t\Delta}h}_{L^1+L^2}
		\leq\norm{h_1-e^{-t\Delta}h_1}_1+
		\norm{h_2-e^{-t\Delta}h_2}_2\longrightarrow0.
		$$
		Since $e^{-t\Delta}h=0$, we obtain $h=0$.
	\end{proof}
	
	\begin{proof}[Proof of Theorem \ref{thm:spectral-multiplier}]
		For an atomic representation $f=\sum_\nu c_\nu a_\nu$, define
		$$
		Tf=\sum_\nu c_\nu\mathfrak m(\Delta)a_\nu.
		$$
		The series converges absolutely in $L^1$ by
		\eqref{eq:mihlin-5}. The band operator
		$$
		T_j=(\mathfrak m(\lambda)\psi(s_j^2\lambda))(\Delta)
		$$
		is bounded on $L^1$ by Proposition \ref{prop:mihlin-bands} and agrees
		with its spectral definition by \eqref{eq:hw-band-consistency}.
		For each atom, spectral calculus gives
		$Q_j\mathfrak m(\Delta)a_\nu=T_ja_\nu$. Passing through the
		$L^1$ convergent series,
		$$
		Q_jTf=\sum_\nu c_\nu T_ja_\nu
		=T_j\sum_\nu c_\nu a_\nu=T_jf.
		$$
		If two representations give $g,g'$ as their images, then
		$Q_j(g-g')=0$ for all $j$. Lemma \ref{lem:band-separation} gives
		$g=g'$. Taking the infimum over representations yields
		$$
		\norm{Tf}_1\leq C_{J,N,\delta}A_N(\mathfrak m)
		\norm f_{H^1_{\mathrm{at},\Delta^J,\delta}}.
		$$
		Density of finite atomic sums gives uniqueness. If also
		$f\in L^2$, then
		$$
		Q_j\bigl(Tf-\mathfrak m(\Delta)f\bigr)=T_jf-T_jf=0
		\quad\text{in }L^1+L^2.
		$$
		The same lemma proves $Tf=\mathfrak m(\Delta)f$.
		
		For $J=1$, Theorem \ref{thm:hardy-characterization}(i) gives
		$H^1_{\mathrm{at},\Delta,\delta}= H^1_{S_c}(M) $ with equivalent norms.
		Part (iii) gives \eqref{eq:interpolation}. With $A_1=C_{N,c,\delta}A_N(\mathfrak m)$
		and $A_2=A_0(\mathfrak m)$, the two consistent bounds are
		$$
		\norm{Tf}_1\leq A_1\norm f_{ H^1_{S_c}(M) },\qquad
		\norm{Tf}_2\leq A_2\norm f_2.
		$$
		Consistency defines $T$ on $ H^1_{S_c}(M) +L^2$ by
		$T(f_1+f_2)=Tf_1+\mathfrak m(\Delta)f_2$.
		If $A_2=0$, then $T=0$. Otherwise, for
		$\vartheta=2-2/p$, $1<p<2$, the $K$-functional satisfies
		$$
		K(t,Tf;L^1,L^2)
		\leq A_1K(tA_2/A_1,f; H^1_{S_c}(M) ,L^2).
		$$
		Consequently,
		\begin{align*}
			\norm{Tf}_p
			&\leq C_p\left(\int_0^\infty
			[t^{-\vartheta}K(t,Tf;L^1,L^2)]^p\frac{\dd t}{t}\right)^{1/p}\\
			&\leq C_p A_1^{1-\vartheta}A_2^{\vartheta}
			\norm f_{( H^1_{S_c}(M) ,L^2)_{\vartheta,p}}\\
			&\leq C_{p,N,c,\delta}A_N(\mathfrak m)^{2/p-1}
			A_0(\mathfrak m)^{2-2/p}\norm f_p.
		\end{align*}
		For $p>2$, put $p'=p/(p-1)<2$ and apply the preceding bound
		to $\overline{\mathfrak m}$. For $f\in L^p\cap L^2$,
		\begin{align*}
			\norm{\mathfrak m(\Delta)f}_p
			&=\sup_{\substack{g\in L^{p'}\cap L^2\\\norm g_{p'}\leq1}}
			|\langle\mathfrak m(\Delta)f,g\rangle|
			=\sup_{\substack{g\in L^{p'}\cap L^2\\\norm g_{p'}\leq1}}
			|\langle f,\overline{\mathfrak m}(\Delta)g\rangle|
			\leq C_{p,N,c,\delta}A_N(\mathfrak m)\norm f_p.
		\end{align*}
		Density extends the operator to $L^p$. The case $p=2$ is the
		spectral theorem. This proves Theorem \ref{thm:spectral-multiplier}.
	\end{proof}
	
	\subsection{The assumptions and the derivative count}
	The facts used in the proof are the real heat kernel bounds,
	the mass bound, finite propagation, the exact product metrics,
	and the compactly supported energy identity
	\eqref{eq:hw-elliptic}. The last identity is proved by integration
	by parts. It requires no gradient estimate for the heat kernel.
	
	The restriction $n>2$ enters as follows:
	$$
	s^{2-n}\leq1\quad(s\geq1,n\geq2),\qquad
	\sum_{\nu_h\leq4}\nu_h^{n/2-1}\leq C_n\quad(n>2).
	$$
	At $n=2$, the second sum has a number of terms comparable to
	$1+\log s$. The present proof gives no uniform estimate in that
	case.
	
	The Fourier moment required in \eqref{eq:hw-Duhamel-L1} is
	$$
	\int_\R|\widehat F(u)|\,|u|(1+|u|)^{m/2+1/6}\dd u,
	\qquad
	N>(m+3)/2\quad\Longrightarrow\quad N-m/2\geq2.
	$$
	Thus one may take $N=\lfloor(m+3)/2\rfloor+1$.

	\subsection{Imaginary powers and cancellation}
	\label{subsec:imaginary-cancellation}
	
	\begin{proposition}\label{prop:imaginary-cancellation}
		For every real $\tau\ne0$, there are nonnegative functions
		$a_R\in C_c^\infty(M)$, supported in unit balls $B(x_R,1)$
		in the exact Euclidean end, such that
		\begin{equation}\label{eq:imaginary-cancellation-failure}
			\|a_R\|_2\leq\mu(B(x_R,1))^{-1/2},\qquad
			\|a_R\|_1\leq1,\qquad
			\|\Delta^{i\tau}a_R\|_1\longrightarrow\infty
			\quad(R\to\infty).
		\end{equation}
		The last norm may be infinite.
	\end{proposition}
	
	\begin{proof}
		Write $\Delta_{\mathrm E}=-\sum_{\ell=1}^m\partial_\ell^2$ on
		$\R^m$. Fix $\psi\in C_c^\infty(B_{\R^m}(0,1))$ such that
		$$
		\psi\geq0,\qquad \int_{\R^m}\psi>0,\qquad
		\|\psi\|_2\leq |B_{\R^m}(0,1)|^{-1/2}.
		$$
		For $R>2$, choose $x_R$ in the exact Euclidean end so that
		$B_{\R^m}(x_R,R)$ lies in that end. Define
		$$
		a_R(x_R+X)=\psi(X)\quad(|X|<1),\qquad
		a_R=0\quad\text{elsewhere on }M.
		$$
		The functions are smooth and nonnegative. Since the metric on
		their supports is Euclidean,
		$$
		\|a_R\|_2=\|\psi\|_2\leq\mu(B(x_R,1))^{-1/2},\qquad
		\|a_R\|_1\leq\mu(B(x_R,1))^{1/2}\|a_R\|_2\leq1.
		$$
		We will prove, for every fixed $D>1$,
		\begin{equation}\label{eq:imaginary-translation}
			\|\Delta^{i\tau}a_R(x_R+\cdot)-\Delta_{\mathrm E}^{i\tau}\psi\|
			_{L^2(B_{\R^m}(0,D))}
			\longrightarrow0.
		\end{equation}
		
		Let $q\in C_c^\infty(\R)$ be even and put
		$\widehat q(s)=\int_\R e^{-isv}q(v)\dd v$. The spectral theorem
		gives the absolutely convergent $L^2$ integral
		$$
		q(\sqrt \Delta)a_R
		=\frac1{2\pi}\int_\R\widehat q(s)
		\cos(s\sqrt \Delta)a_R\dd s.
		$$
		For $|s|<R-1$, finite propagation confines both wave solutions
		to $B_{\R^m}(x_R,|s|+1)$. The metric and the Laplacian there
		are Euclidean. The translated Euclidean solution, extended
		by zero outside that ball, solves the same wave equation on
		$M$ with the same initial values. Uniqueness of the $L^2$
		wave solution therefore gives
		$$
		\cos(s\sqrt \Delta)a_R(x_R+X)
		=\cos(s\sqrt{\Delta_{\mathrm E}})\psi(X)
		\qquad(|s|<R-1,\ |X|<R).
		$$
		For $R>D+1$, the $L^2$ contraction of each cosine operator
		now gives
		\begin{equation}\label{eq:imaginary-smooth-translation}
			\begin{split}
				&\|q(\sqrt \Delta)a_R(x_R+\cdot)
				-q(\sqrt{\Delta_{\mathrm E}})\psi\|_{L^2(B(0,D))}\\
				&\hspace{12mm}\leq
				\frac{\|\psi\|_2}{\pi}
				\int_{|s|\geq R-1}|\widehat q(s)|\dd s
				\longrightarrow0.
			\end{split}
		\end{equation}
		
		\medskip
		\noindent\emph{The spectral values near zero and infinity.}
		Choose $\chi\in C_c^\infty(\R)$ with $0\leq\chi\leq1$,
		$\chi=1$ on $[-1,1]$, and $\chi=0$ outside $[-2,2]$.
		For $0<\varepsilon<1$, spectral positivity gives
		$$
		\|\mathbf1_{[0,\varepsilon]}(\Delta)a_R\|_2^2
		\leq\langle\chi(\Delta/\varepsilon)a_R,a_R\rangle.
		$$
		Apply \eqref{eq:imaginary-smooth-translation} with
		$q(v)=\chi(v^2/\varepsilon)$ and $D=2$. Since $a_R$ is
		supported in its unit ball,
		\begin{align}
			\limsup_{R\to\infty}
			\|\mathbf1_{[0,\varepsilon]}(\Delta)a_R\|_2^2
			&\leq\langle\chi(\Delta_{\mathrm E}/\varepsilon)\psi,\psi\rangle
			\notag\\
			&=\frac1{(2\pi)^m}\int_{\R^m}
			\chi(|\xi|^2/\varepsilon)|\widehat\psi(\xi)|^2\dd\xi
			\notag\\
			&\leq C_m\varepsilon^{m/2}\|\psi\|_1^2.
			\label{eq:imaginary-low-spectral}
		\end{align}
		Here $|\widehat\psi|\leq\|\psi\|_1$. The same Fourier estimate
		bounds $\|\mathbf1_{[0,\varepsilon]}(\Delta_{\mathrm E})\psi\|_2^2$.
		For an integer $J\geq1$ and $K>1$, the spectral theorem gives
		\begin{equation}\label{eq:imaginary-high-spectral}
			\|\mathbf1_{[K,\infty)}(\Delta)a_R\|_2
			\leq K^{-J}\|\Delta^Ja_R\|_2
			=K^{-J}\|\Delta_{\mathrm E}^J\psi\|_2.
		\end{equation}
		The equality holds because $a_R$ is smooth and its support
		lies in the exact Euclidean end. This bound also holds with
		$\Delta_{\mathrm E}$ and $\psi$.
		
		Fix $0<\varepsilon<1<K$. Choose
		$\eta\in C_c^\infty((\varepsilon/2,2K))$ with
		$0\leq\eta\leq1$ and $\eta=1$ on $[\varepsilon,K]$.
		Put $g(\lambda)=\lambda^{i\tau}\eta(\lambda)$ for
		$\lambda>0$ and $g(0)=0$. Then $q(v)=g(v^2)$ is even and
		belongs to $C_c^\infty(\R)$. Moreover,
		$$
		\|(\Delta^{i\tau}-g(\Delta))a_R\|_2
		\leq\|\mathbf1_{[0,\varepsilon]}(\Delta)a_R\|_2
		+\|\mathbf1_{[K,\infty)}(\Delta)a_R\|_2,
		$$
		because $|\lambda^{i\tau}|=1$ for $\lambda>0$.
		The same estimate holds on $\R^m$. Combining
		\eqref{eq:imaginary-smooth-translation}--\eqref{eq:imaginary-high-spectral}
		therefore gives
		\begin{align*}
			\limsup_{R\to\infty}
			&\|\Delta^{i\tau}a_R(x_R+\cdot)-\Delta_{\mathrm E}^{i\tau}\psi\|
			_{L^2(B(0,D))}
			\leq C_m\varepsilon^{m/4}\|\psi\|_1
			+2K^{-J}\|\Delta_{\mathrm E}^J\psi\|_2.
		\end{align*}
		Let $\varepsilon\downarrow0$ and $K\to\infty$. This proves
		\eqref{eq:imaginary-translation}. The values assigned
		to the multipliers at zero do not matter: both Laplacians
		have zero $L^2$ kernel.
		
		\medskip
		\noindent\emph{The Euclidean tail.}
		We derive the kernel needed here. For $\lambda>0$,
		$$
		\lambda^{i\tau}
		=\frac1{\Gamma(1-i\tau)}
		\int_0^\infty s^{-i\tau}\lambda e^{-s\lambda}\dd s.
		$$
		For our smooth compactly supported $\psi$, the Euclidean
		heat kernel and Young's inequality give
		$$
		\|\Delta_{\mathrm E}e^{-s\Delta_{\mathrm E}}\psi\|_2
		\leq
		\begin{cases}
			\|\Delta_{\mathrm E}\psi\|_2,&0<s\leq1,\\
			C_ms^{-1-m/4}\|\psi\|_1,&s\geq1.
		\end{cases}
		$$
		Thus the following integral converges absolutely in $L^2$:
		$$
		\Delta_{\mathrm E}^{i\tau}\psi
		=\frac1{\Gamma(1-i\tau)}\int_0^\infty
		s^{-i\tau}\Delta_{\mathrm E}e^{-s\Delta_{\mathrm E}}\psi\dd s.
		$$
		To justify the identity, first integrate over a finite
		interval. Its scalar multiplier has modulus at most
		$|\Gamma(1-i\tau)|^{-1}$, since
		$\int_0^\infty\lambda e^{-s\lambda}\dd s=1$.
		Dominated convergence in the spectral theorem then gives
		the displayed $L^2$ identity.
		
		For $|X|>1$, the pointwise integral also converges locally
		uniformly. At zero the heat kernel has exponential decay
		away from $\supp\psi$, and at infinity it is at most
		$C_ms^{-m/2}\|\psi\|_1$. Since
		$\Delta_{\mathrm E}e^{-s\Delta_{\mathrm E}}=-\partial_s e^{-s\Delta_{\mathrm E}}$,
		integration by parts has zero boundary terms and gives
		\begin{align*}
			\Delta_{\mathrm E}^{i\tau}\psi(X)
			&=\frac1{\Gamma(-i\tau)}\int_0^\infty
			s^{-i\tau-1}e^{-s\Delta_{\mathrm E}}\psi(X)\dd s\\
			&=\frac1{\Gamma(-i\tau)}\int_{\R^m}\psi(Z)
			\int_0^\infty
			\frac{s^{-i\tau-1}}{(4\pi s)^{m/2}}
			e^{-|X-Z|^2/(4s)}\dd s\dd Z\\
			&=c_{m,\tau}\int_{\R^m}|X-Z|^{-m-2i\tau}\psi(Z)\dd Z,
		\end{align*}
		where absolute convergence justifies Fubini, and the change
		of variable $v=|X-Z|^2/(4s)$ gives
		$$
		c_{m,\tau}
		=\frac{4^{i\tau}\Gamma(m/2+i\tau)}
		{\pi^{m/2}\Gamma(-i\tau)}\ne0.
		$$
		The last inequality follows because the Gamma function has
		no zeros and neither argument is a pole when $\tau\ne0$.
		For $|X|\geq2$ and $|Z|\leq1$, differentiation of
		$|X-Z|^{-m-2i\tau}$ gives
		$$
		\bigl||X-Z|^{-m-2i\tau}-|X|^{-m-2i\tau}\bigr|
		\leq C_{m,\tau}|Z||X|^{-m-1}.
		$$
		Consequently,
		\begin{equation}\label{eq:Euclidean-imaginary-tail}
			\Delta_{\mathrm E}^{i\tau}\psi(X)
			=c_{m,\tau}\Bigl(\int\psi\Bigr)|X|^{-m-2i\tau}
			+O_{m,\tau,\psi}(|X|^{-m-1}).
		\end{equation}
		There are constants $D_0>2$ and $c_0>0$ such that
		$$
		|\Delta_{\mathrm E}^{i\tau}\psi(X)|\geq c_0|X|^{-m}
		\qquad(|X|\geq D_0).
		$$
		Polar coordinates therefore give
		$$
		\int_{B(0,D)}|\Delta_{\mathrm E}^{i\tau}\psi(X)|\dd X
		\geq c_1\log(D/D_0)\longrightarrow\infty.
		$$
		For every fixed $D>D_0$, \eqref{eq:imaginary-translation}
		and Cauchy--Schwarz give
		$$
		\liminf_{R\to\infty}\|\Delta^{i\tau}a_R\|_1
		\geq\lim_{R\to\infty}\int_{B(0,D)}
		|\Delta^{i\tau}a_R(x_R+X)|\dd X
		=\int_{B(0,D)}|\Delta_{\mathrm E}^{i\tau}\psi(X)|\dd X.
		$$
		Let $D\to\infty$. This proves \eqref{eq:imaginary-cancellation-failure}.
	\end{proof}
	
	The multiplier in this proposition satisfies, for every
	integer $N\geq0$,
	$$
	A_N(\mathfrak m_\tau)
	=\max_{0\leq j\leq N}
	\left|\prod_{\ell=0}^{j-1}(i\tau-\ell)\right|<\infty,
	\qquad \mathfrak m_\tau(\lambda)=\lambda^{i\tau}\quad(\lambda>0),
	$$
	where the empty product equals $1$. It has no limit at zero.
	Every $\Delta$-atom in \eqref{eq:atoms} has zero integral, whereas
	$\int_Ma_R\dd\mu=\int_{\R^m}\psi>0$.

	\section{Further discussions}

\begin{remark}
The restriction of the horizontal radius gives the norm estimates
over all heights proved above. For $0<t<c\rho(x)$, its condition
is still $d(x,y)<t$, so the Poisson extension has the usual
non-tangential boundary limit by Corollary~\ref{cor:restricted-boundary}.
\end{remark}

	On manifolds with ends, which may be nondoubling, Hassell and Sikora \cite{HS} studied the $L^p$ boundedness of the Riesz
	transform. Furthermore, Hassell, Nix and Sikora \cite{HNS} proved that the Riesz
	transform is of weak type $(1,1)$.

	Very recently, Ouyang, Spector and Stockdale \cite{OSS} obtained the dimension-free weak-type $(1, 1)$
	bound for the vector Riesz transform on ${\mathbb R}^n$, which addressed the question raised by E.M. Stein.
	Then, Chen, Jiang, B. Li and H.-Q. Li used the method in \cite{OSS} to prove that the Riesz
	transform is of weak type $(1,1)$ on every complete noncompact
	Riemannian manifold \cite{CJLL2026}. Consequently,
	$$
	\|\nabla \Delta^{-1/2}f\|_p\leq C_p\|f\|_p,
	\qquad 1<p\leq2,
	$$
	where $\Delta$ is the nonnegative Laplace--Beltrami operator and $C_p$
	depends only on $p$. Their proof requires no curvature, doubling,
	or heat kernel assumption.
	
	\begin{remark}\label{rem:Riesz-Hardy-question}
		We leave open the question of whether the Riesz transform
		characterizes the Hardy space $ H^1_{S_c}(M) $ in \eqref{eq:H1}.
		More precisely, for $f\in L^2(M)$, do we have
		$$
		f\in H^1_{S_c}(M)
		\quad\Longleftrightarrow\quad
		f\in L^1(M)\ \text{and}\ |\nabla \Delta^{-1/2}f|\in L^1(M),
		$$
		with
		$$
		\|f\|_{ H^1_{S_c}(M) }
		\approx\|f\|_1+\|\nabla \Delta^{-1/2}f\|_1?
		$$
		Here the Riesz transform is understood as its bounded extension
		on $L^2$, and $|\nabla \Delta^{-1/2}f|$ is its pointwise length.
	\end{remark}
	
	\vskip 0.5cm
	
	\noindent
{\bf Acknowledgments}:  During preparation of the manuscript, the authors used GPT-6 for assistance with
	concrete computations and language editing. The authors take responsibility for the mathematical
	arguments and the final content of the paper. 	 P. Chen was supported by National Key R$\&$D Program of China 2022YFA1005700
	and by NNSF of China No. 12671127.  J. Li is supported by the Australian Research
	Council (DP260100485).  L. Yan was supported by National Key R$\&$D Program of China
	2022YFA1005700 and by NNSF of China No. 12571111.
	
\bigskip
	
\noindent 
Data availability No data was used for the preparation of this manuscript.

\bigskip
\noindent
{\bf Declarations}

\noindent
Conflict of interest On behalf of all authors, the corresponding author states that there is no conflict of
interest.

\end{document}